\documentclass[]{amsart}
\usepackage{color,latexsym,amssymb,amsmath,amsthm}
\usepackage{listings}
\usepackage{comment}
\usepackage{graphicx}

\definecolor{brightred}{rgb}{1.0,0.9,0.9}
\definecolor{brightblue}{rgb}{0.9,0.9,1.0}
\definecolor{brightyellow}{rgb}{1.0,1.0,0.8}
\definecolor{brightgreen}{rgb}{0.9,1.0,0.9}
\def\absatz#1{\begin{center}\fcolorbox{brightgreen}{brightgreen}{\parbox{9cm}{#1}}\end{center}}
\def\appendix#1{\fcolorbox{brightblue}{brightblue}{#1}}
\theoremstyle{definition}

\title{Illustrating Mathematics using 3D Printers}
\author{Oliver Knill and Elizabeth Slavkovsky}
\date{June 20, 2013}
\address{
        Department of Mathematics, Harvard University, Cambridge, MA, 02138. 
        knill@math.harvard.edu, writetoliz@yahoo.com
        }


\begin{document}
\maketitle

\begin{abstract}
3D printing technology can help to visualize proofs in mathematics. 
In this document we aim to illustrate how 3D printing can help to visualize 
concepts and mathematical proofs. As already known to educators
in ancient Greece, models allow to bring mathematics closer to the public.
The new 3D printing technology makes the realization of such tools more
accessible than ever. This is an updated version of a paper included in \cite{CFZ}. 
\end{abstract}

\section{Visualization}
Visualization has always been an important ingredient for communicating mathematics \cite{Visualization}. 
Figures and models have helped to express ideas even before formal mathematical language was able 
to describe the structures. Numbers have been recorded as marks on bones, represented with pebbles, then
painted onto stone, inscribed into clay, woven into talking knots, 
written onto papyrus or paper, then printed on paper or displayed on computer screens. 
While figures extend language and pictures allow to visualize concepts, realizing objects in 
space has kept its value. Already in ancient Greece, wooden models of Apollonian cones were used
to teach conic sections. Early research in mathematics was often visual: figures on Babylonian Clay 
tablets illustrate Pythagorean triples, the Moscow mathematical papyrus features a picture which 
helps to derive the volume formula for a frustum. Al-Khwarizmi drew figures to solve the
quadratic equation. Visualization is not only illustrative, educational or heuristic,
it has practical value: Pythagorean triangles realized by ropes helped measuring and dividing up of 
land in Babylonia. Ruler and compass, introduced to construct mathematics on paper, 
can be used to build plans for machines. Greek mathematicians like Apollonius, Aristarchus,
Euclid or Archimedes mastered the art of representing mathematics with figures \cite{NetzDeduction}. 
Visual imagination plays an important role in extending geometrical knowledge \cite{Giaquinto}.
While pictures do not replace proofs - \cite{kline} gives a convincing visual proof that 
all triangles are equilateral - they help to transmit intuition about results and ideas \cite{Tufte,Bender}.
The visual impact on culture is well documented \cite{Emmer}. 
Visualization is especially crucial for education and can lead to new insight. Many examples of mechanical 
nature are in the textbook \cite{Levi}. As a pedagogical tool, it assists teachers on any level 
of mathematics, from elementary and high school over higher education to modern research
\cite{HannaSidoli,Posamentier,Palais}.
A thesis of \cite{Slavkovsky} has explored the feasibility of the technology 
in the classroom. We looked at work of Archimedes \cite{knillslavkovsky} using this technology.
Visualizations helps also to showcase the beauty of mathematics and to promote the field to a
larger public. Figures can inspire new ideas, generate new theorems or assist in computations;
examples are  Feynman or Dynkin diagrams or Young tableaux. Most mathematicians draw creative ideas and 
intuition from pictures, even so these figures often do not make it into papers or textbooks.
Artists \cite{Ferguson}, architects \cite{Hahn}, film makers, engineers and designers draw inspiration from visual mathematics. 
Well illustrated books like \cite{bookofnumbers,symmetries,Pickover2009,Jackson2012,Fomenko,symmetries,BergerPanorama,ApostolHorizons} 
advertise mathematics with figures and illustrations. Such publications help to counterbalance the impression 
that mathematics is difficult to communicate to non-mathematicians. Mathematical exhibits like at the 
science museum in Boston or the Museum of Math in New York play an important role in making
mathematics accessible. They all feature visual or even hands-on realizations of mathematics.
While various technologies have emerged which allow to display spacial and dynamic content
on the web, like Javascript, Java, Flash, WRML, SVG or WebGl, the possibility to 
{\bf manipulate an object with our bare hands} is still unmatched. 3D printers allow us to do that 
with relatively little effort.

\section{3D printing}

The industry of rapid prototyping and 3D printing in particular emerged about 30 years ago
\cite{Cooper,CLL,KamraniNasr,Brain,GRS} and is by some considered part of an 
{\bf industrial revolution} in which manufacturing becomes digital, personal, and affordable 
\cite{Rifkin,worldfinancialreview,economist2012}. 
First commercialized in 1994 with printed wax material,
the technology has moved to other materials like acrylate photopolymers or metals and is now entering
the range of consumer technology. Printing services can print in color, with various materials and
in high quality. The development of 3D printing is the latest piece in a chain of 
visualization techniques. We live in an exciting time, because we experience not only 
one revolution, but two revolutions at the same time: an information revolution and an industrial revolution.
These changes also affect mathematics education \cite{BorweinRochaRodgrigues}.
3D printing is now used in the medical field, the airplane industry,
to prototype robots, to create art and jewelry, to build nano structures, bicycles, ships, circuits,
to produce art, robots, weapons, houses and even used to decorate cakes. 
Its use in education was investigated in \cite{Slavkovsky}. Since physical models are important for hands-on active learning  \cite{Toolsteaching},
3D printing technology in education has been used since a while
\cite{Lipson} and considered for \cite{sustainabledevelopement}
sustainable development, for K-12 education in STEM projects \cite{rapman} as well as
elementary mathematics education \cite{contemporary}. 
There is no doubt that it will have a huge impact 
in education \cite{Cliff,BullGroves}. Printed models allow to illustrate 
concepts in various mathematical fields like calculus, geometry or topology.
It already has led to new prospects in mathematics education.
The literature about 3D printing explodes, similar as in computer literature expanded,
when PCs entered the consumer market. Examples of books are
\cite{Evans,Singh,Kelly}. As for any emerging technology, these publications might
be outdated quickly, but will remain a valuable testimony of the exciting time we live in.

\section{Bringing mathematics to life}
The way we think about mathematics influences our teaching \cite{Gold}. Images and objects can 
influence the way we think about mathematics. 
To illustrate visualizations using 3D printers, our focus is on mathematical
models generated with the help of {\bf computer algebra systems}. Unlike
{\bf 3D modelers}, mathematical software has the advantage that the source code is short
and that programs used to illustrate mathematics for research or the classroom can be reused. Many of the 
examples given here have been developed for classes or projects and redrawn so that it can be printed.
In contrast to ``modelers", software  which generate a large list of triangles, 
computer algebra systems describe and display three dimensional objects mathematically. 
While we experimented also with other software like ``123D Design" from Autodesk, ``Sketchup" from Trimble,
the modeler ``Free CAD", ``Blender",  or ``Rhinoceros" by McNeel Accociates, we worked mostly with computer 
algebra systems and in particular with Mathematica \cite{BorweinSkerritt,Trott,Wagon,Maeder,WellinKamin}. 
To explain this with a concrete example, lets look at a {\bf theorem of Newton on sphere packing} 
which tells that {\bf the kissing number of spheres in three dimensional space is 12.}
The theorem tells that the maximal number of spheres we can be placed around a given sphere
is twelve, if all spheres have the same radius, touch the central sphere and do not overlap.
While Newton's contemporary Gregory thought that one can place a thirteenth sphere, Newton believed
the kissing number to be 12. The theorem was only proved in 1953 \cite{CS}.
To show that the kissing number is at least 12, take an icosahedron with side length $2$ and place 
unit spheres at each of the 12 vertices then they kiss the unit sphere centered at the origin.
The proof that it is impossible to place 13 spheres \cite{Leech} uses an elementary estimate
\cite{OosteromStrackee} for the area of a spherical triangle, 
the Euler polyhedron formula, the discrete Gauss-Bonnet theorem 
assuring that the sum of the curvatures is $2$ and some combinatorics to check
 through all cases of polyhedra, which are allowed by these constraints.  
In order to visualize the use of Mathematica, we plotted 12 spheres which kiss a central sphere. 
While the object consists of 13 spheres only, the entire solid is made of 8640 triangles. 
The Mathematica code is very short because we only need to compute the
vertex coordinates of the icosahedron, generate the object and then export the STL file.
By displaying the source code, we have {\bf illustrated the visualization}, similar than 
communicating proof. If fed into the computer, the code generates a printable "STL" or "X3D" 
file.

\section{Sustainability considerations}
Physical models are important for hands-on active learning. 
Repositories of 3D printable models for education have emerged \cite{Lipson}. 
3D printing technology has been used for K-12 education in STEM projects \cite{rapman},
and elementary mathematics education \cite{contemporary}. There is optimism that it
will have a large impact in education \cite{Cliff}. 
The new technology allows everybody to build models for the classroom - in principle. 
To make it more accessible, many hurdles still have to be taken.
There are some good news: the STL files can be generated easily because the format is simple 
and open. STL files can also be exported to other formats. Mathematica for example allows to import it
and convert it to other forms. Programs like ``Meshlab" allow to manipulate it. Terminal 
conversions like ``admesh" allow to deal with STL files from the command line.
Other stand-alone programs like ``stl2pov" allow to convert it into
a form which can be rendered in a ray tracer like Povray.
One major point is that good software to generate the objects is not cheap.
The use of a commercial computer algebra system like Mathematica can be 
costly, especially if a site licence is required. There is no
free computer algebra software available now, which is able to export STL or 3DS 
or WRL files with built in routines. The computer algebra system SAGE, which is the 
most sophisticated open source system, has only export in experimental stage \cite{Olah}. 
It seems that a lot of work needs to be done there. Many resources are available
however \cite{Hart,Thingiverse}.
The following illustrations consist of Mathematica graphics which could be
printed. This often needs adaptation because a printer can not print objects
of zero thickness. 
\vfill
\pagebreak

\section{Illustrations}
\vspace{-5mm}
\parbox{16.8cm}{ \parbox{10cm}{
\scalebox{0.121}{\includegraphics{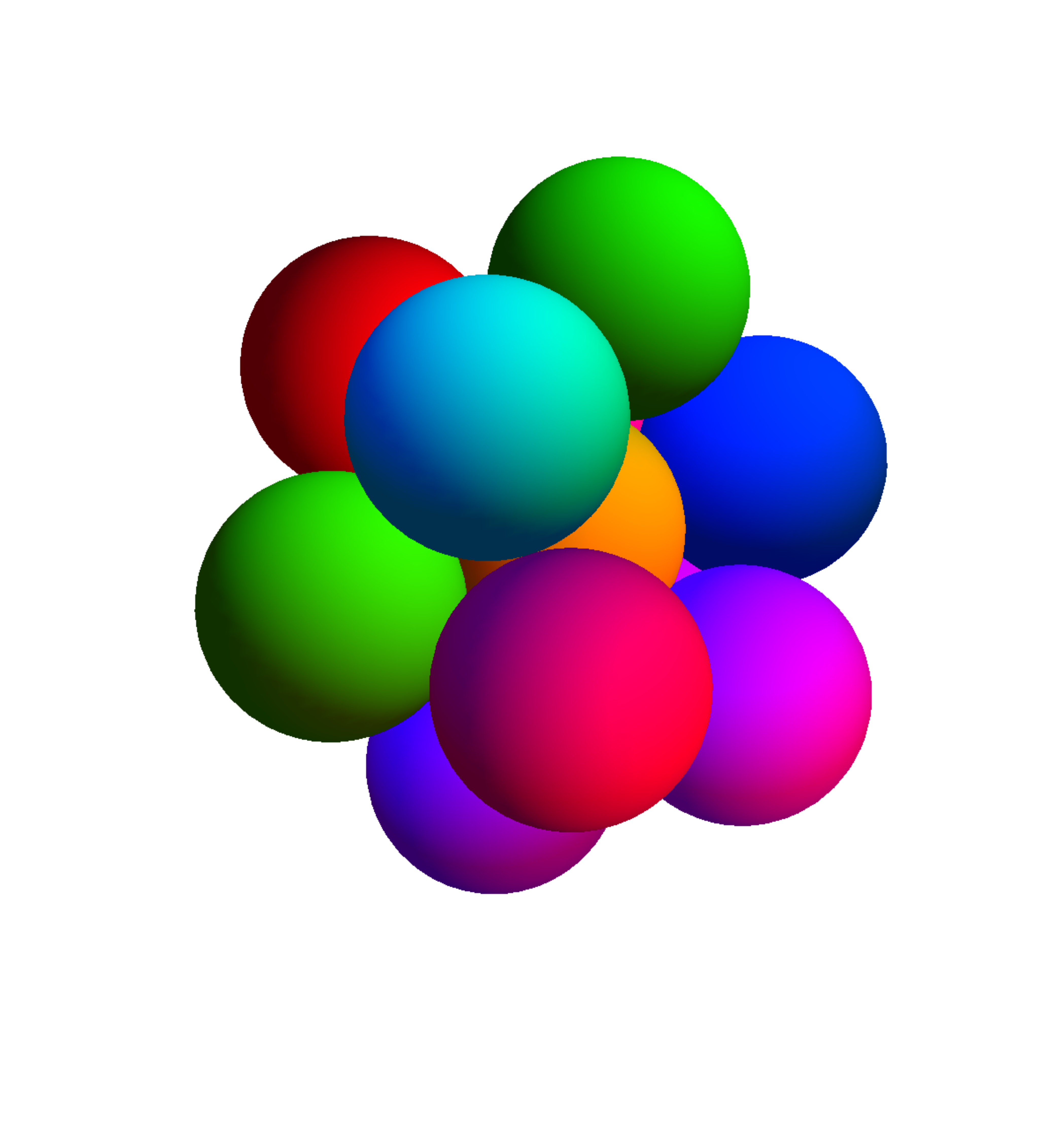}}
\scalebox{0.121}{\includegraphics{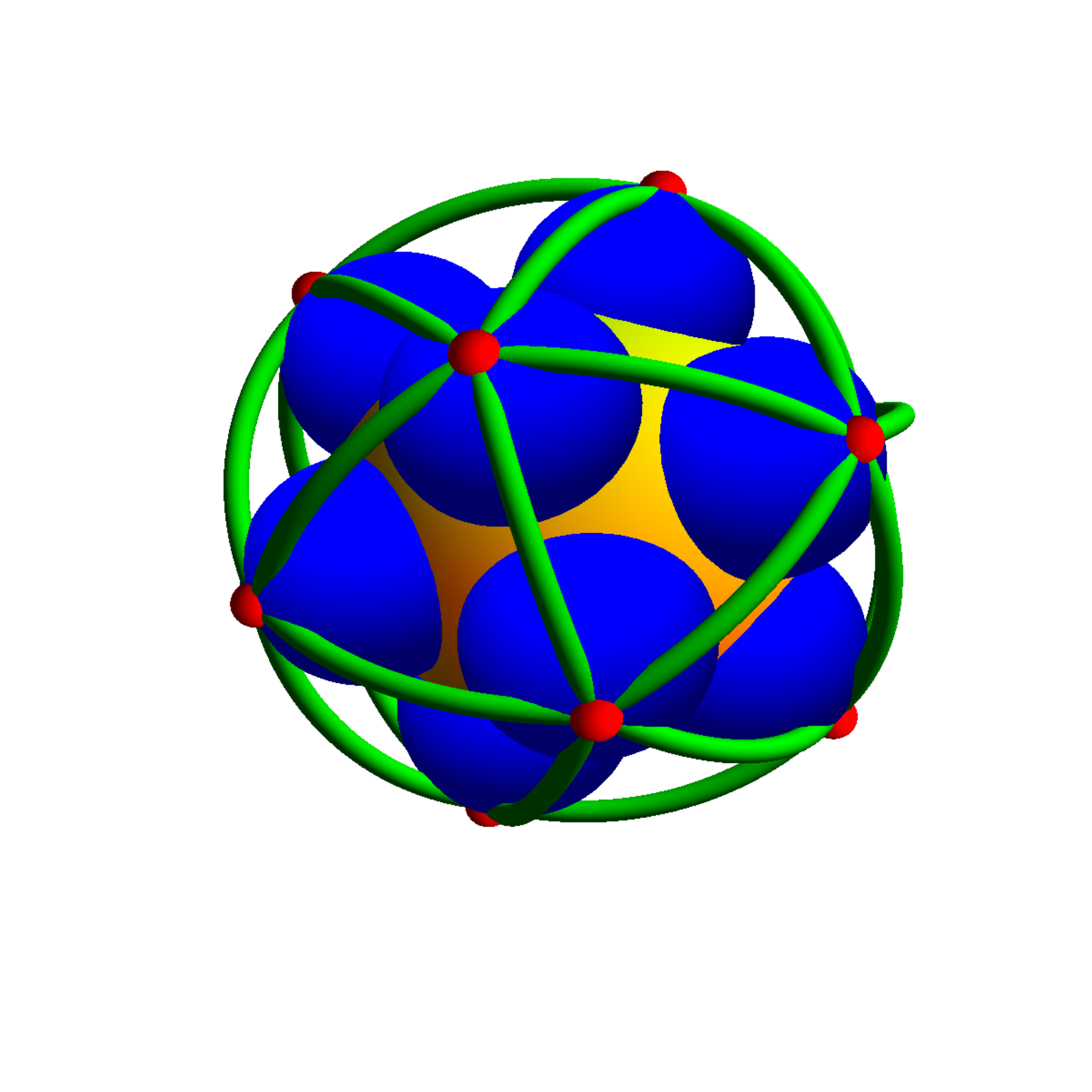}}
} \parbox{6cm}{
This figure aims to visualize that the kissing number
of a sphere is $\geq 12$. The Mathematica code producing this
object is given in the text. It produces a file containing tens of thousands
of triangles which the 3D printer knows to bring to live. The printed
object visualizes that there is still quite a bit of space left on the sphere.
Newton and his contemporary Gregory had a disagreement over whether this is 
enough to place a thirteenth sphere. 
}}

\parbox{16.8cm}{ \parbox{10cm}{
\scalebox{0.121}{\includegraphics{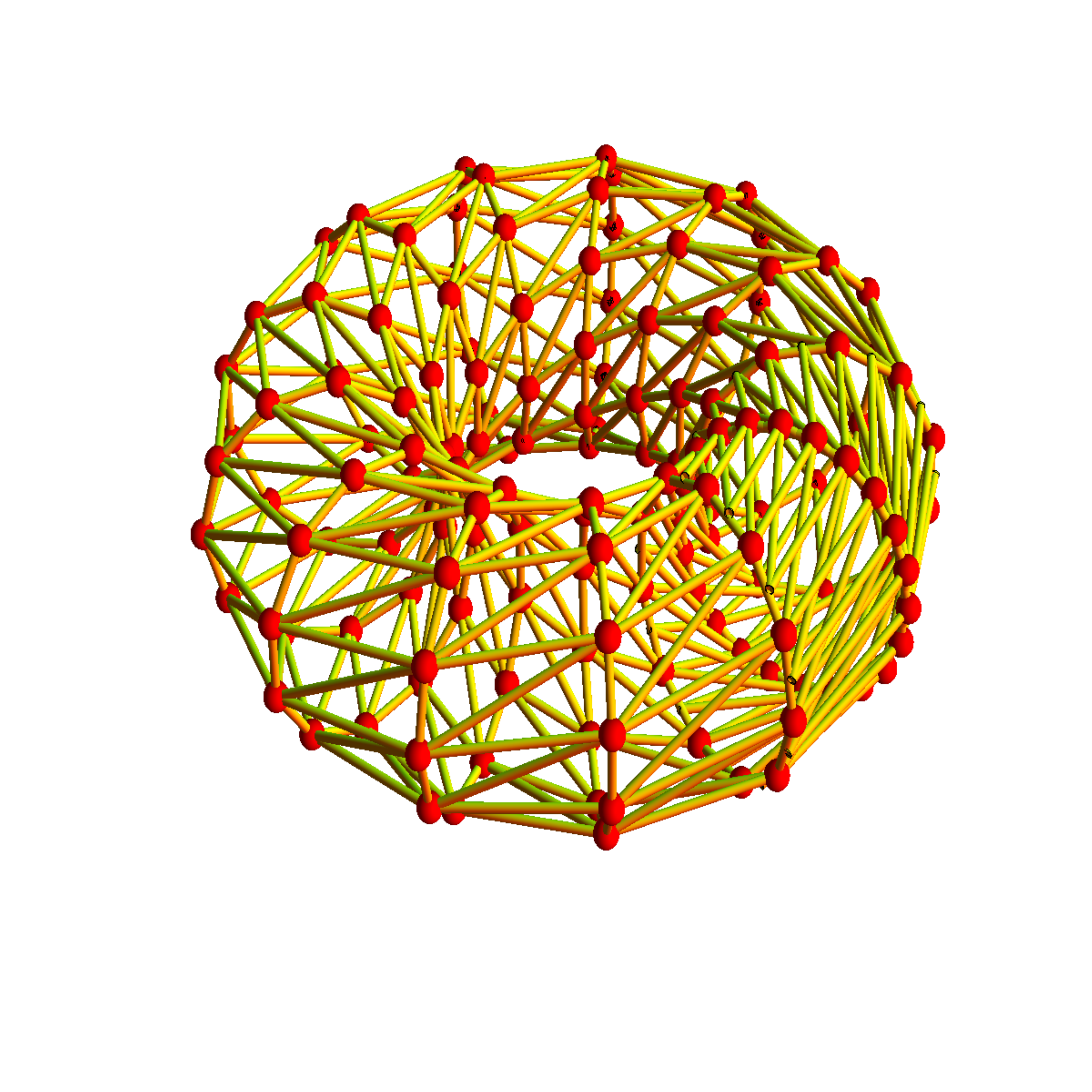}}
\scalebox{0.121}{\includegraphics{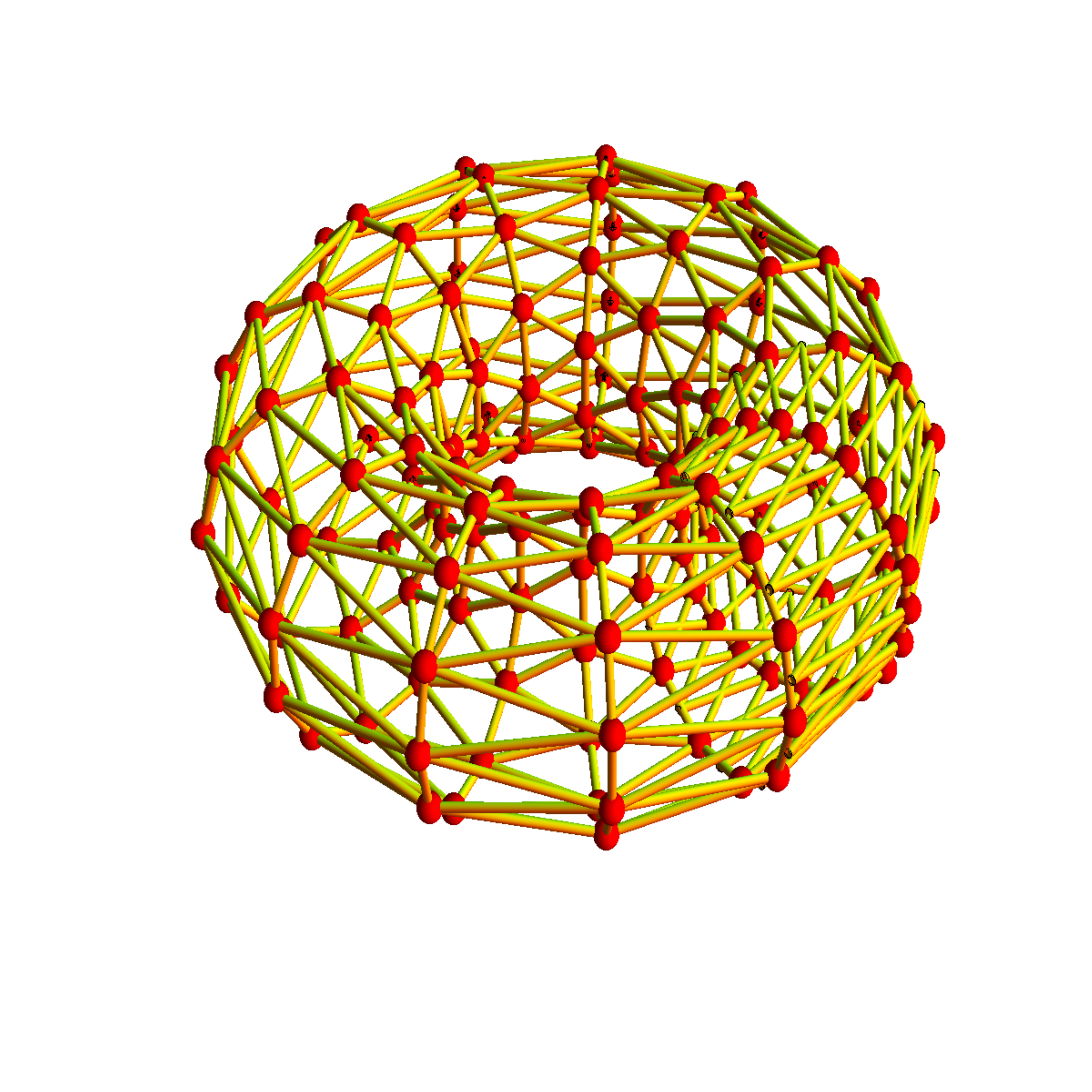}}
} \parbox{6cm}{
A Dehn twisted flat torus and an untwisted torus. The left
and right picture show two non-isomorphic graphs, but they
have the same topological properties and are isospectral
for the Laplacian as well as for the Dirac operator.
It is the easiest example of a pair of nonisometric but
Dirac isospectral graphs. }}

\parbox{16.8cm}{ \parbox{10cm}{
\scalebox{0.121}{\includegraphics{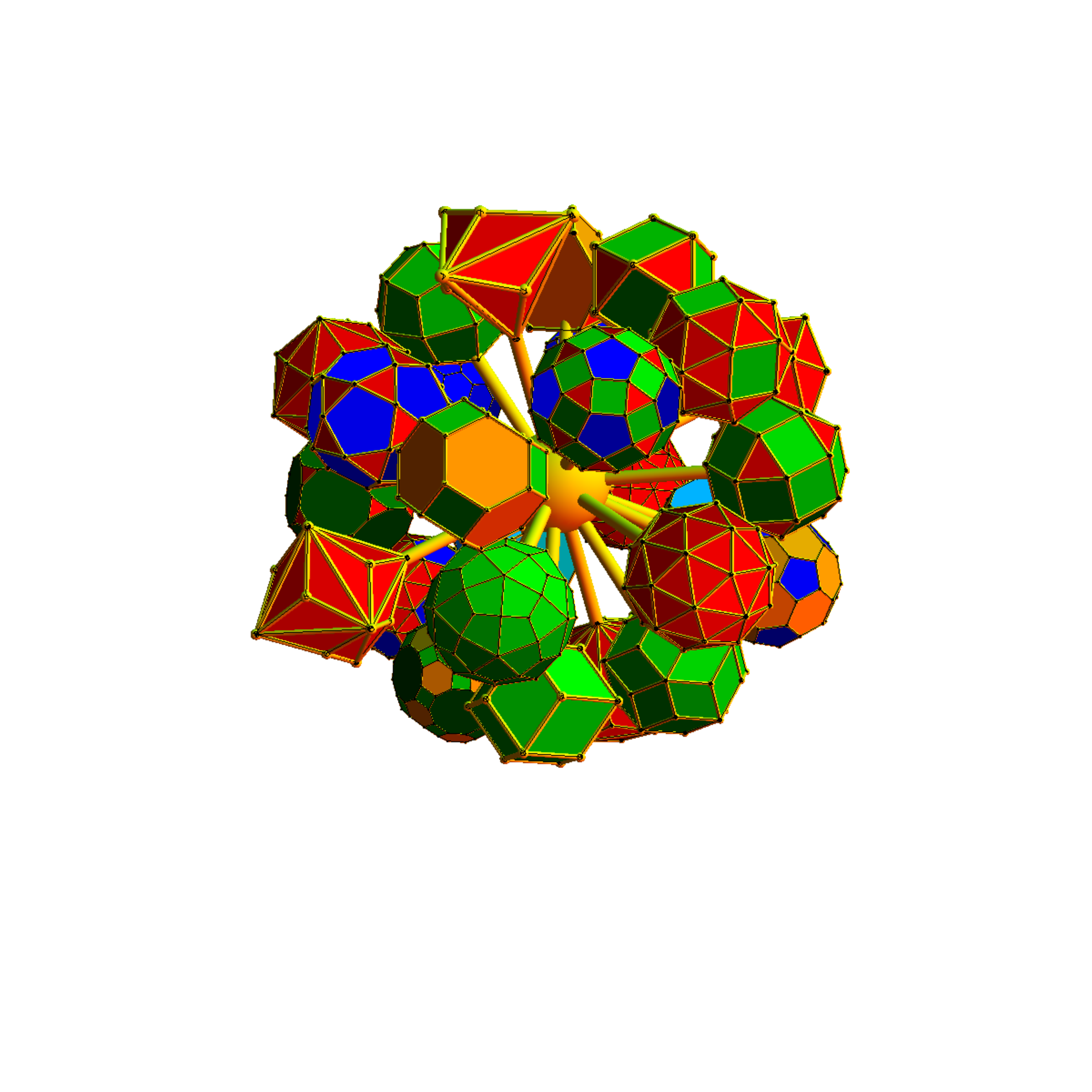}}
\scalebox{0.121}{\includegraphics{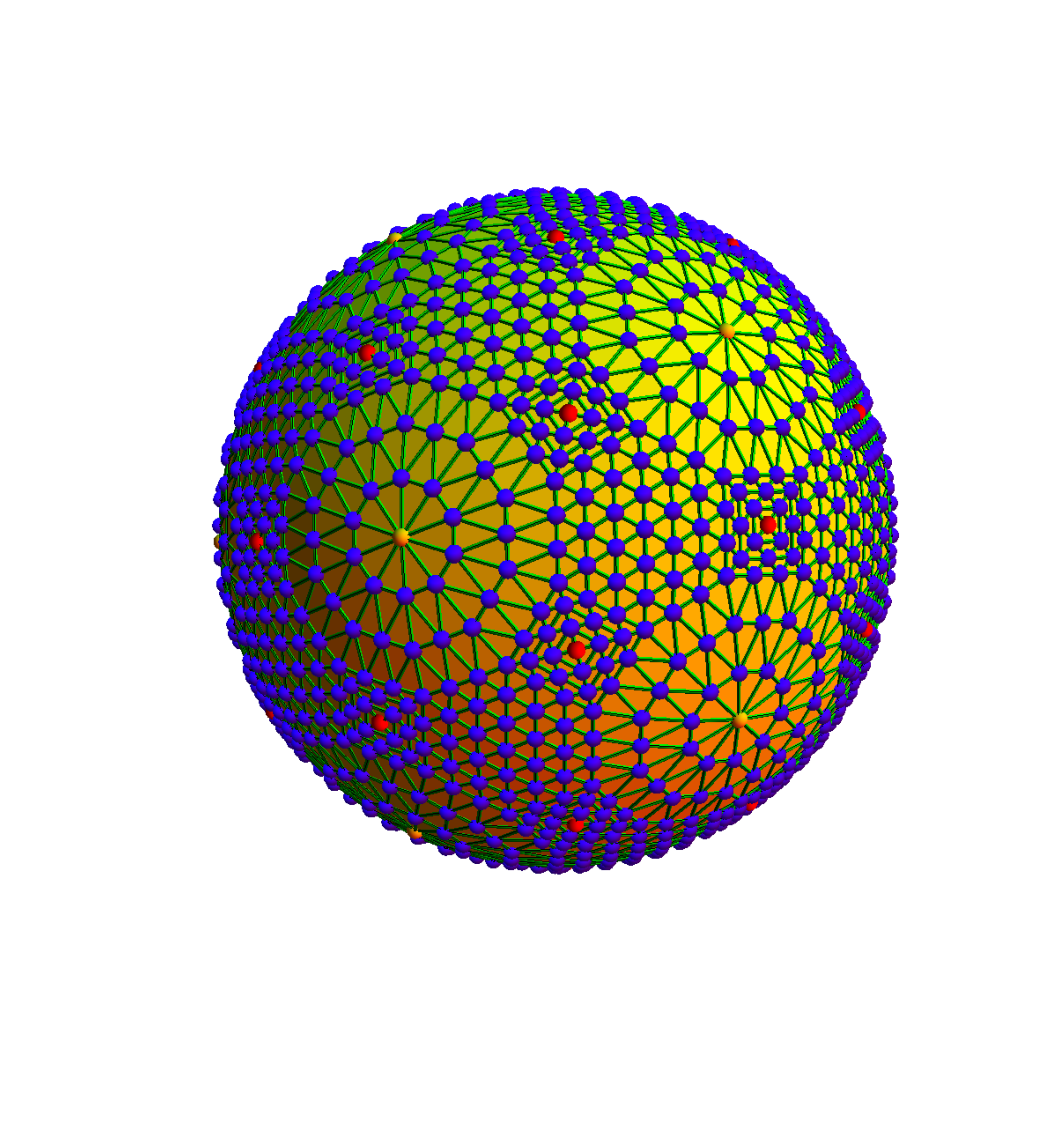}}
} \parbox{6cm}{
All 26 Archimedean and Catalan solids joined to a 
"gem" in the form of a DisdyakisDodecahedron. The gem has since been in the process
of being printed (http://sdu.ictp.it/3d/gem/, http://www.3drucken.ch). 
The right figure shows a Great Rhombicosidodecahedron with 30 points of curvature 
$1/3$ and 12 points of curvature $-2/3$. The total curvature is $2$ and 
agrees with the Euler characteristic. This illustrates a discrete Gauss-Bonnet
theorem \cite{elemente11}. }}

\parbox{16.8cm}{ \parbox{10cm}{
\scalebox{0.121}{\includegraphics{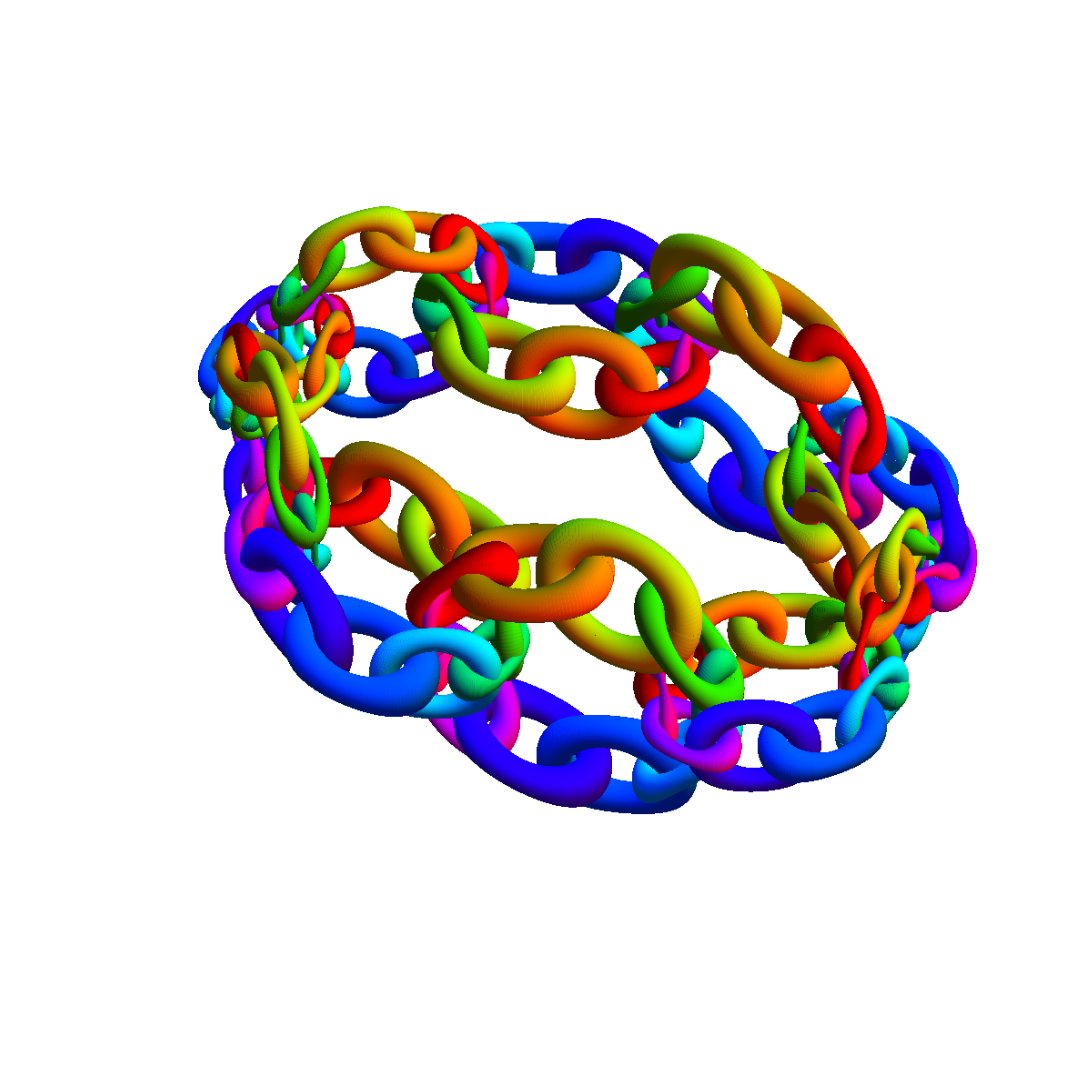}}
\scalebox{0.121}{\includegraphics{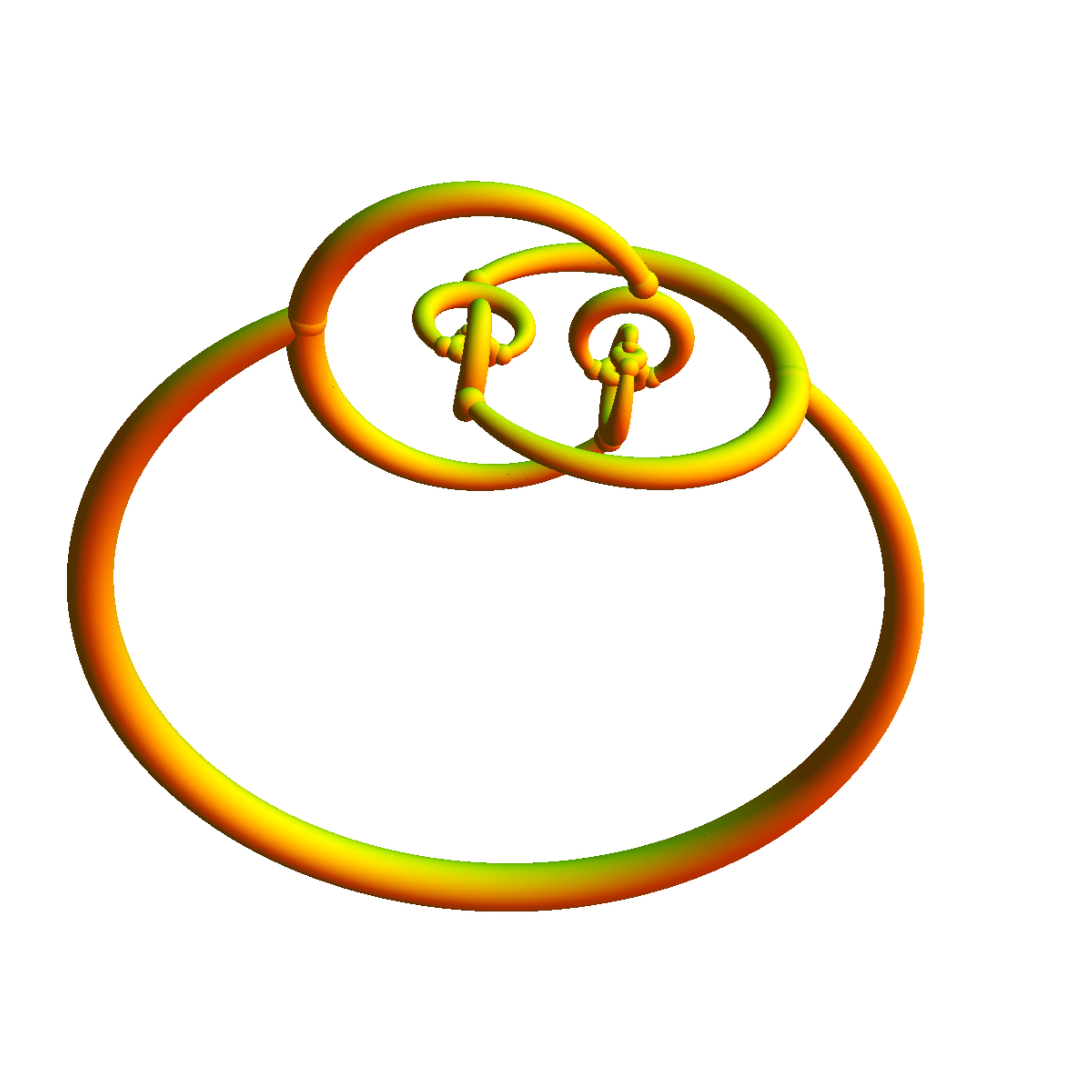}}
} \parbox{6cm}{
The Antoine necklace is a Cantor set in space whose
complement is not simply connected. The Alexander sphere
seen to the right is a topological 3 ball which is simply 
connected but which has an exterior which is not simply connected. 
Alexander spheres also make nice ear rings when printed.
}}

\parbox{16.8cm}{ \parbox{10cm}{
\scalebox{0.12}{\includegraphics{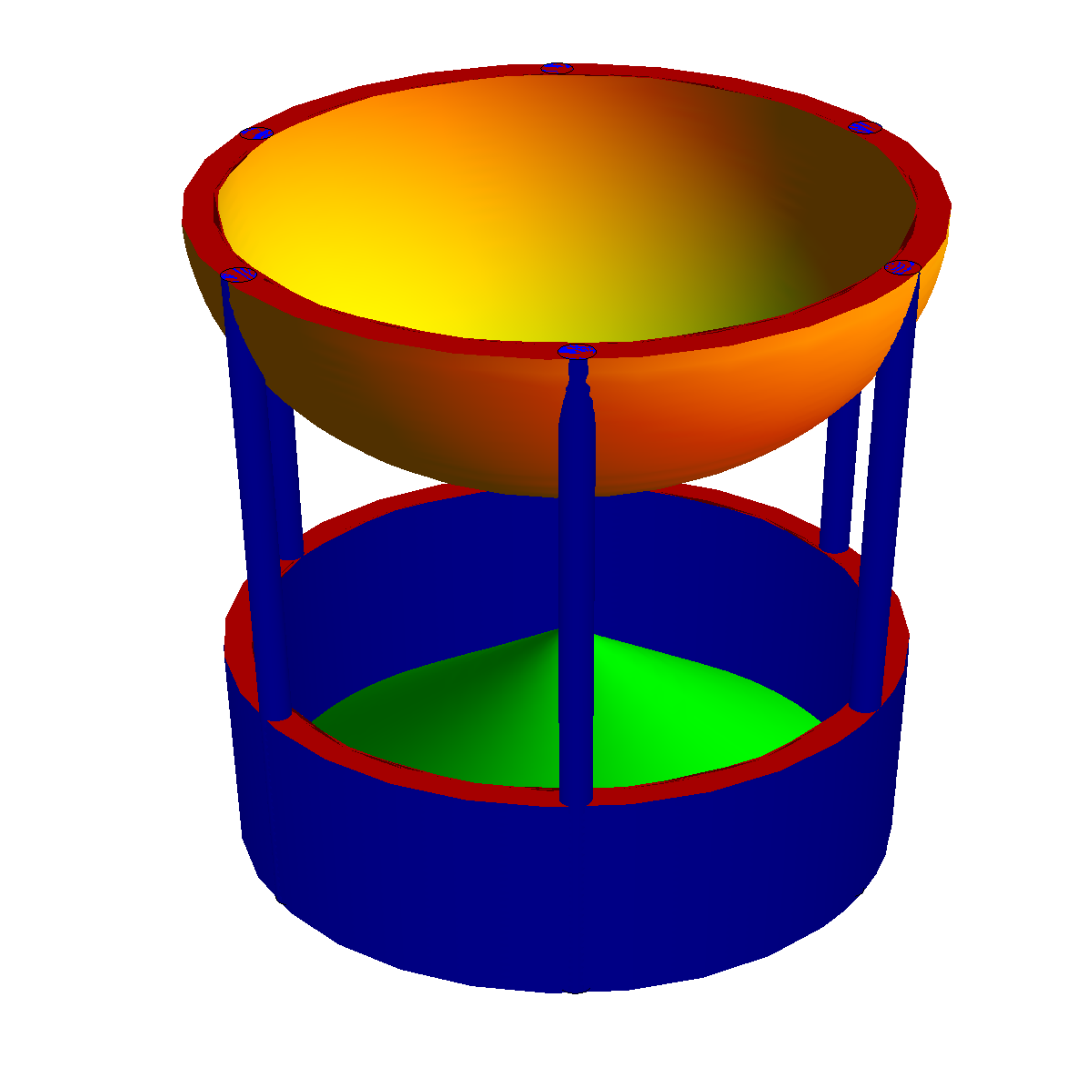}}
\scalebox{0.12}{\includegraphics{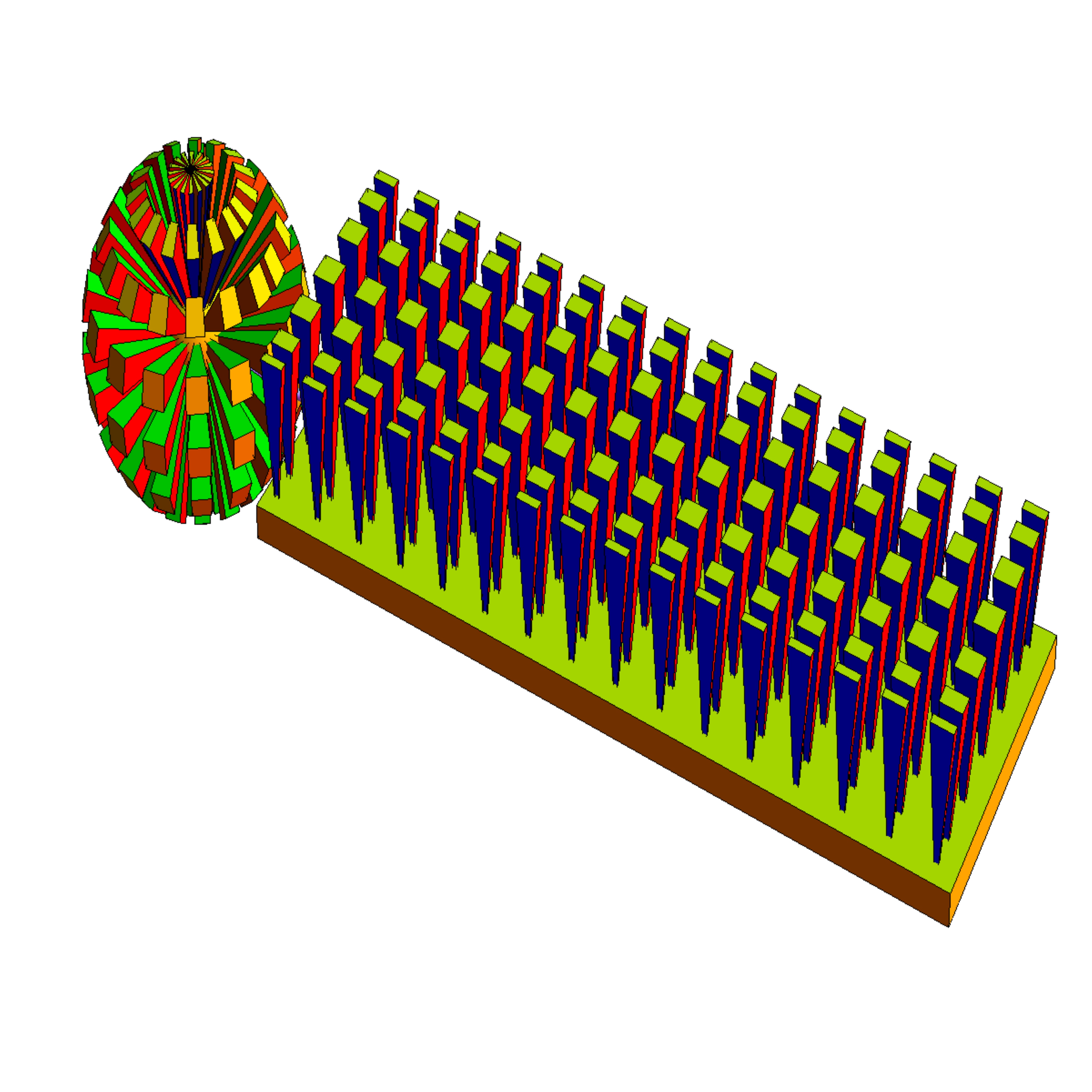}}
} \parbox{6cm}{
Two Archimedes type proofs that the volume of the sphere is $4\pi/3$. \cite{Heath3,Heath4,Thomas}.
The first one assumes that the surface area $A$ is known. The formula $V=A r/3$ can be
seen by cutting up the sphere into many small tetrahedra of volume $dA r/3$.
When summing this over the sphere, we get $A r/3$. The second proof compares the 
half sphere volume with the complement of a cone in a cylinder.  }}

\parbox{16.8cm}{ \parbox{10cm}{
\scalebox{0.12}{\includegraphics{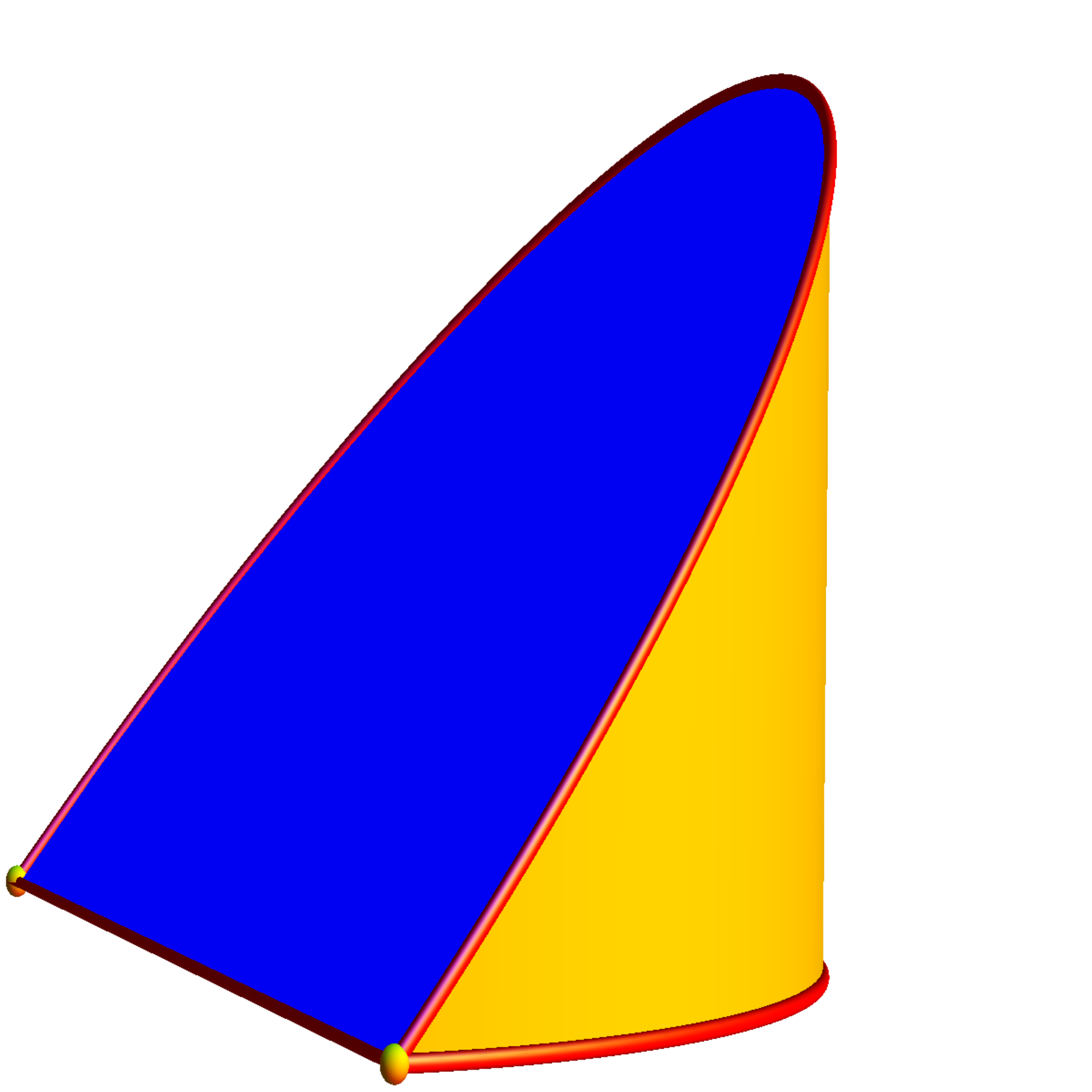}}
\scalebox{0.12}{\includegraphics{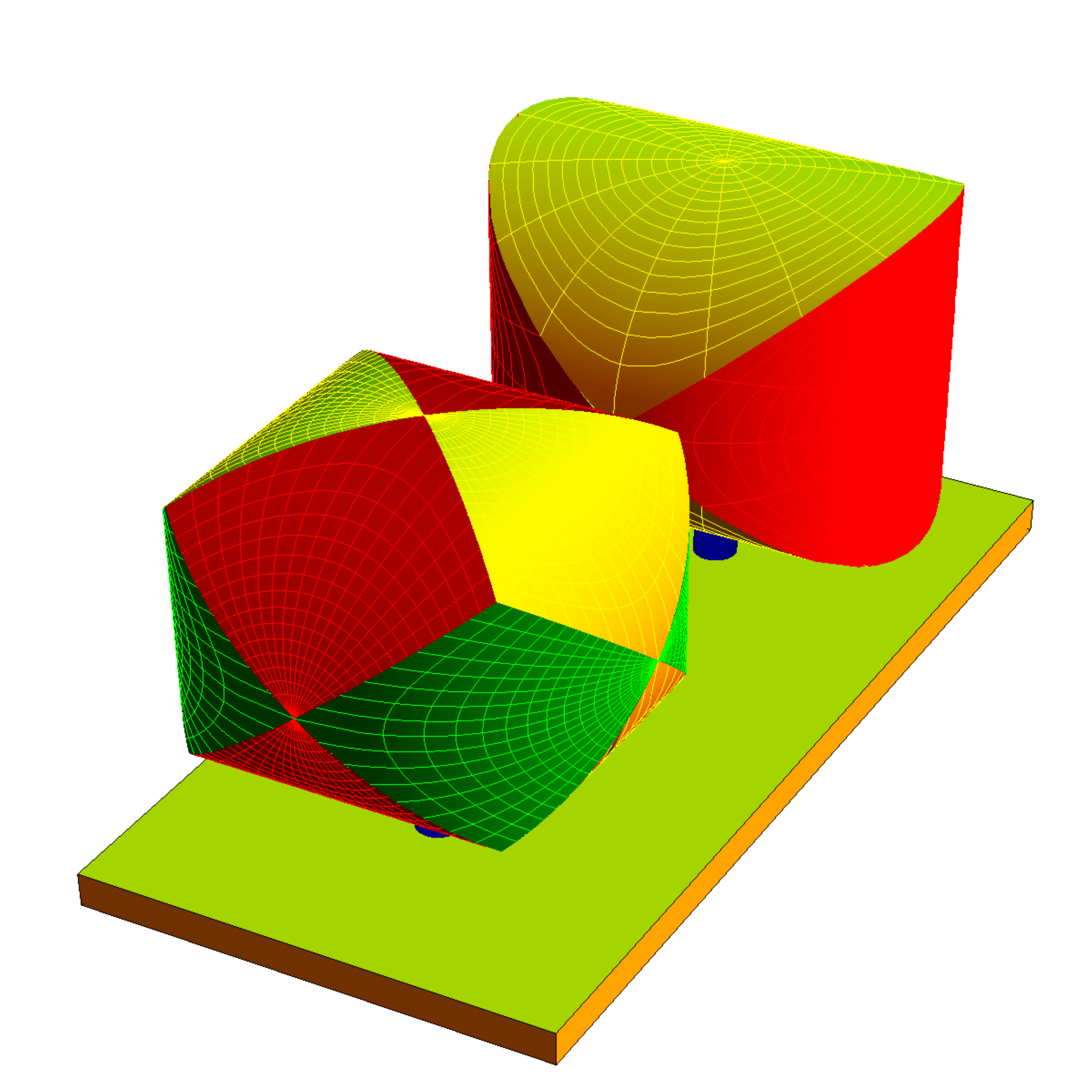}}
} \parbox{6cm}{
The hoof of Archimedes, the Archimedean dome, the intersection 
of cylinders are solids for which Archimedes could compute the volume 
with comparative integration methods \cite{Apostol}. The hoof 
is also an object where Archimedes had to
use a limiting sum, probably the first in the history of humankind
\cite{archimedescodex}.
}}

\parbox{16.8cm}{ \parbox{10cm}{
\scalebox{0.12}{\includegraphics{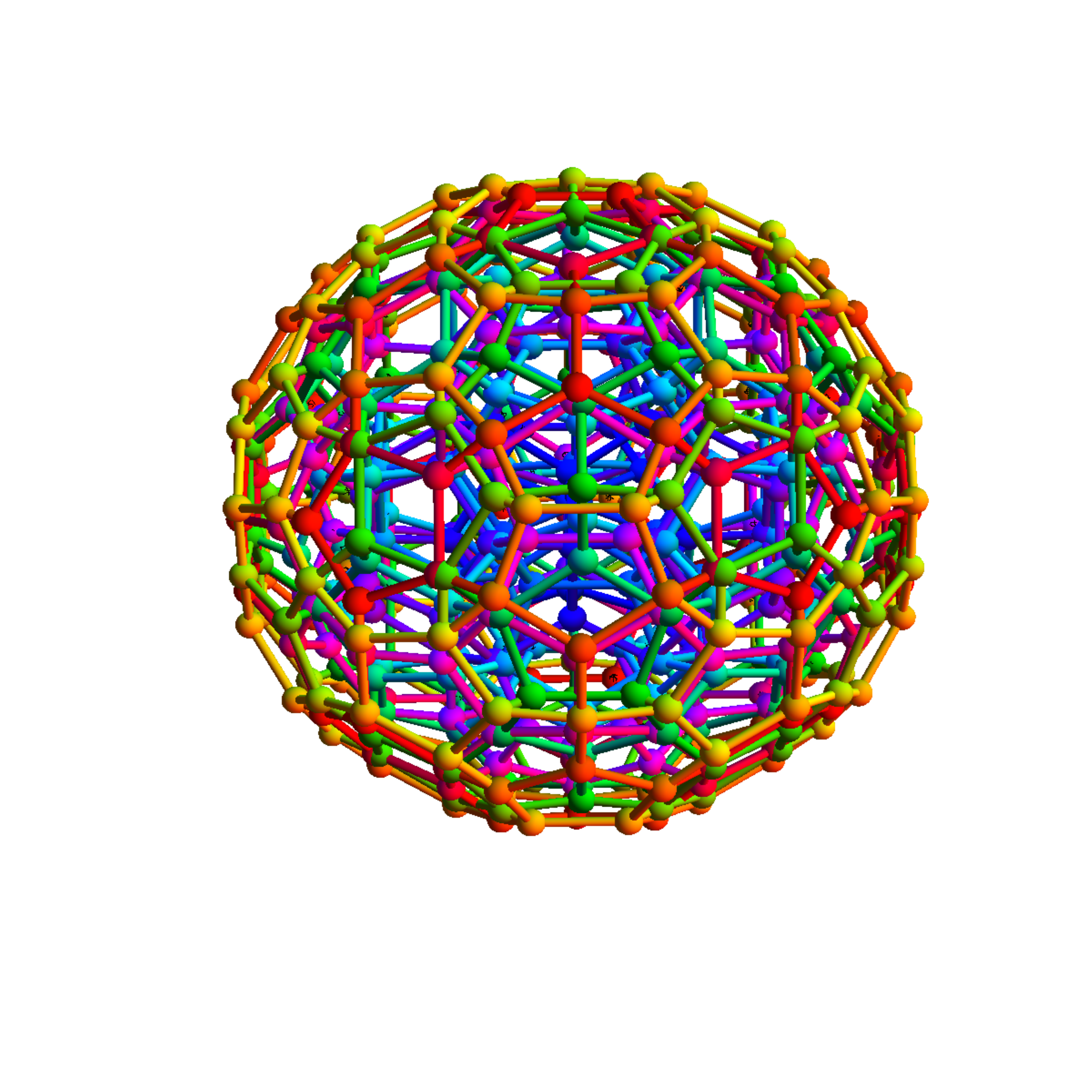}}
\scalebox{0.12}{\includegraphics{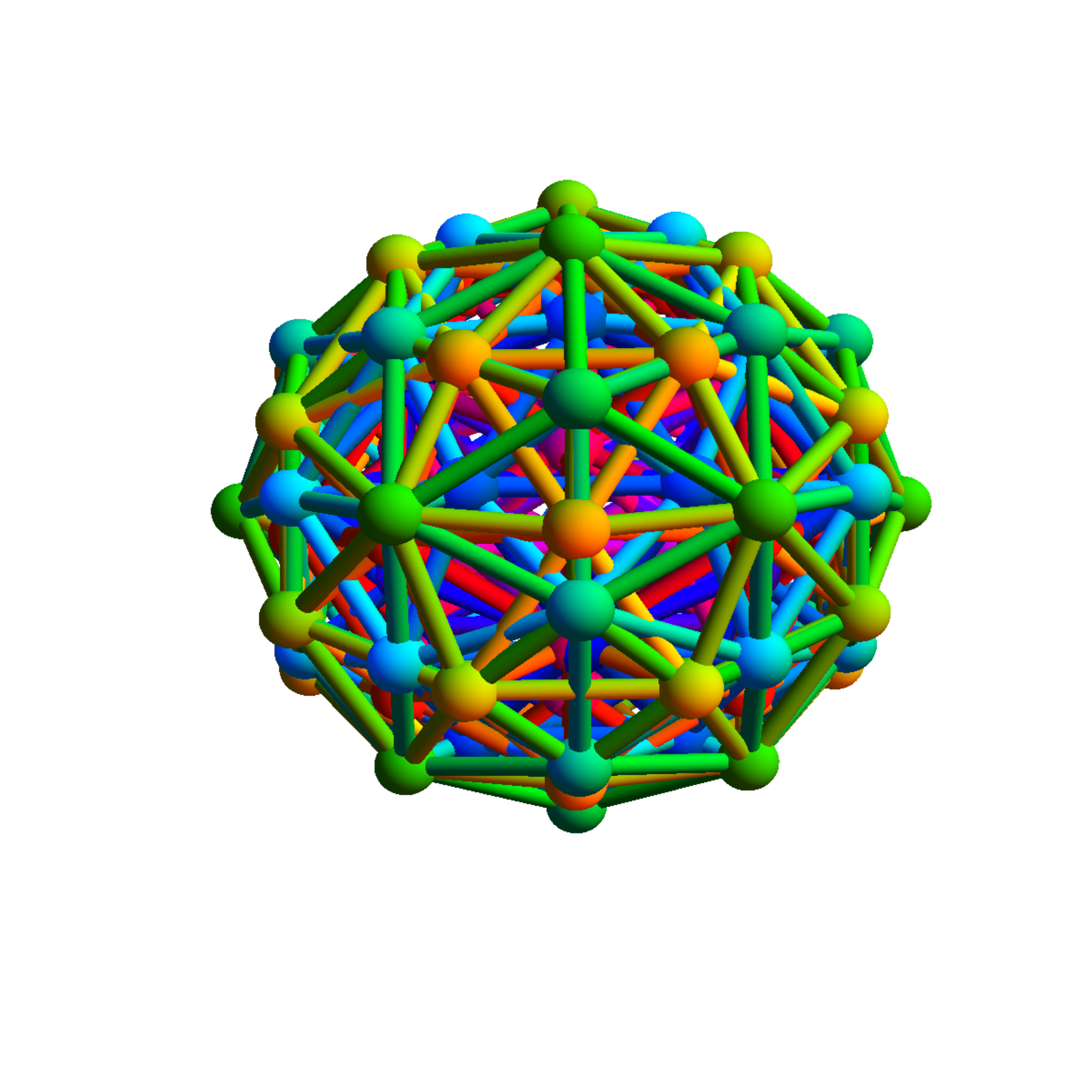}}
} \parbox{6cm}{
Two of the 6 regular convex polytopes in 4 dimensions. 
The color is the height in the four-dimensional space.
We see the 120 cell and the 600 cell. The color of a node
encodes its position in the forth dimension. 
We can not see the fourth dimension but only its shadow
on the cave to speak with an allegory of Plato. The game of 
projecting higher dimensional objects into lower dimensions
has been made a theme of \cite{Flatland}.
}}

\parbox{16.8cm}{ \parbox{10cm}{
\scalebox{0.12}{\includegraphics{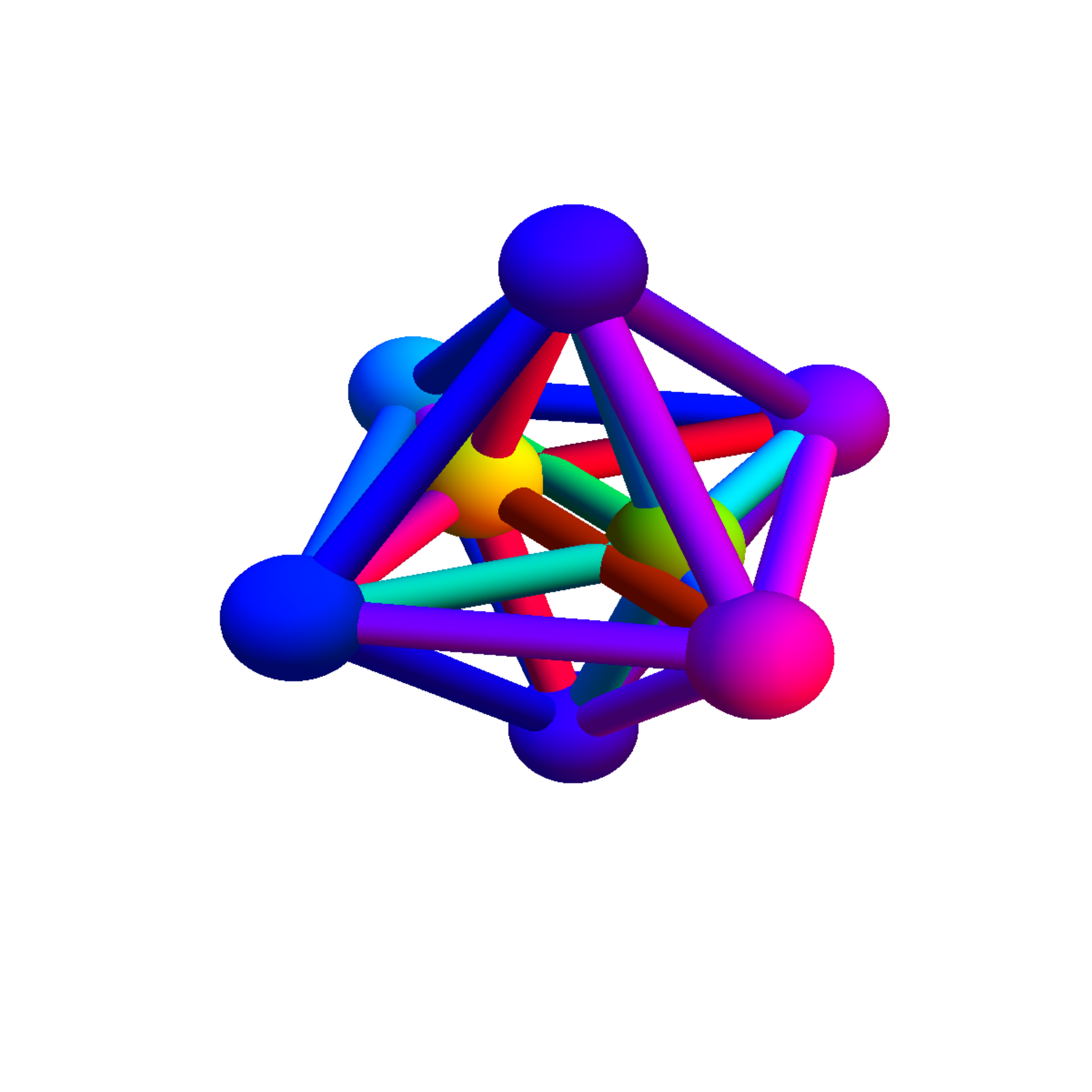}}
\scalebox{0.12}{\includegraphics{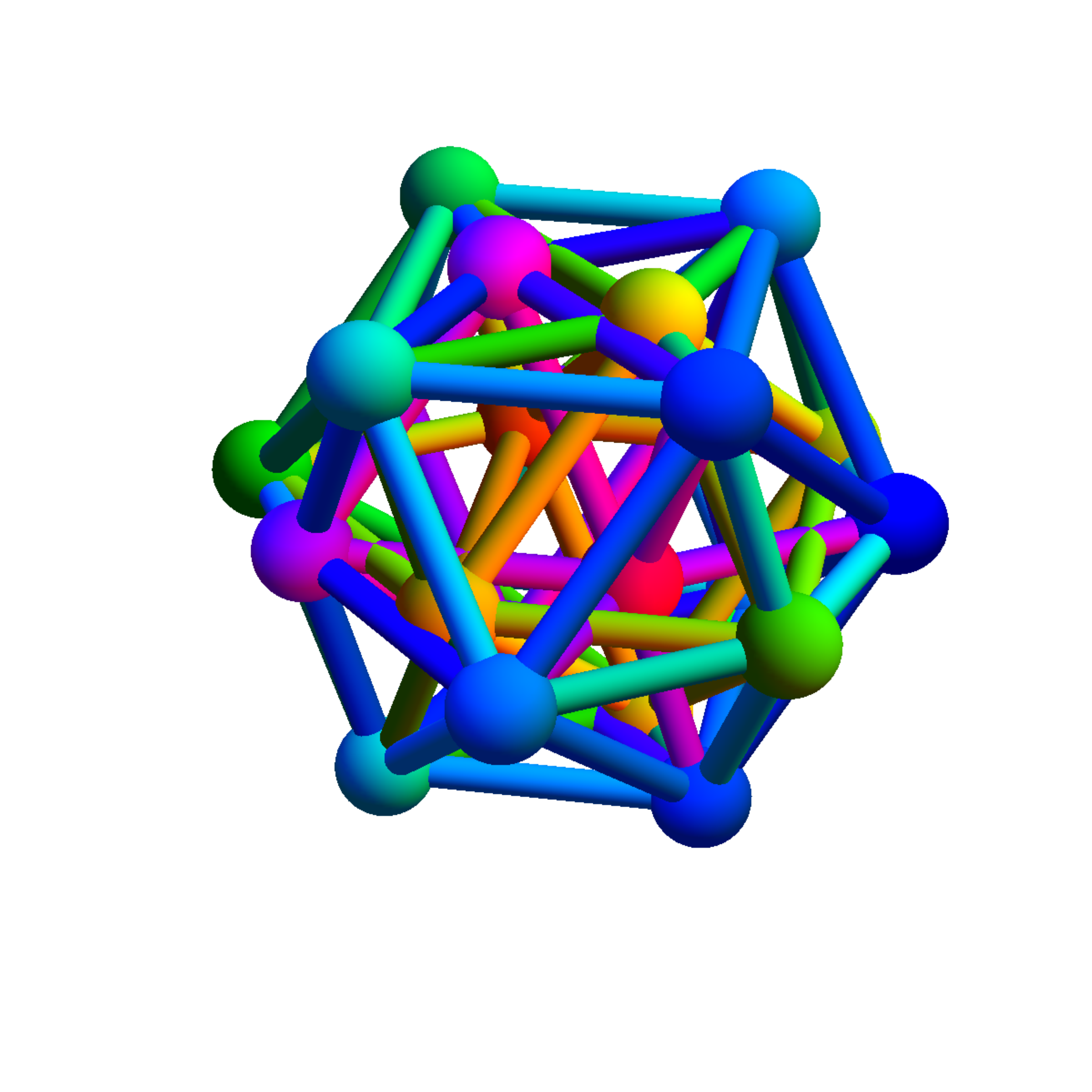}}
} \parbox{6cm}{
An other pair of the 6 regular convex polytopes in 4 dimensions.
The color is the height in the four-dimensional space.
We see the 16 cell (the analogue of the octahedron) and the 24 cell. 
The later allows to tessellate 4 dimensional Euclidean space similarly
as the octahedron can tessellate 3 dimensional Euclidean space. 
}}

\parbox{16.8cm}{ \parbox{10cm}{
\scalebox{0.12}{\includegraphics{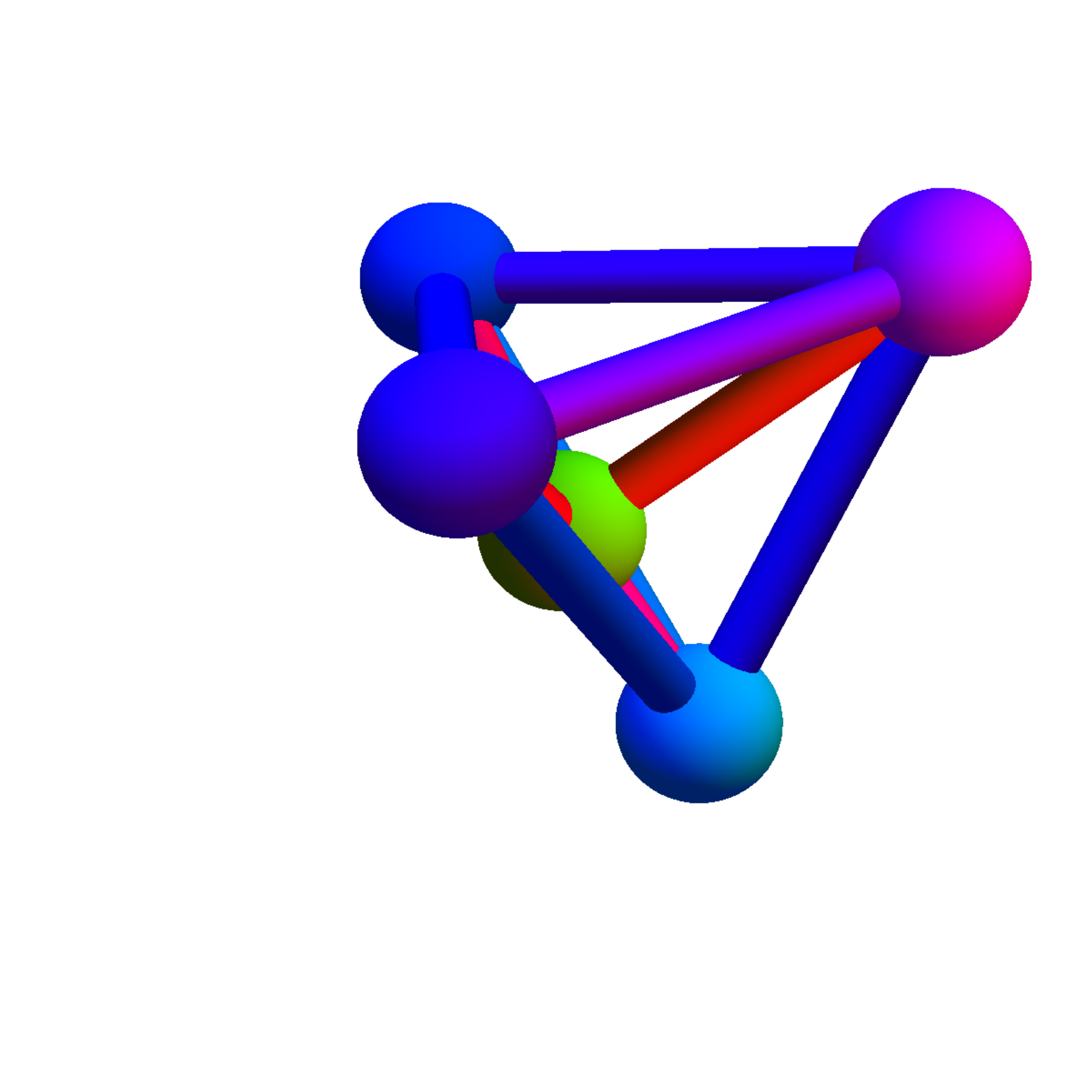}}
\scalebox{0.12}{\includegraphics{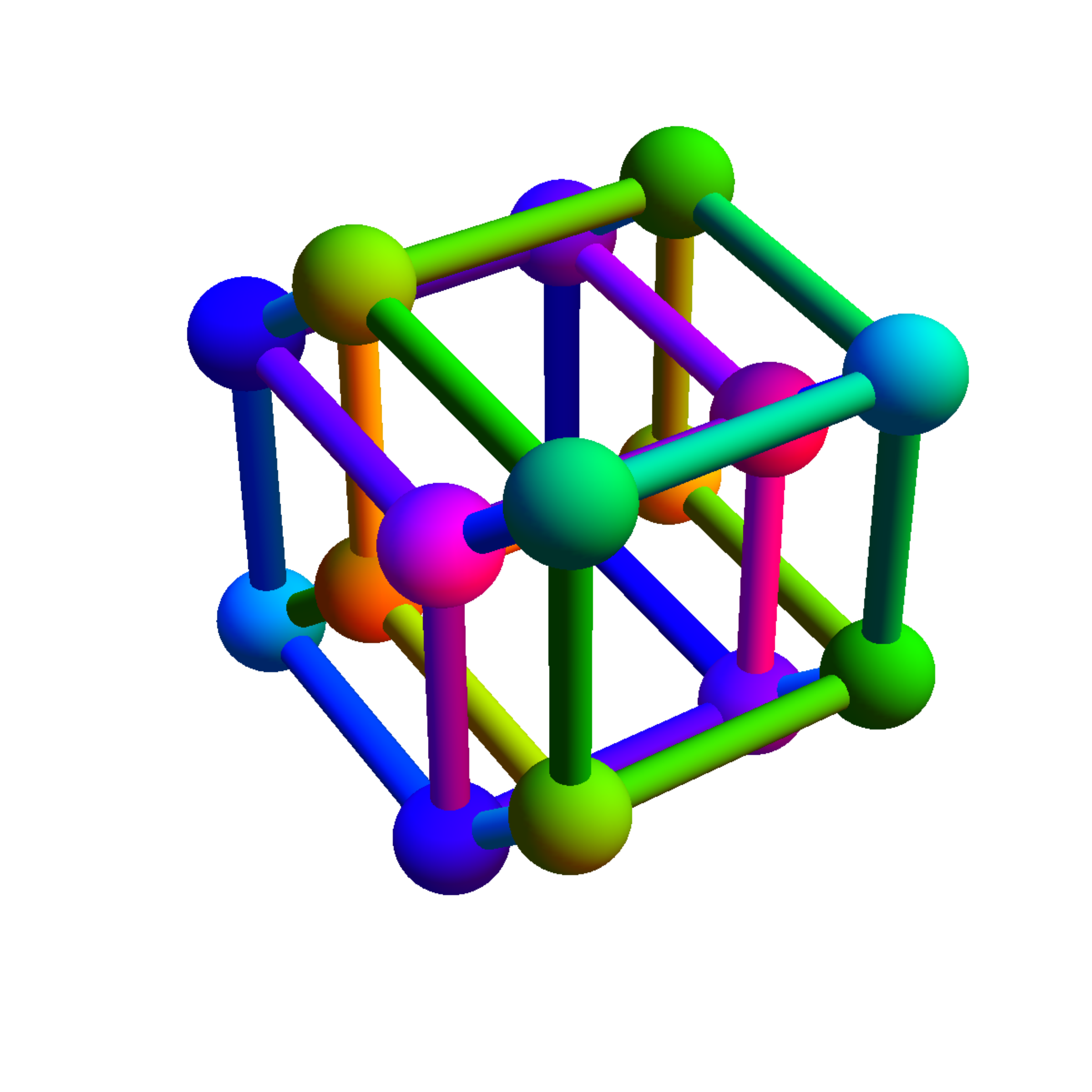}}
} \parbox{6cm}{
The 5 cell is the complete graph with $5$ vertices and the simplest 
4 dimensional polytop. The 8 cell to the right is also called the
tesseract or hypercube. It is the 4 dimensional analogue of the cube and has
reached stardom status among all mathematical objects. 
}}

\parbox{16.8cm}{ \parbox{10cm}{
\scalebox{0.12}{\includegraphics{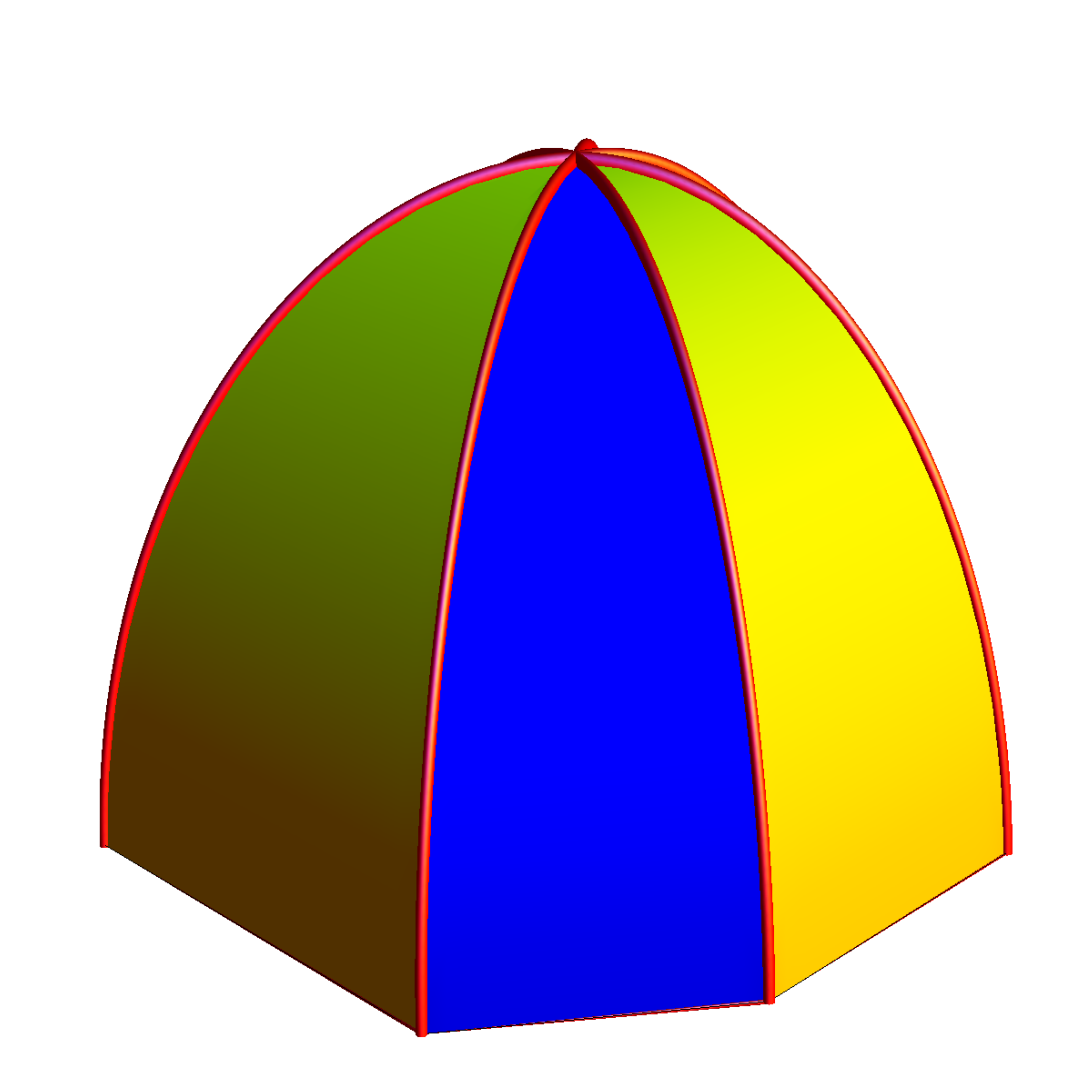}}
\scalebox{0.12}{\includegraphics{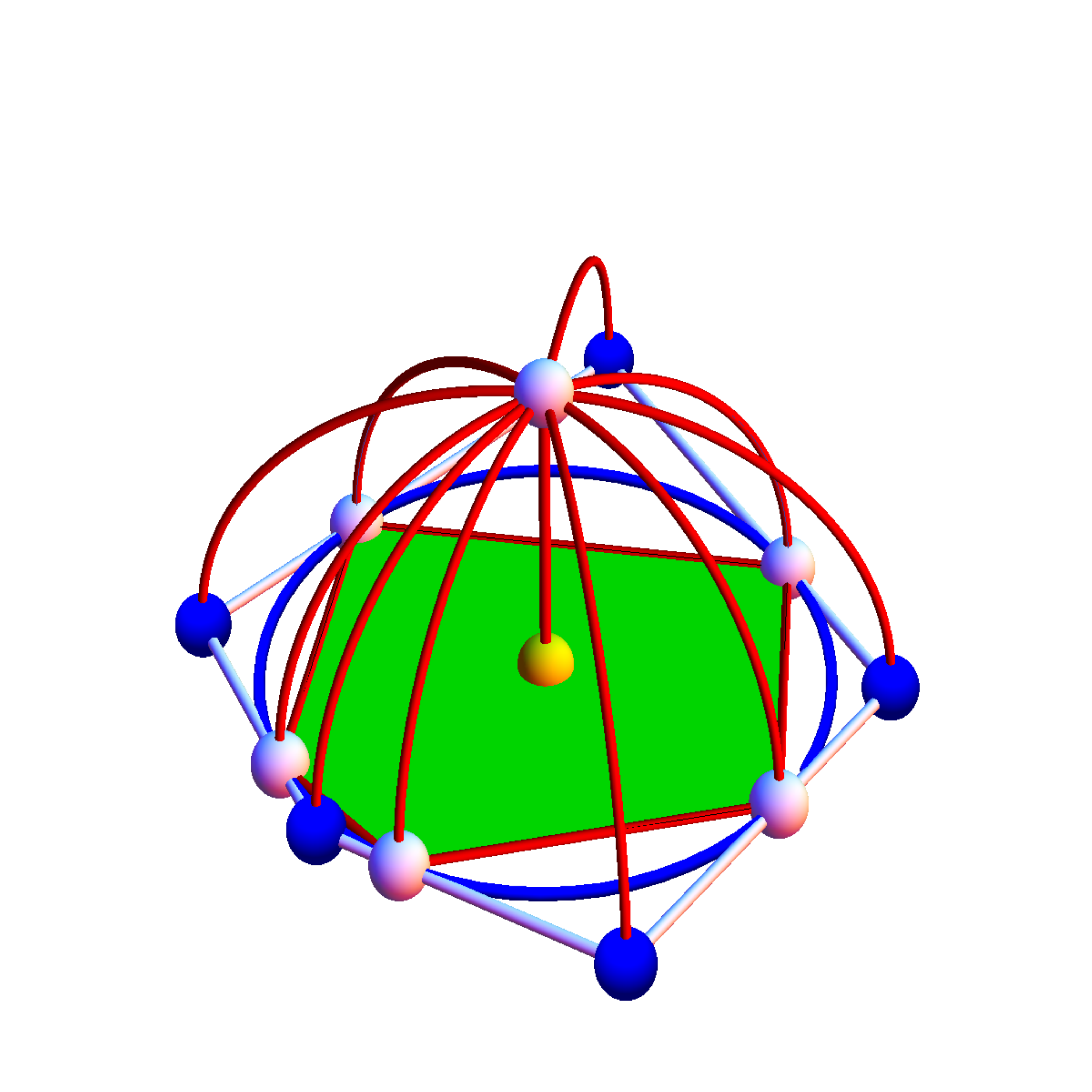}}
} \parbox{6cm}{
An Archimedean dome is half of an Archimedean sphere. 
Archimedean domes have volume equal to {\bf $2/3$} of the prism in which they are
inscribed. It was discovered only later that for Archimedean globes, 
the surface area is {\bf $2/3$} of the surface area of a circumscribing 
prism \cite{Apostol}.
}}

\parbox{16.8cm}{ \parbox{10cm}{
\scalebox{0.12}{\includegraphics{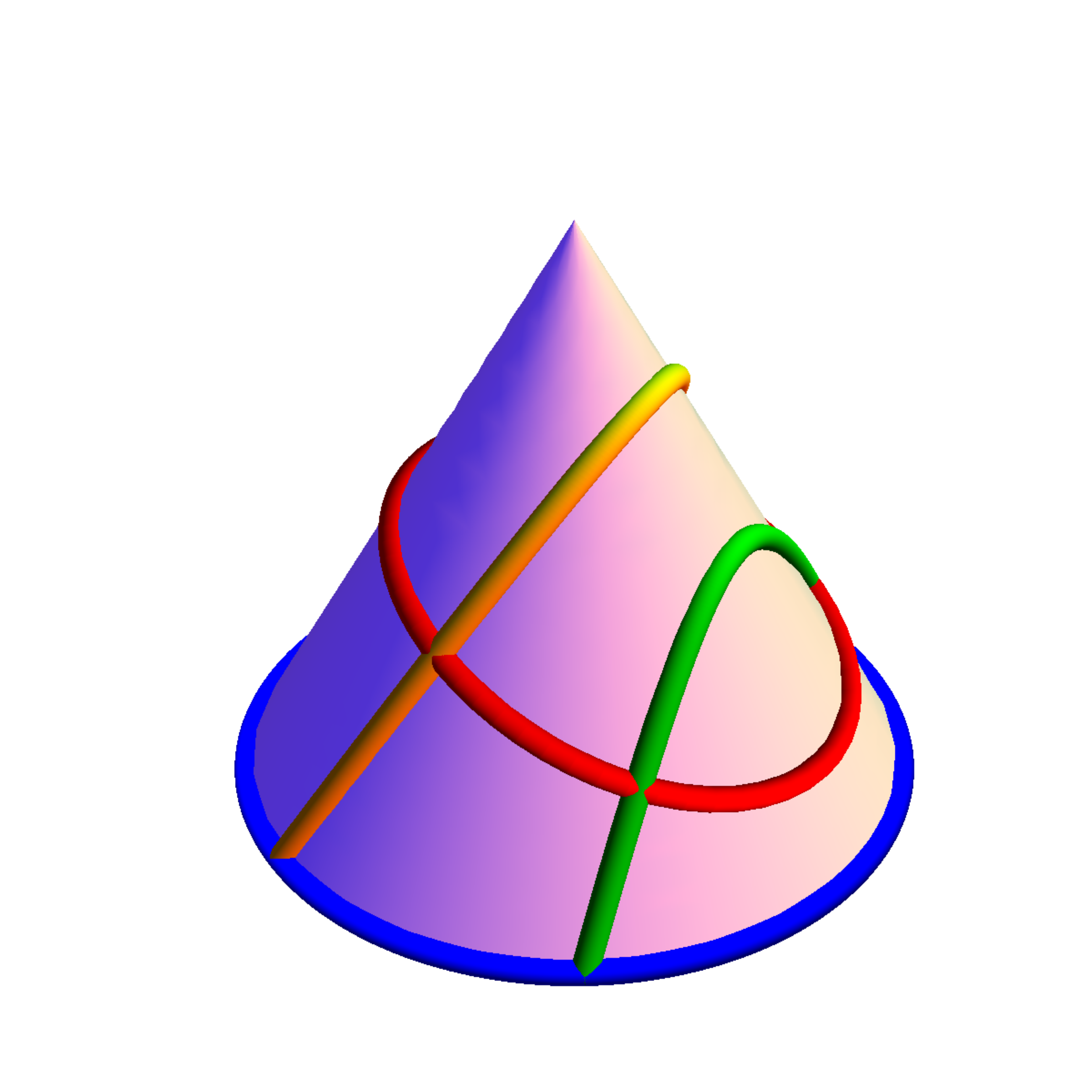}}
\scalebox{0.12}{\includegraphics{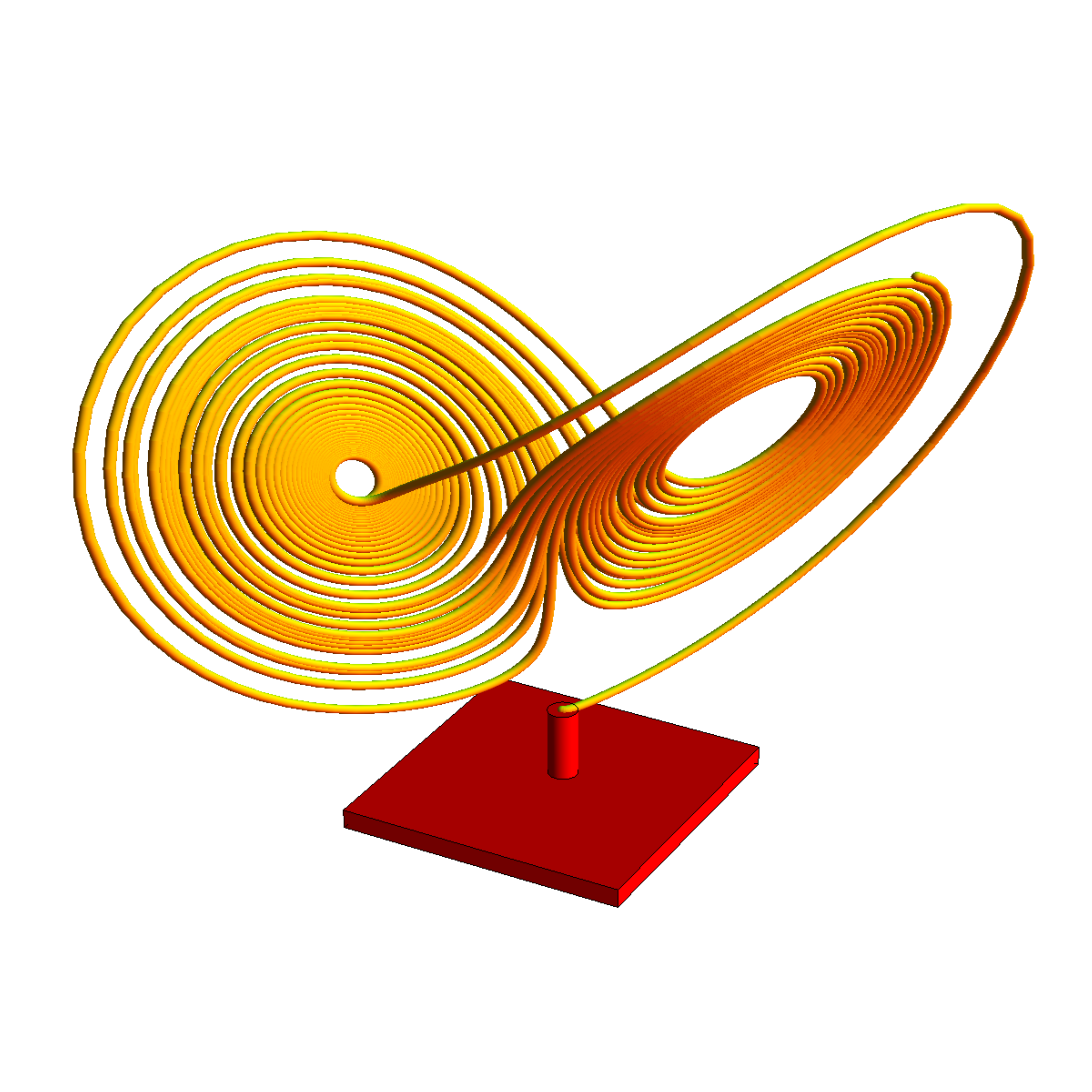}}
} \parbox{6cm}{
An Apollonian cone named after 
Apollonius of Perga is used to visualize the conic sections.
Wooden models appear in school rooms.
The right figure shows an icon of Chaos, the Lorentz attractor
\cite{Sparrow}. It is believed to be a fractal. The dynamics on
this set is chaotic for various parameters. In an appendix, we 
have the source code for computing the attractor. 
}}

\parbox{16.8cm}{ \parbox{10cm}{
\scalebox{0.12}{\includegraphics{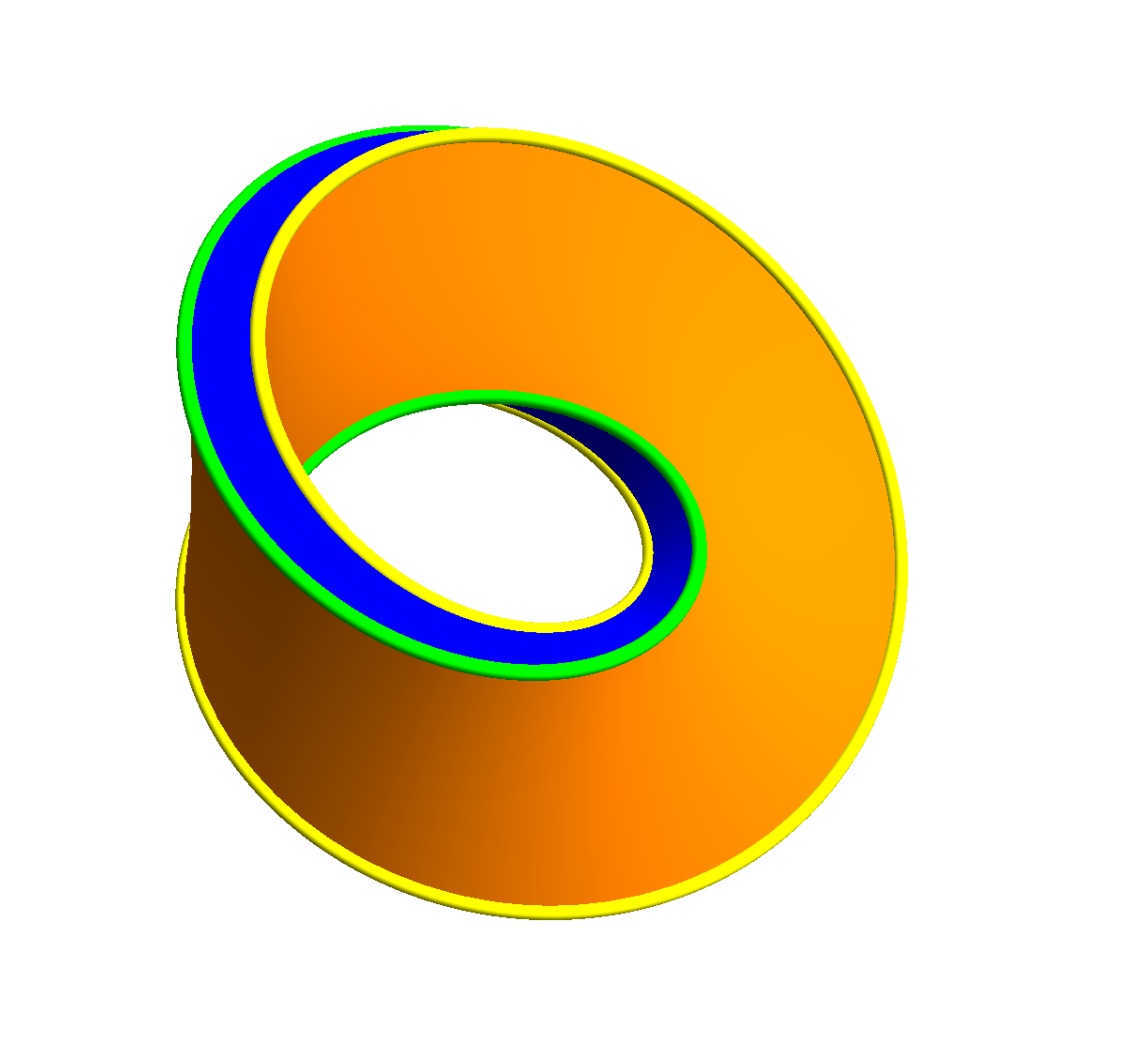}}
\scalebox{0.12}{\includegraphics{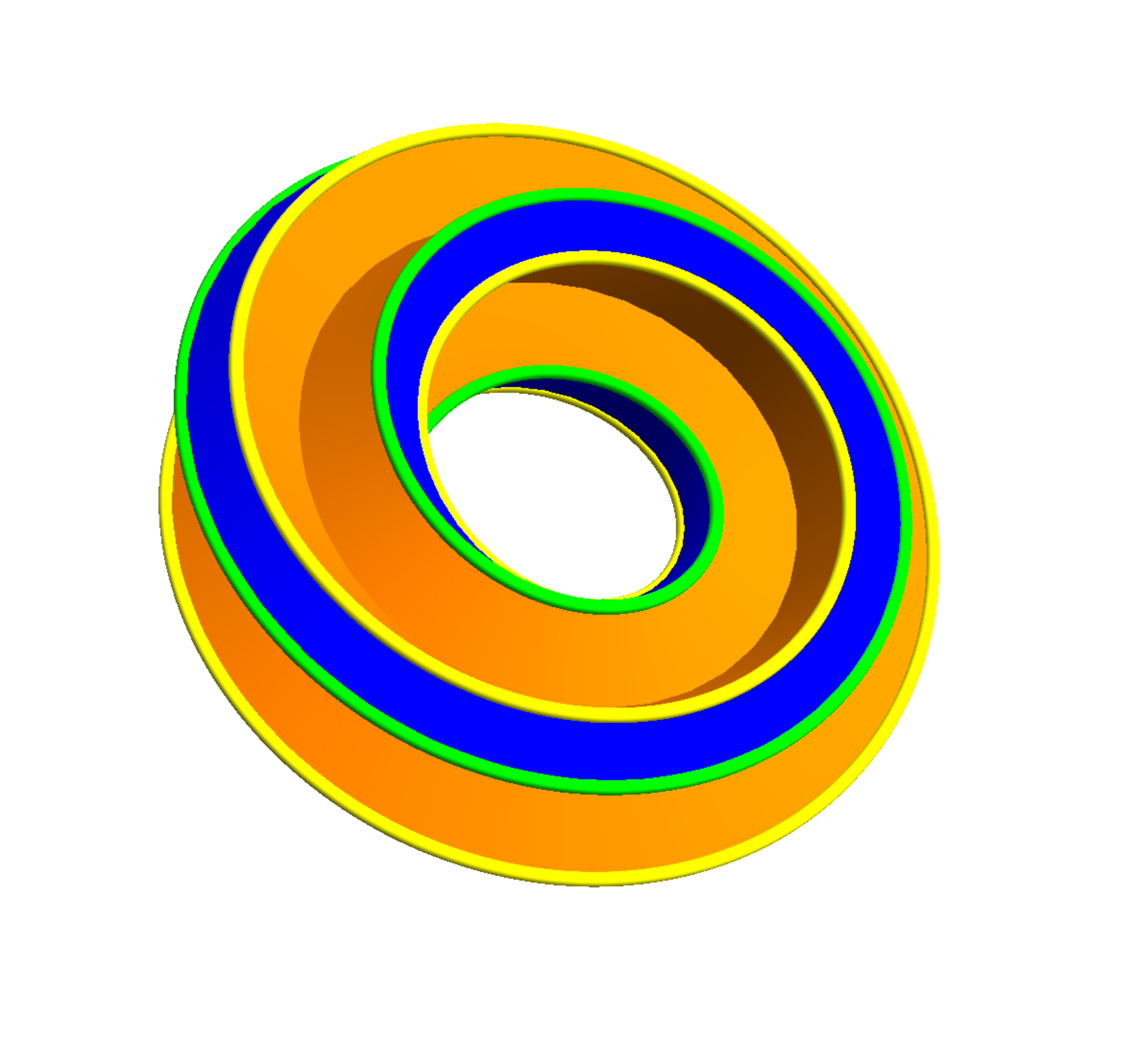}}
} \parbox{6cm}{
The M\"obius strip was thickened so that it can be printed. 
The next picture shows a M\"obius strip with self intersection. 
This is a situation where the computer algebra system shines. 
To make a surface thicker, we have have to compute the normal
vector at every point of the surface. 
}}

\parbox{16.8cm}{ \parbox{10cm}{
\scalebox{0.12}{\includegraphics{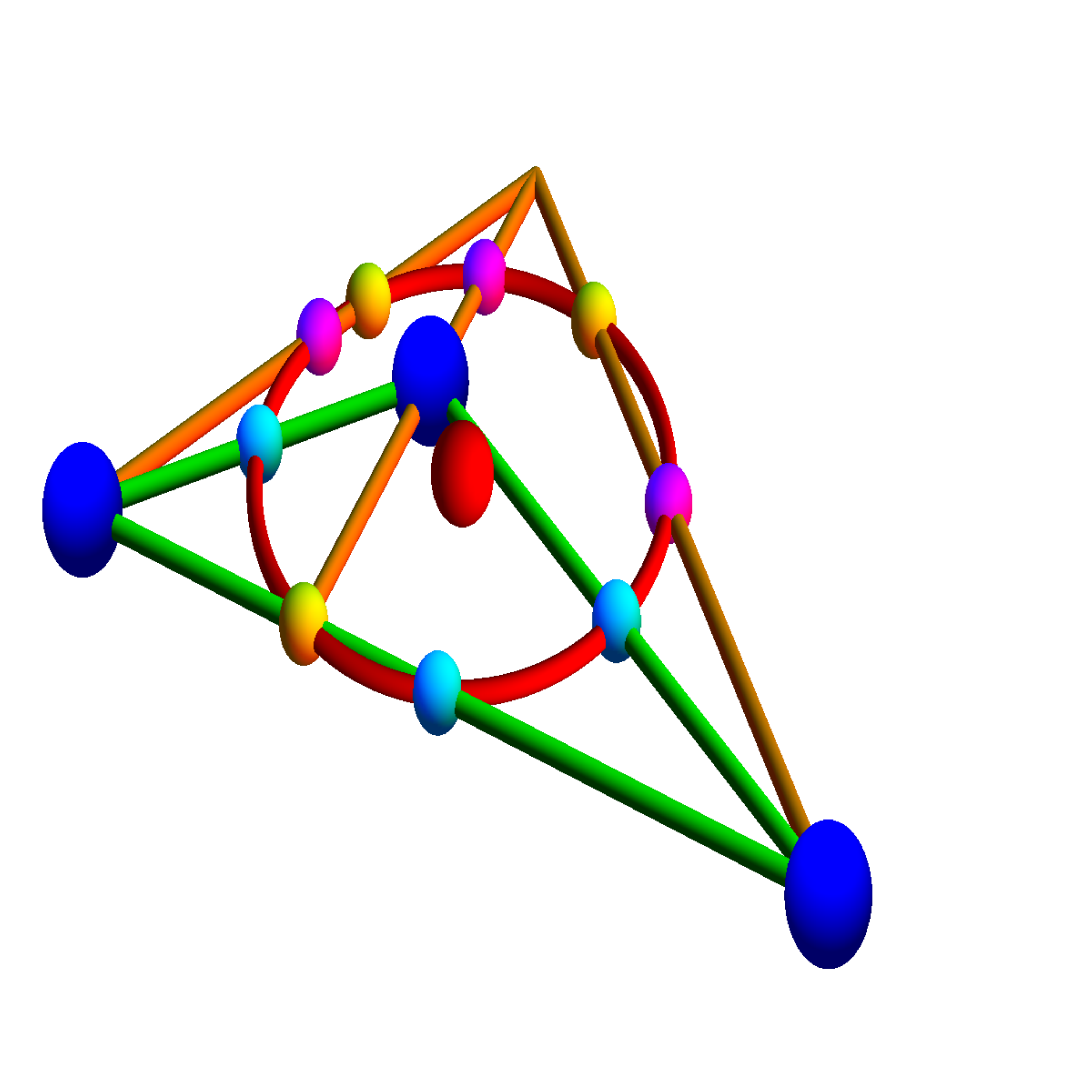}}
\scalebox{0.12}{\includegraphics{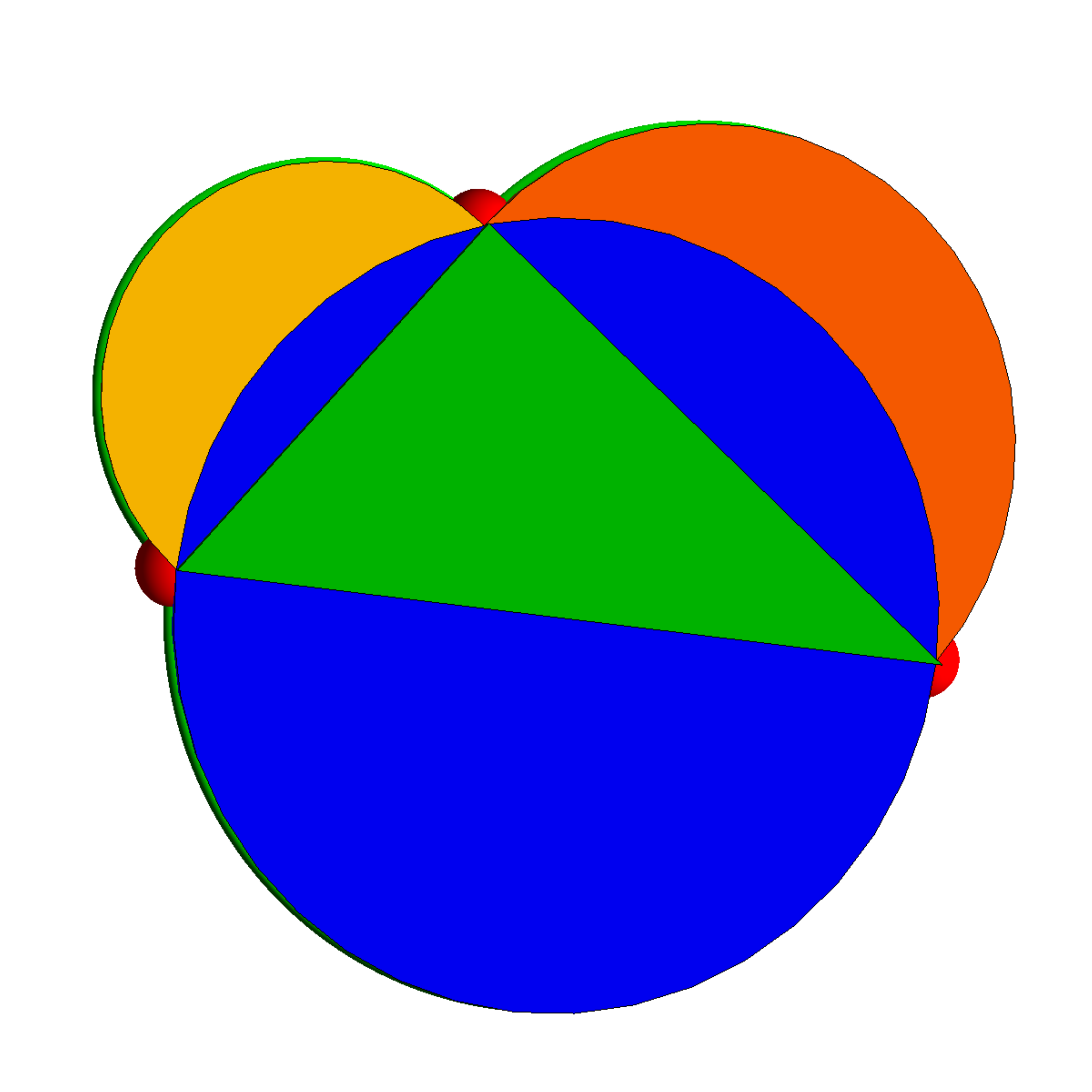}}
} \parbox{6cm}{
The nine point theorem of Feuerbach realized in 3D. 
The right figure illustrates the theorem of Hippocrates, 
an attempt to the quadrature of the circle. The triangle
has the same area than the two moon shaped figures together.
Even so this is planimetry, we have produced three dimensional
objects which could be printed - or worn as a pin.
}}

\parbox{16.8cm}{ \parbox{10cm}{
\scalebox{0.12}{\includegraphics{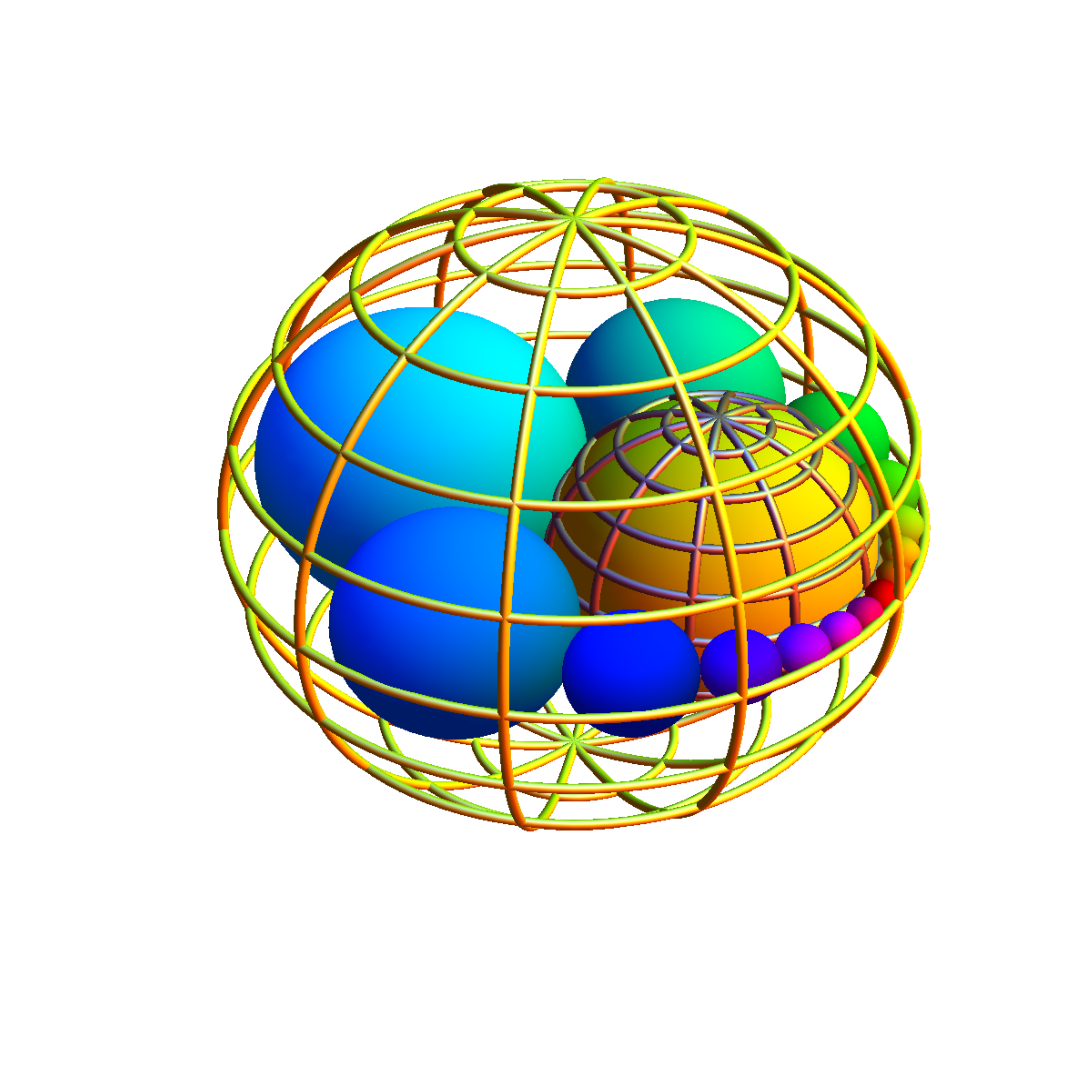}}
\scalebox{0.12}{\includegraphics{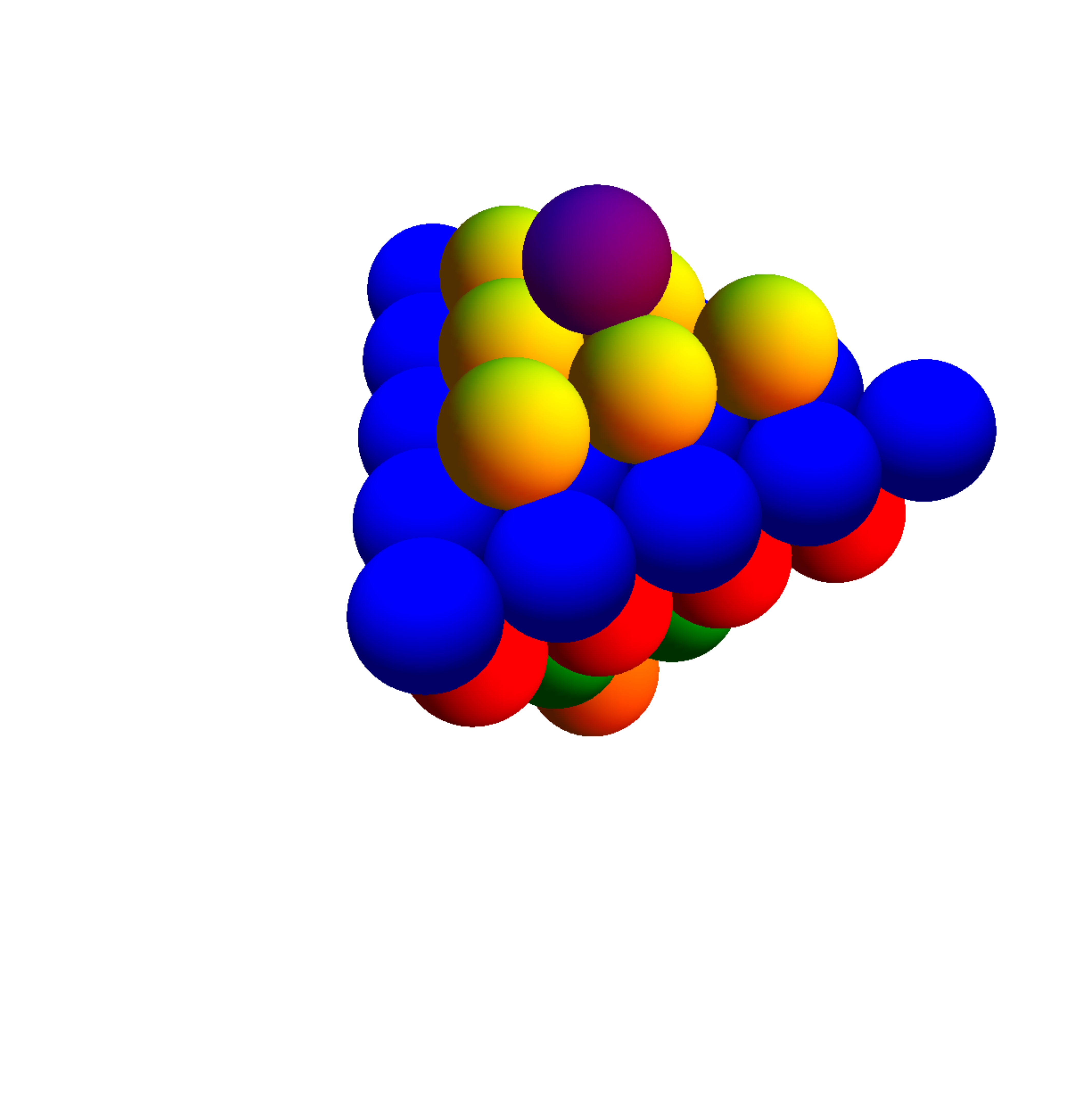}}
} \parbox{6cm}{
The left figure shows Soddy's hexlet. One needs conformal transformations,
M\"obius transformations in particular to construct this solid. 
The right figure hopes to illustrate that there are infinitely many densest packings in
space. While there is a cubic close packing and
hexagonal close packing, these packings can be mixed.
}}

\parbox{16.8cm}{ \parbox{10cm}{
\scalebox{0.12}{\includegraphics{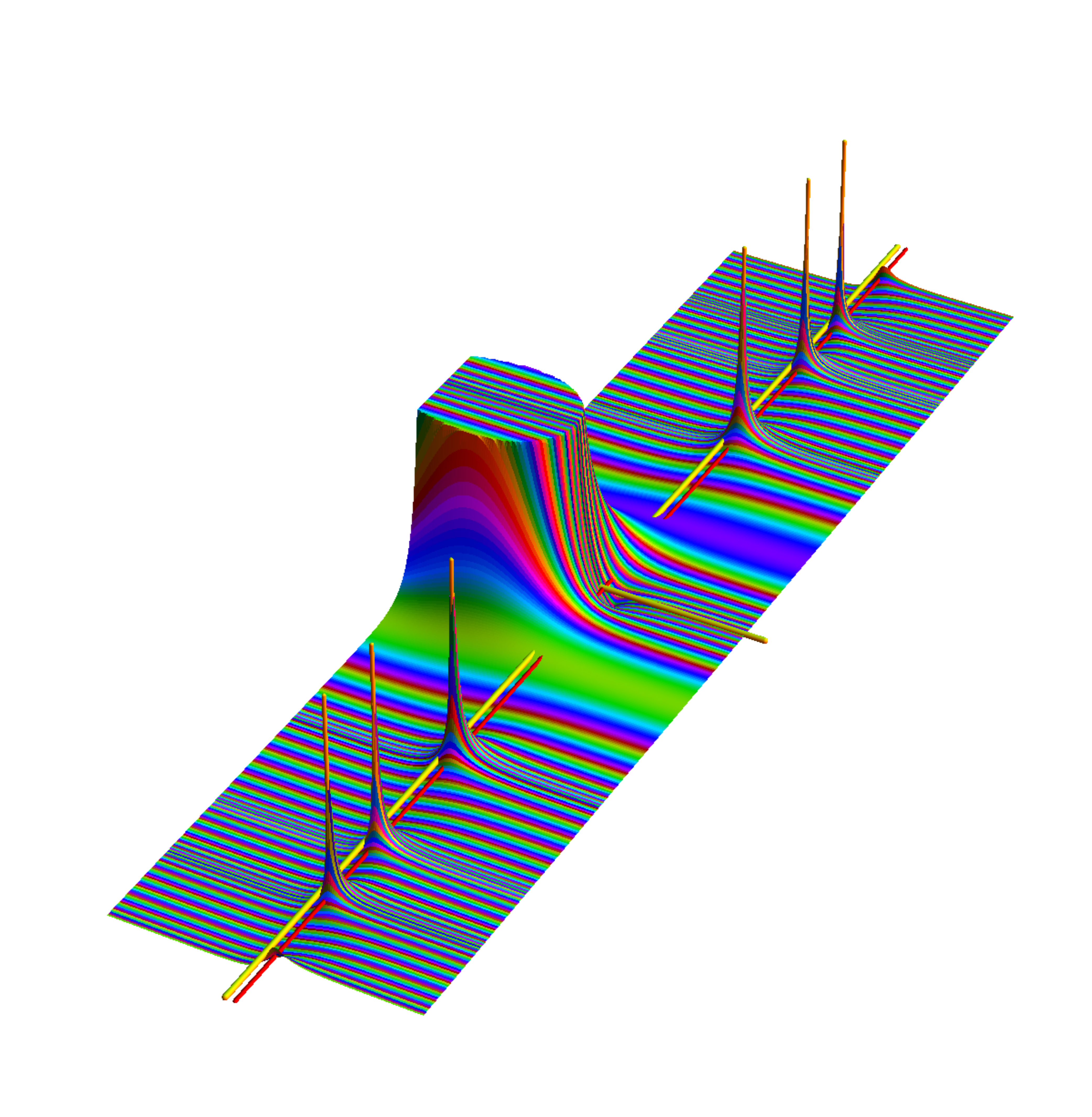}}
\scalebox{0.12}{\includegraphics{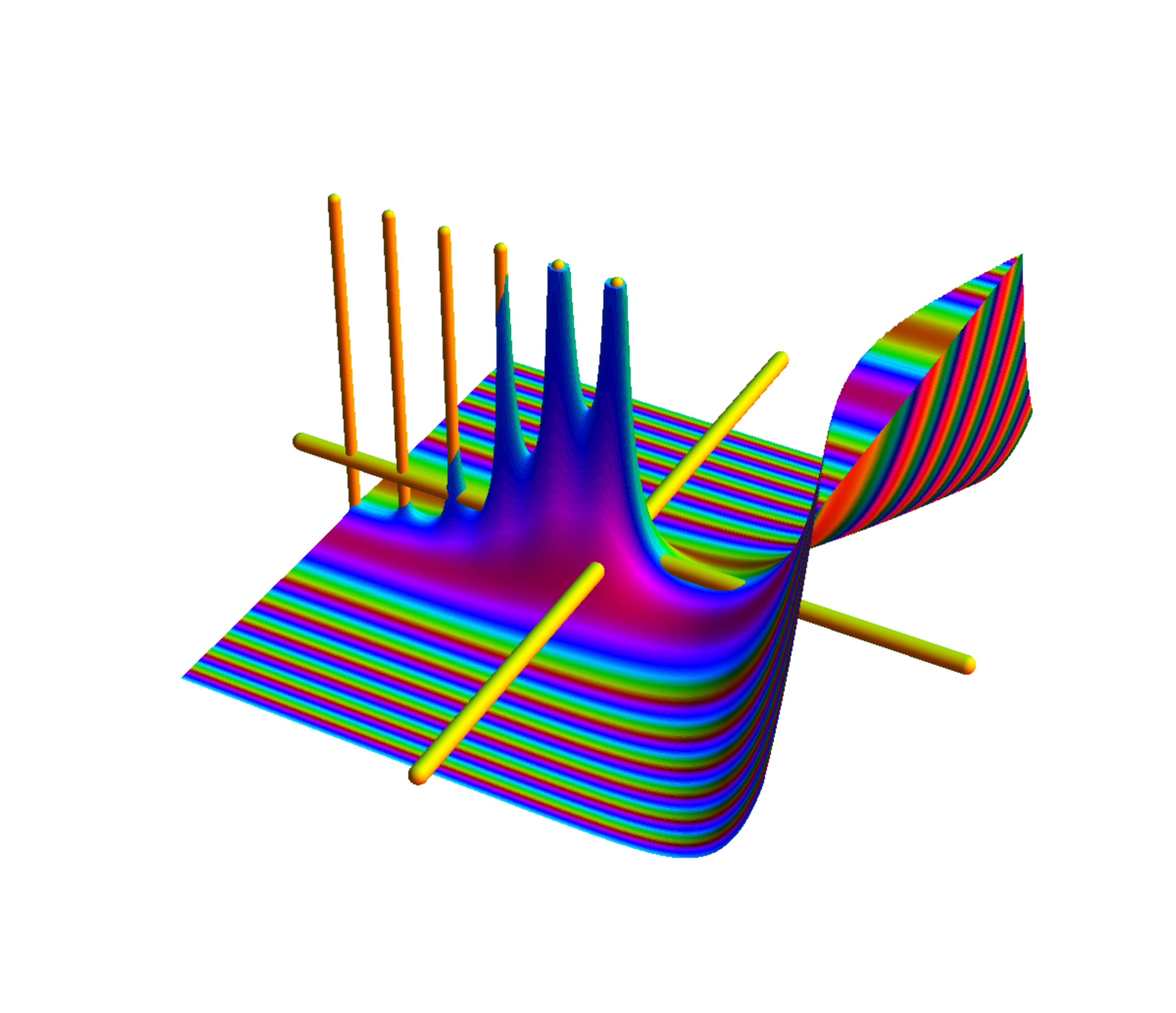}}
} \parbox{6cm}{
The graph of $1/|\zeta(x+i y)|$ shows the zeros of 
the zeta function $\zeta(z)$ as peaks. The Riemann conjecture is that all
these roots are on the line $x=1/2$. The right figure
shows the Gamma function which extends the factorial function
from positive integers to the complex plane $\Gamma(x)=(x-1)!$
for positive $x$. These graphs are produced in a way so that
they can be printed. 
}}

\parbox{16.8cm}{ \parbox{10cm}{
\scalebox{0.121}{\includegraphics{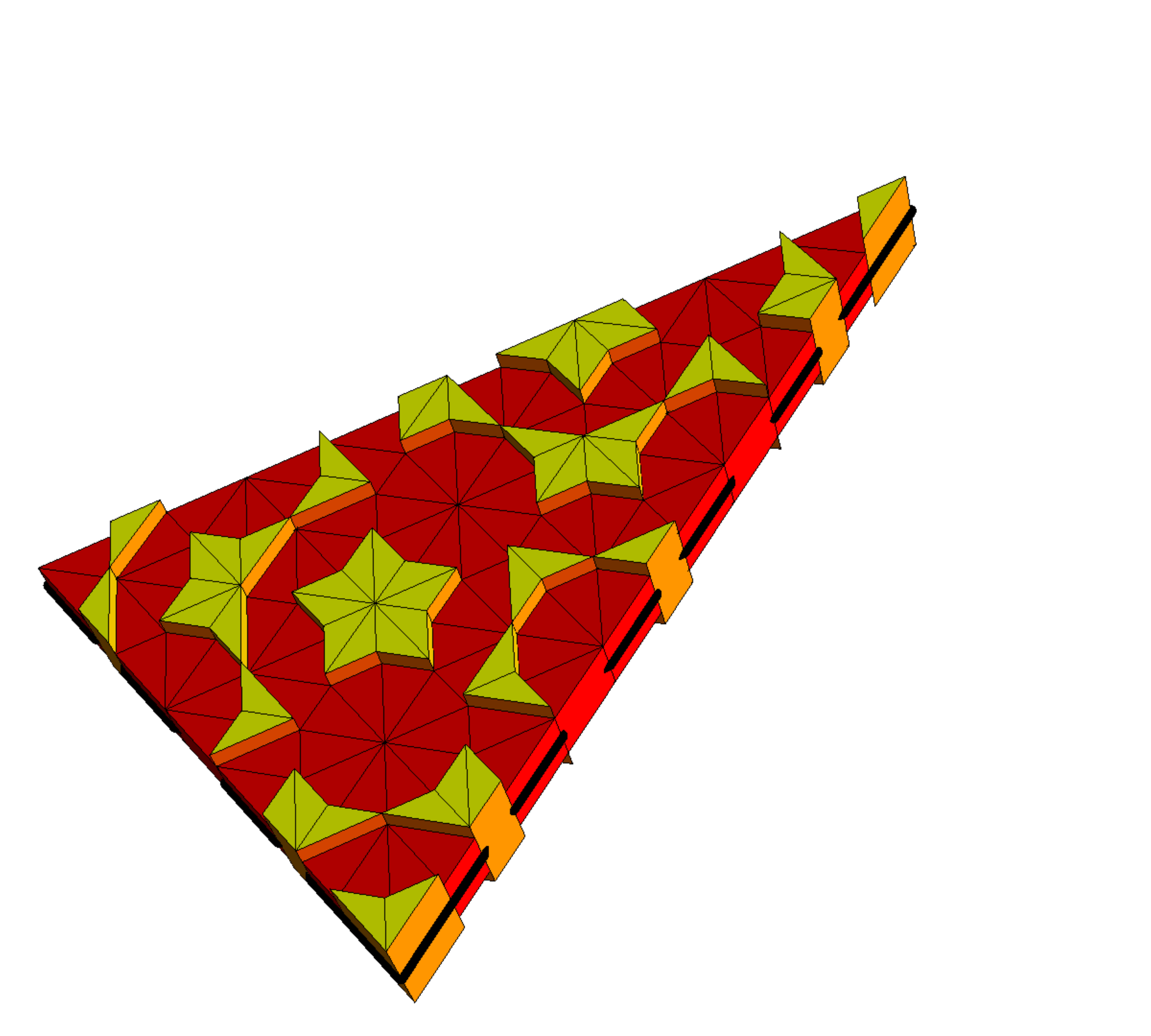}}
\scalebox{0.121}{\includegraphics{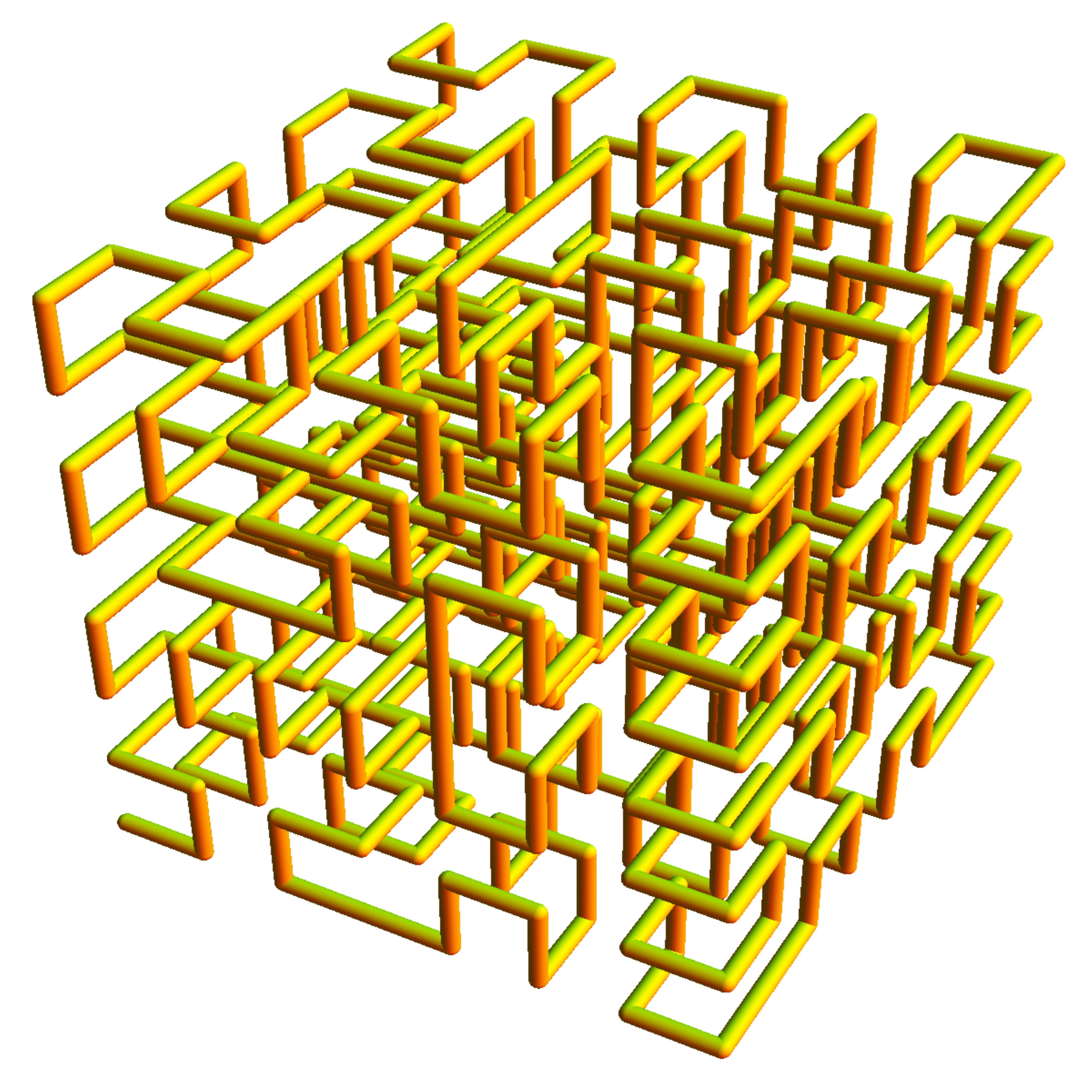}}
} \parbox{6cm}{
Two figures from different areas of geometry. The first picture 
allows printing an aperiodic Penrose tiling consisting of darts and kites.
To construct the tiling in 2D first, we used code from 
\cite{Wagon} section 10.2. The second figure is the third stage of the 
recursively defined Peano curve, a space filling curve. 
}}

\parbox{16.8cm}{ \parbox{10cm}{
\scalebox{0.121}{\includegraphics{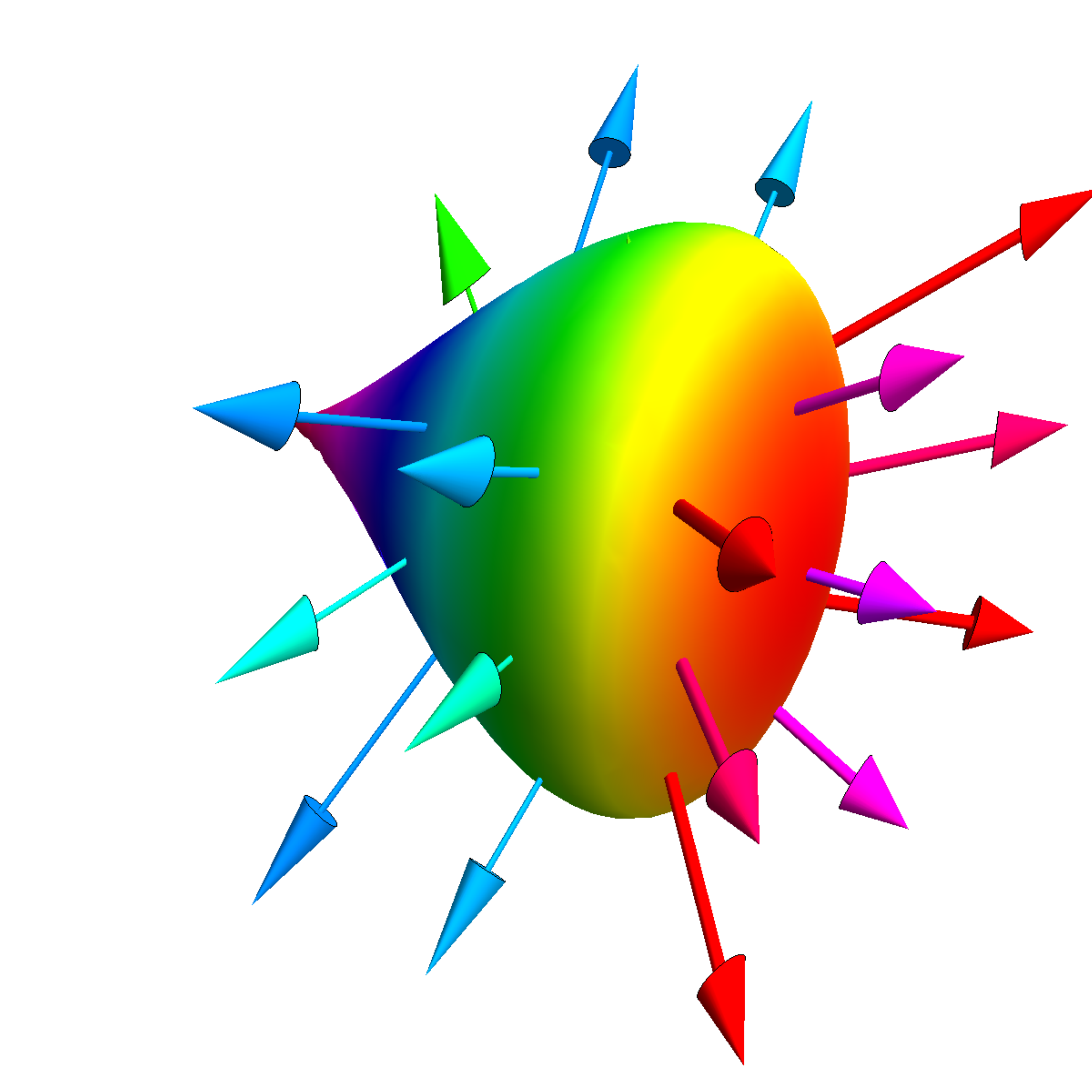}}
\scalebox{0.121}{\includegraphics{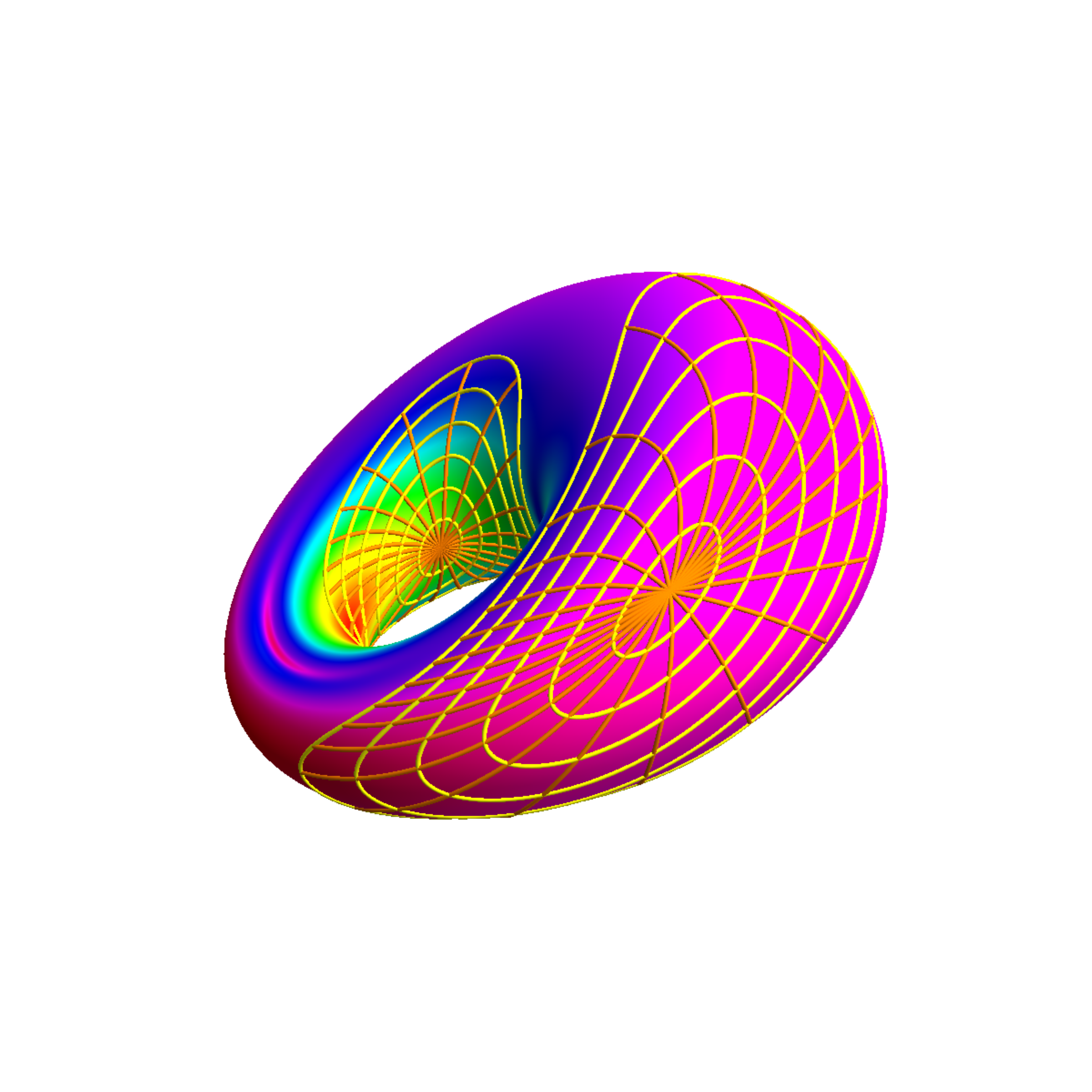}}
} \parbox{6cm}{
An illustration of the theorem in multivariable calculus that
the gradient is perpendicular to the level surface. 
The second picture illustrates the exponential map in Riemannian geometry,
where we see wave fronts at a point of positive curvature and at a point of
negative curvature. The differential equations are complicated but Mathematica
takes care of it. }}

\parbox{16.8cm}{ \parbox{10cm}{
\scalebox{0.121}{\includegraphics{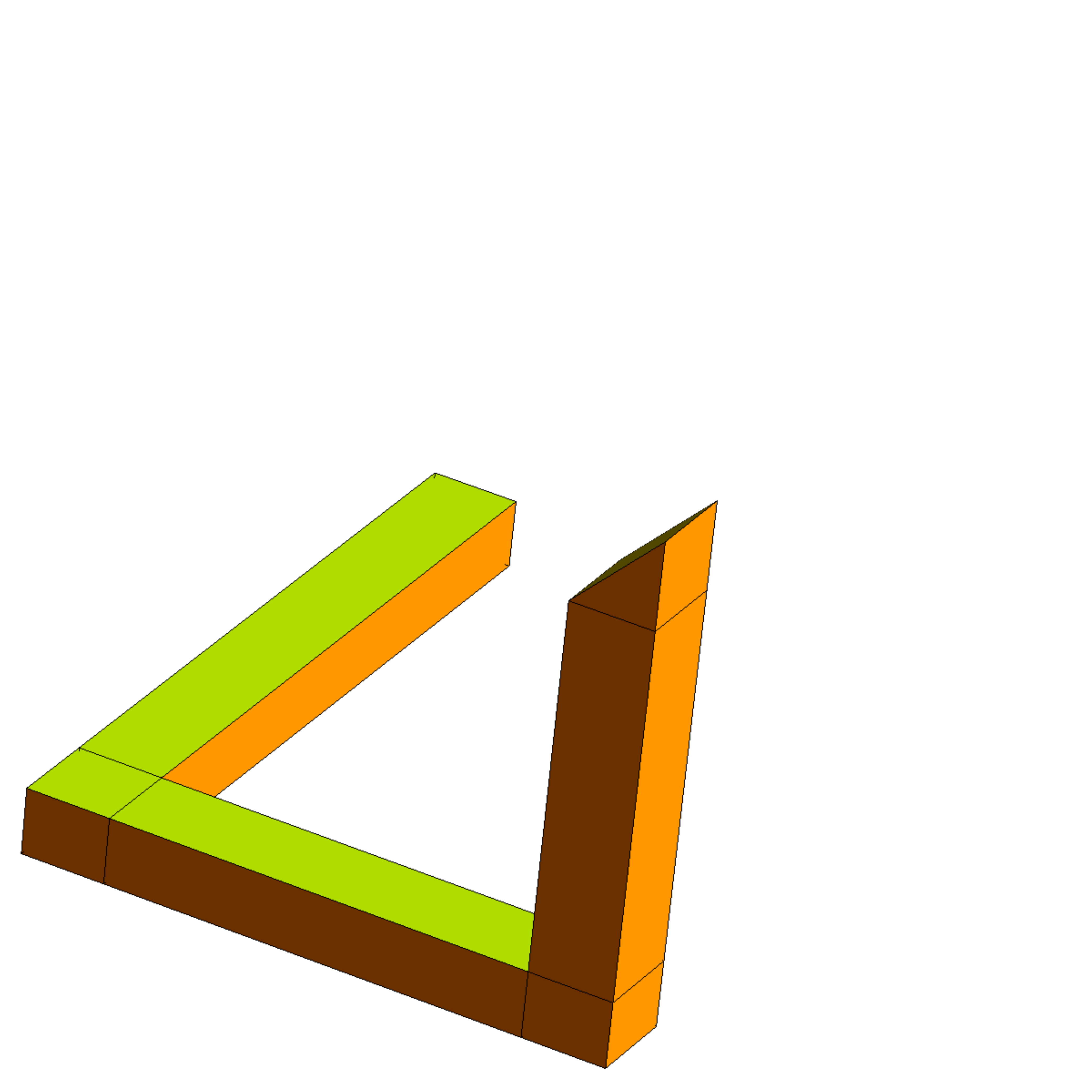}}
\scalebox{0.121}{\includegraphics{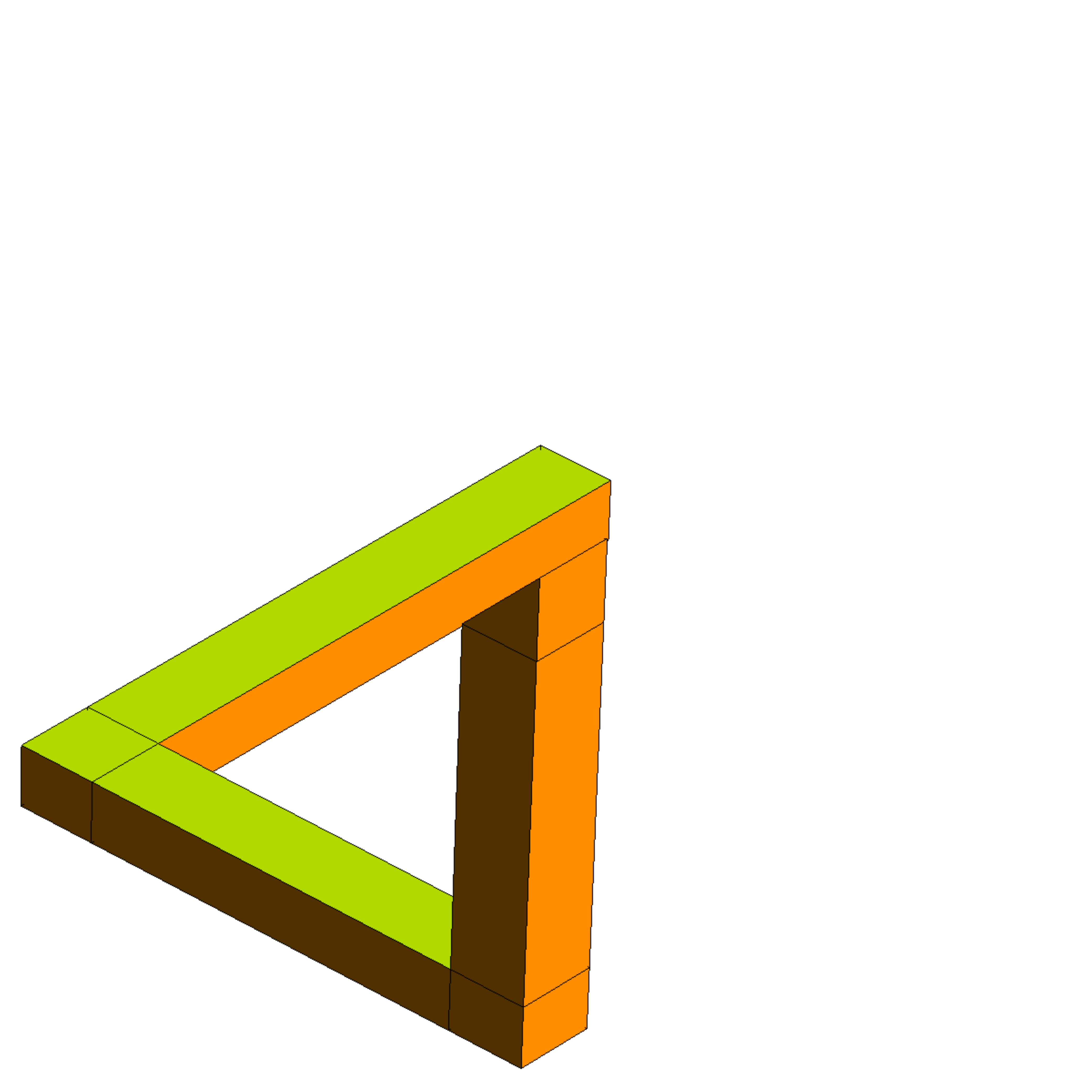}}
} \parbox{6cm}{
Printing the Penrose triangle. The solid was created by Oscar Reutersvard and
popularized by Roger Penrose \cite{Francis}. A Mathematica
implementation has first appeared in \cite{Trott}. The figure is featured on
one of the early editions of \cite{Hofstadter} and of \cite{Barrow}. 
}}

\parbox{16.8cm}{ \parbox{10cm}{
\scalebox{0.121}{\includegraphics{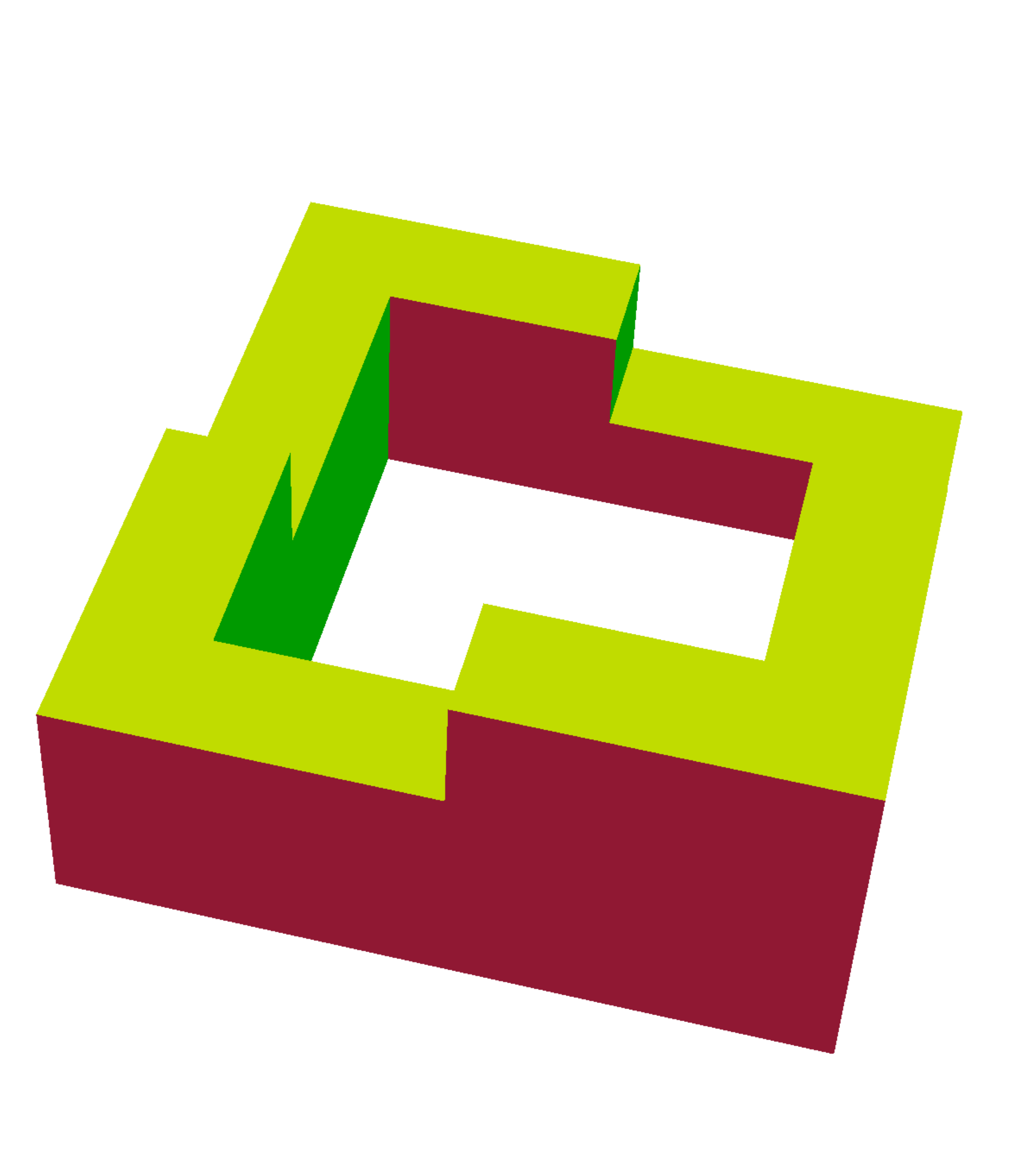}}
\scalebox{0.121}{\includegraphics{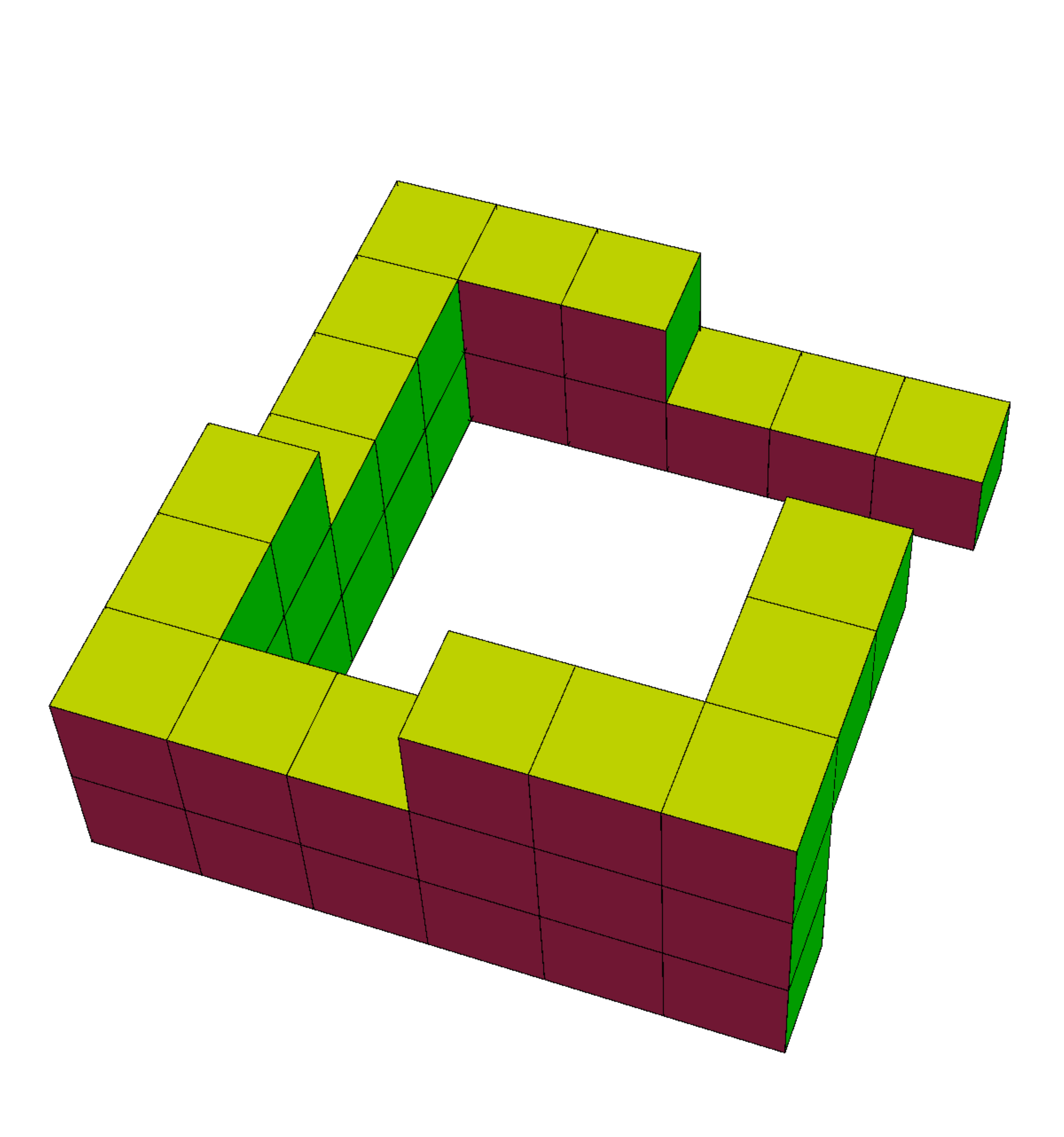}}
} \parbox{6cm}{
Printing a simplified version of the Escher stairs. 
If  the object is turned in the right angle, an impossible
stair is visible. When printed, this object can visualize
the geometry of impossible figures. 
}}

\parbox{16.8cm}{ \parbox{10cm}{
\scalebox{0.12}{\includegraphics{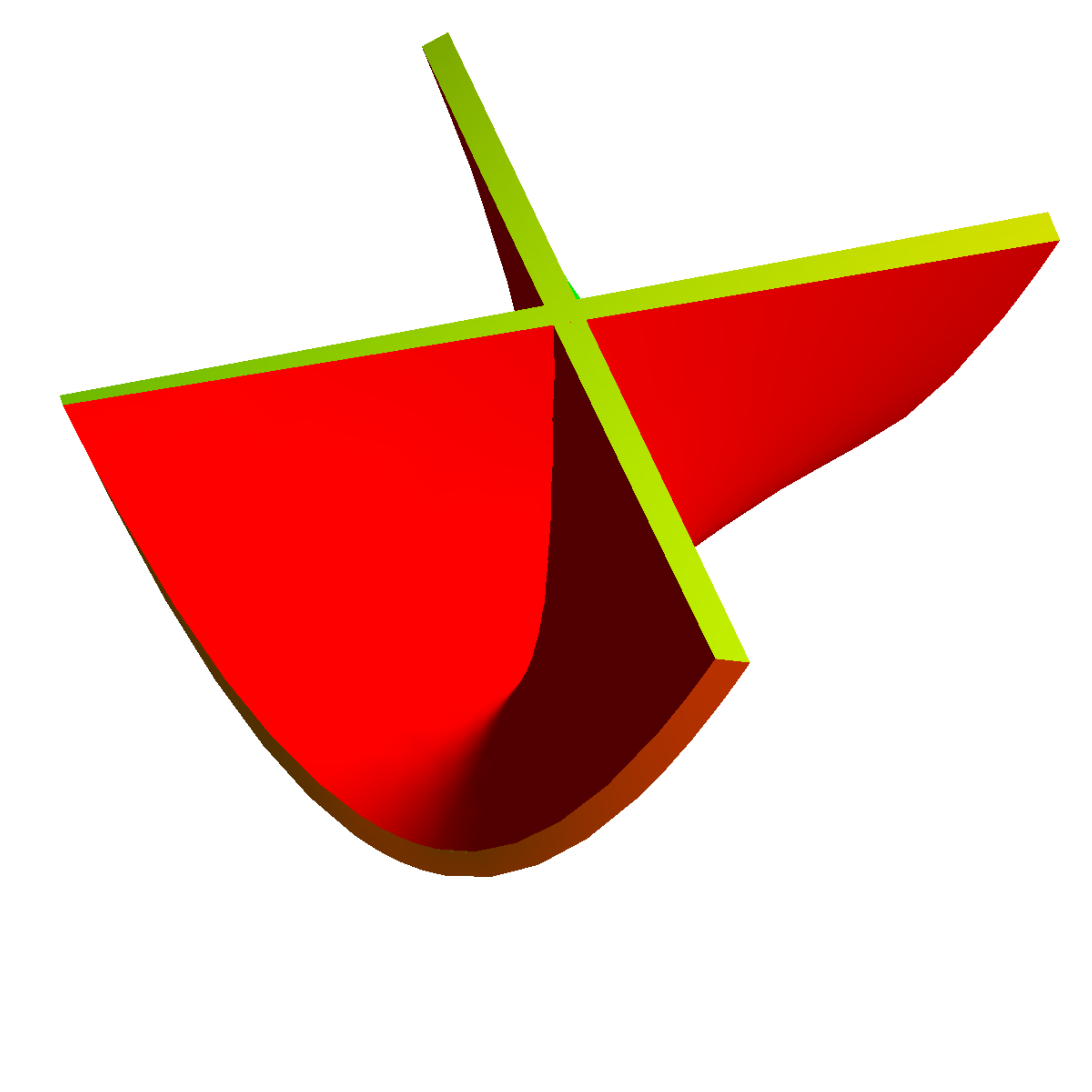}}
\scalebox{0.12}{\includegraphics{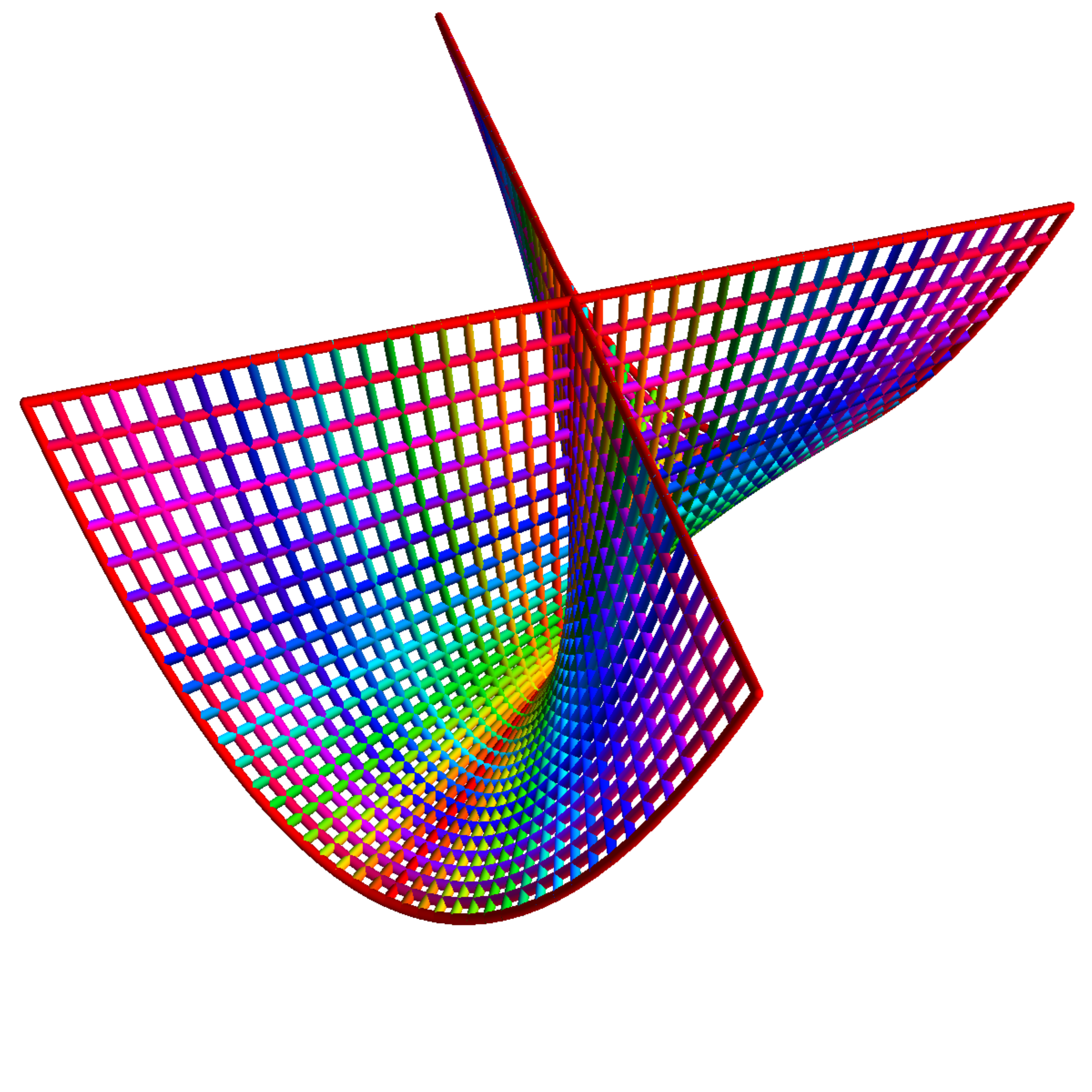}}
} \parbox{6cm}{
The Whitney umbrella is an icon of catastrophe theory. 
This is a typical shape of a caustic of a wave front moving
in space. To the left, we see how the surface has been thickened
to make it printable. On the right, the grid curves are shown 
as tubes. Also this is a technique which is printable.
}}

\parbox{16.8cm}{ \parbox{10cm}{
\scalebox{0.121}{\includegraphics{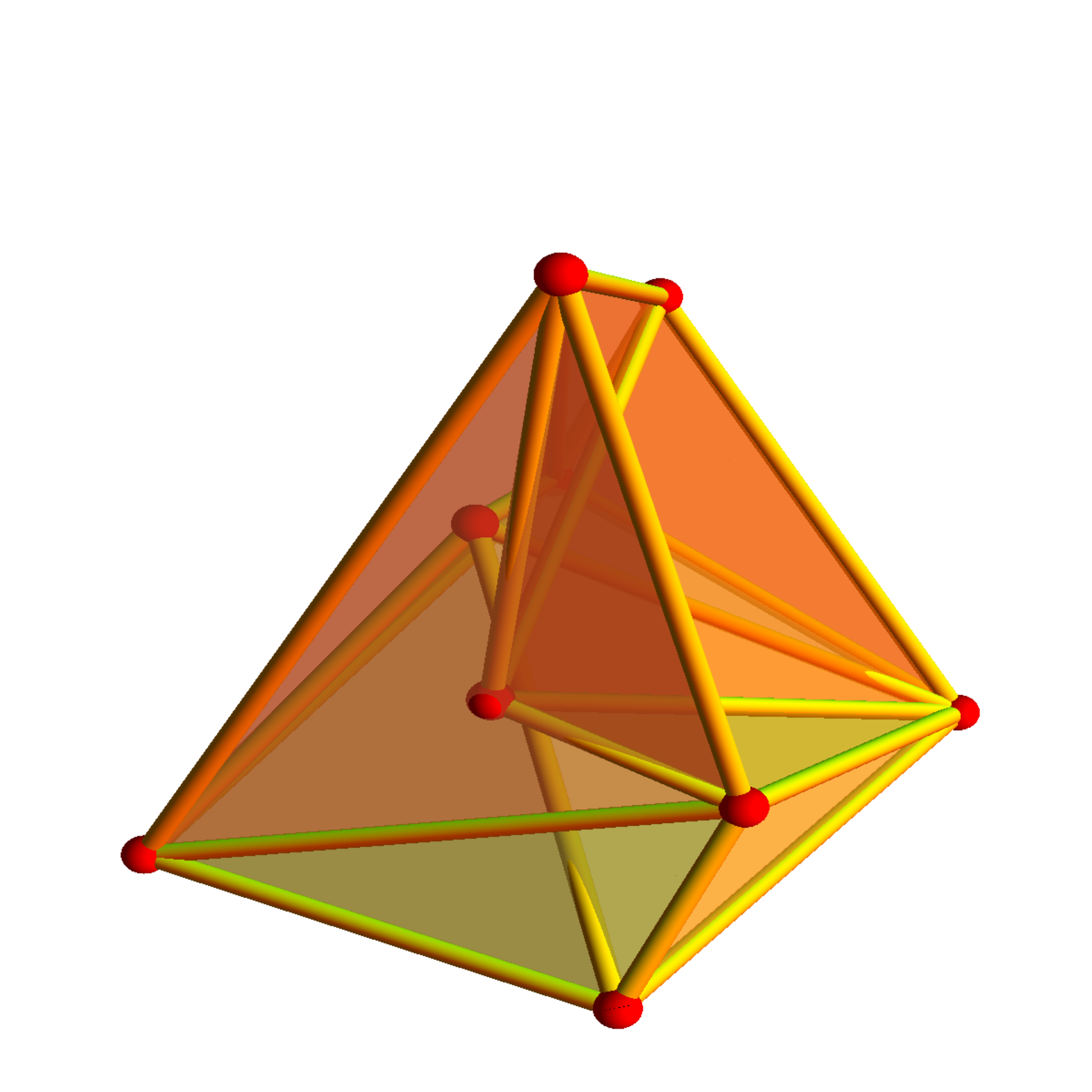}}
\scalebox{0.121}{\includegraphics{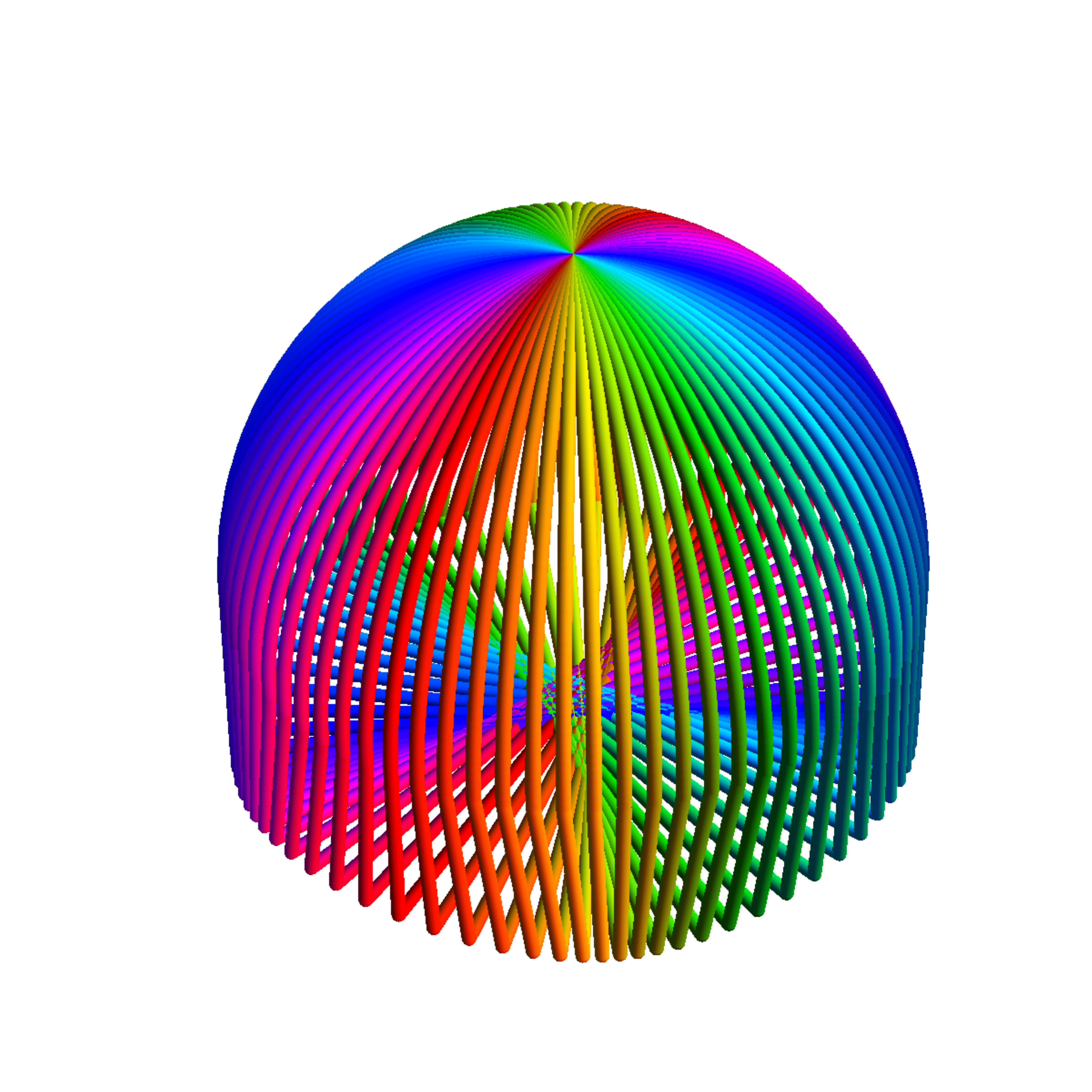}}
} \parbox{6cm}{
The left figure shows the Steffen polyhedron, a flexible
surface. It can be deformed without that the distances
between the points change. This is a surprise, since a
theorem of Cauchy tells that this is not possible for 
convex solids \cite{AigZie}. 
The right picture illustrates
how one can construct caustics on surfaces which have
prescribed shape. 
}}

\parbox{16.8cm}{ \parbox{10cm}{
\scalebox{0.121}{\includegraphics{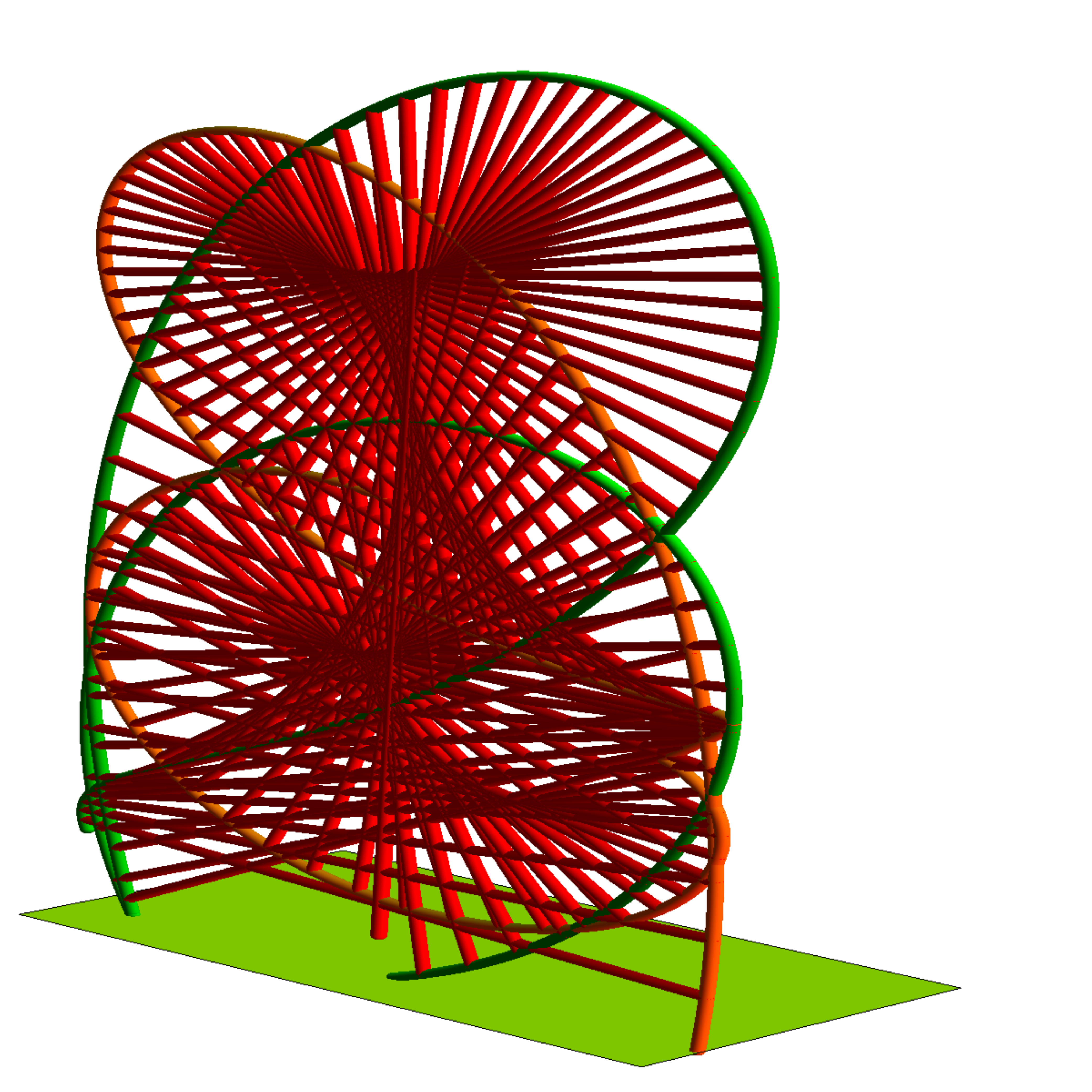}}
\scalebox{0.121}{\includegraphics{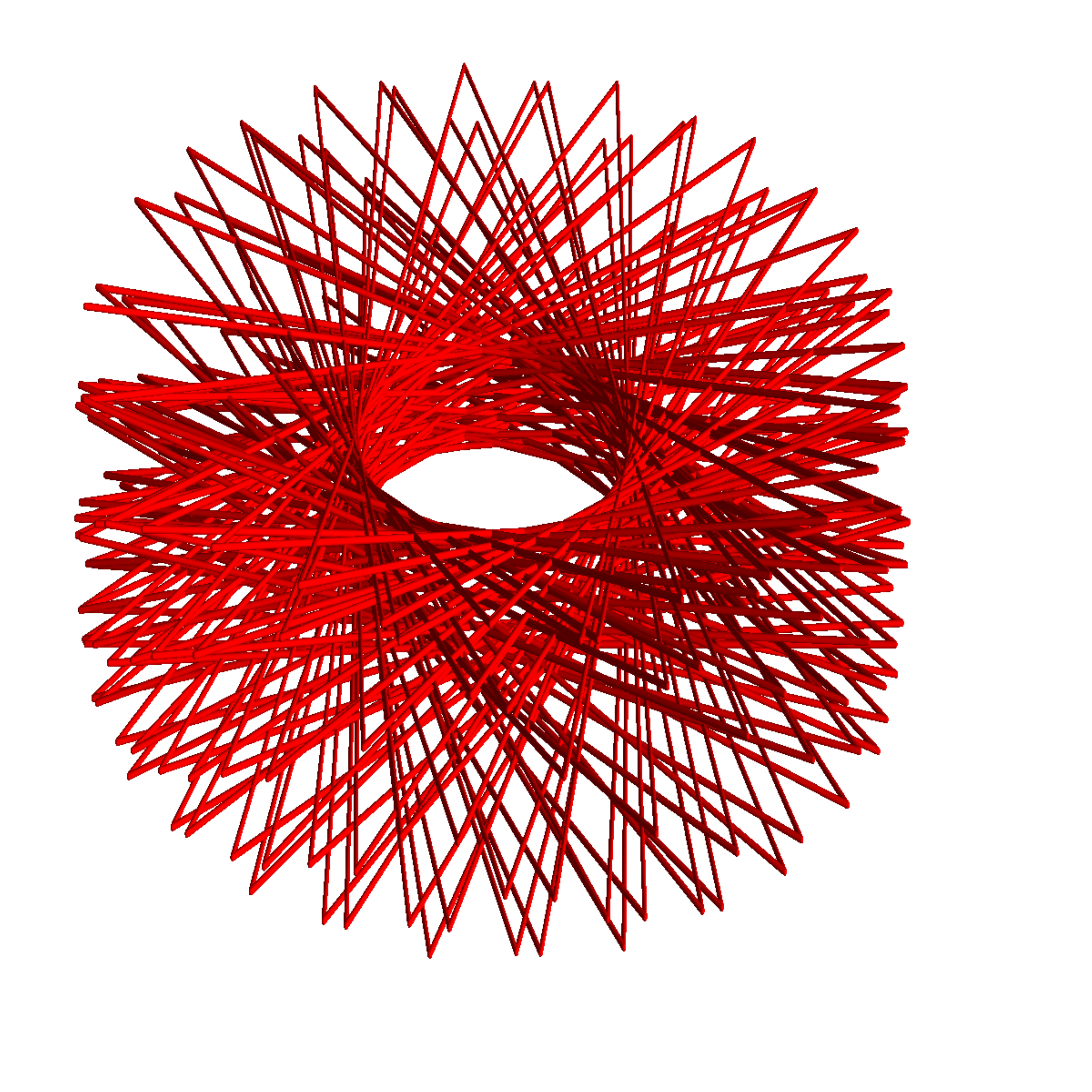}}
} \parbox{6cm}{
The first pictures illustrates a falling stick, bouncing
off a table. We see a stroboscopic snapshot of the trajectory.
The second picture illustrating the orbit of a billiard in a three dimensional 
billiard table \cite{elemente98}.
}}

\parbox{16.8cm}{ \parbox{10cm}{
\scalebox{0.121}{\includegraphics{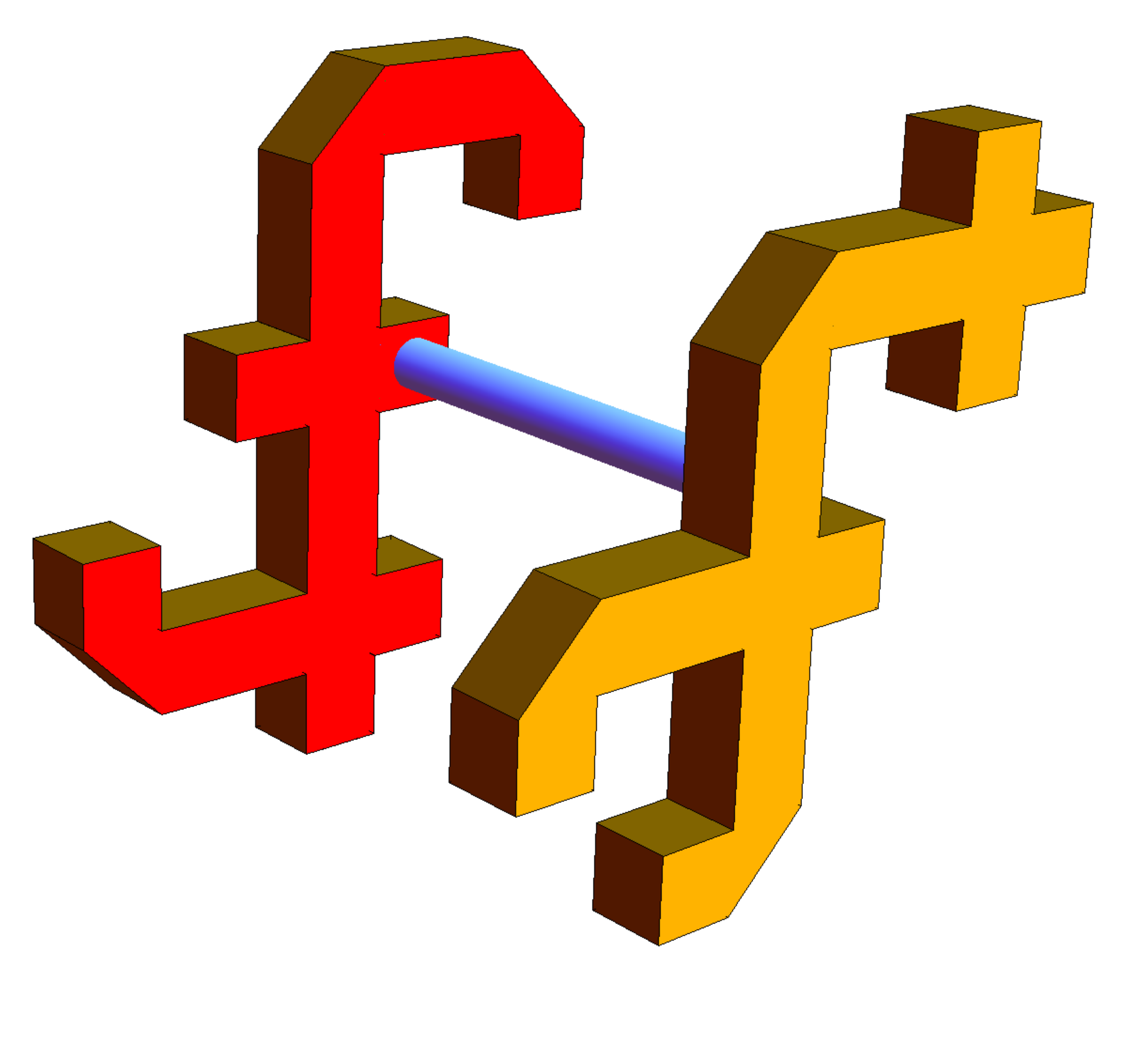}}
\scalebox{0.121}{\includegraphics{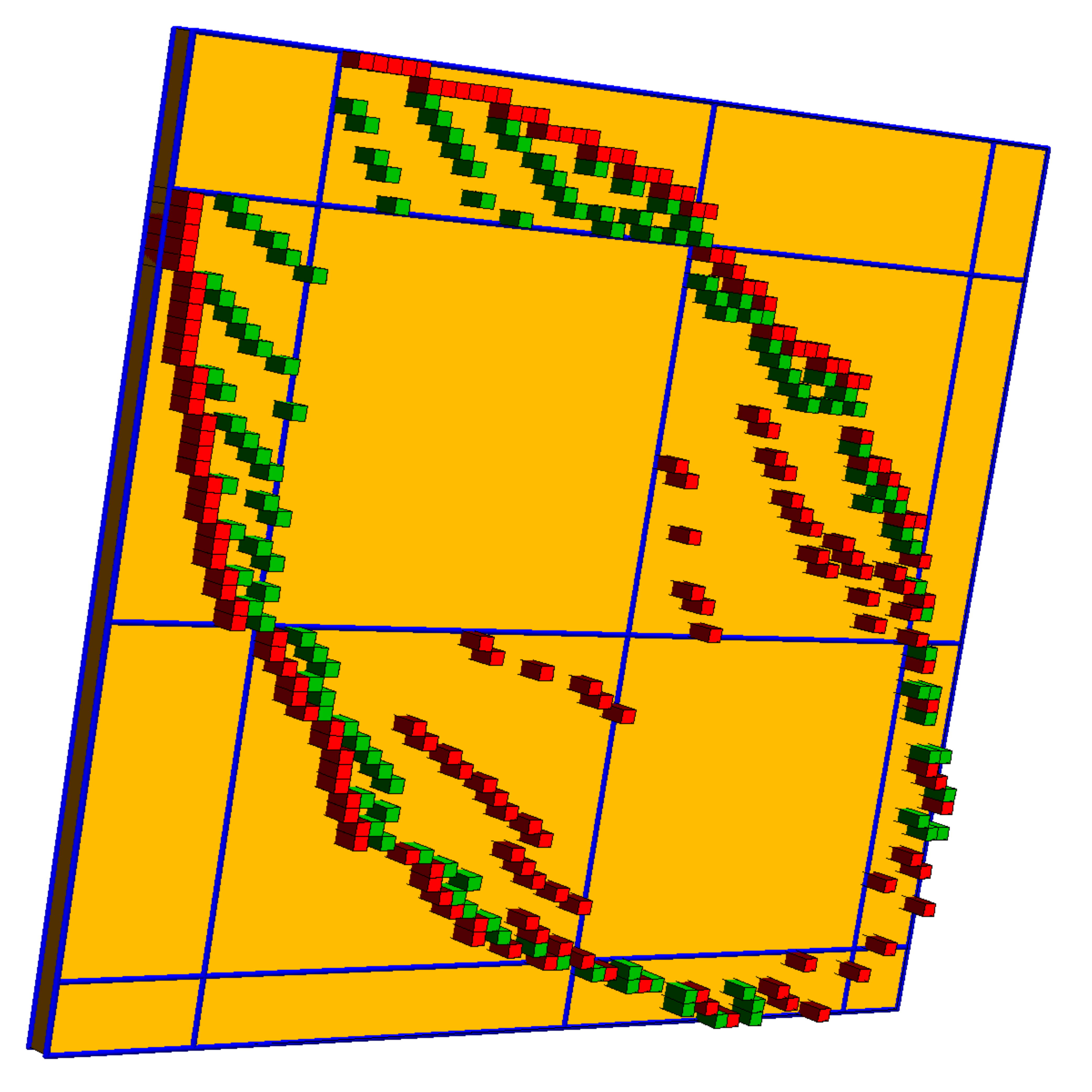}}
} \parbox{6cm}{
The left picture shows two isospectral drums found by Gordon-Webb. 
One does not know whether there are convex isospectral drums yet. 
The right picture shows a printed realization of a Dirac operator
of a graph \cite{knillmckeansinger3D}.  }}

\parbox{16.8cm}{ \parbox{10cm}{
\scalebox{0.121}{\includegraphics{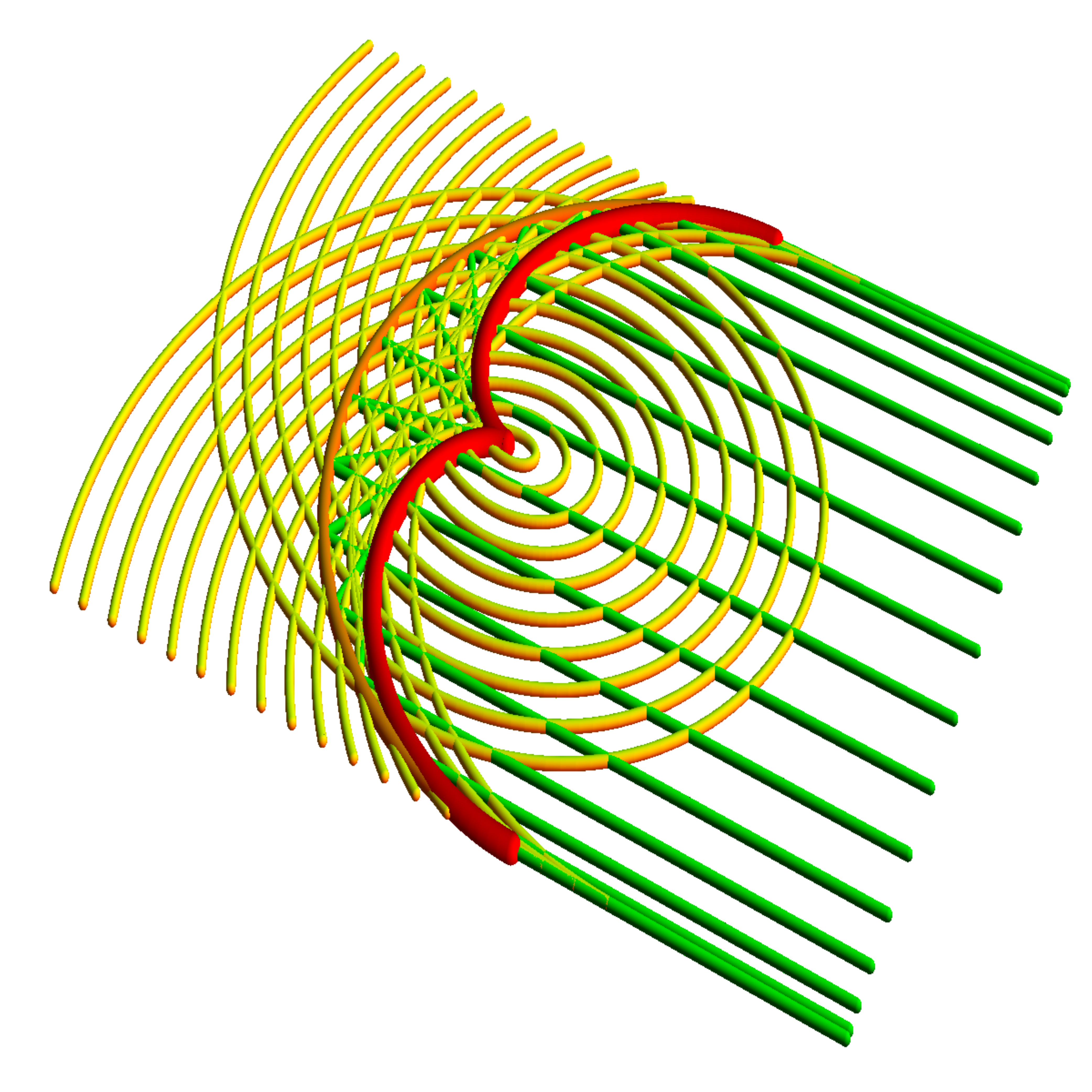}}
\scalebox{0.121}{\includegraphics{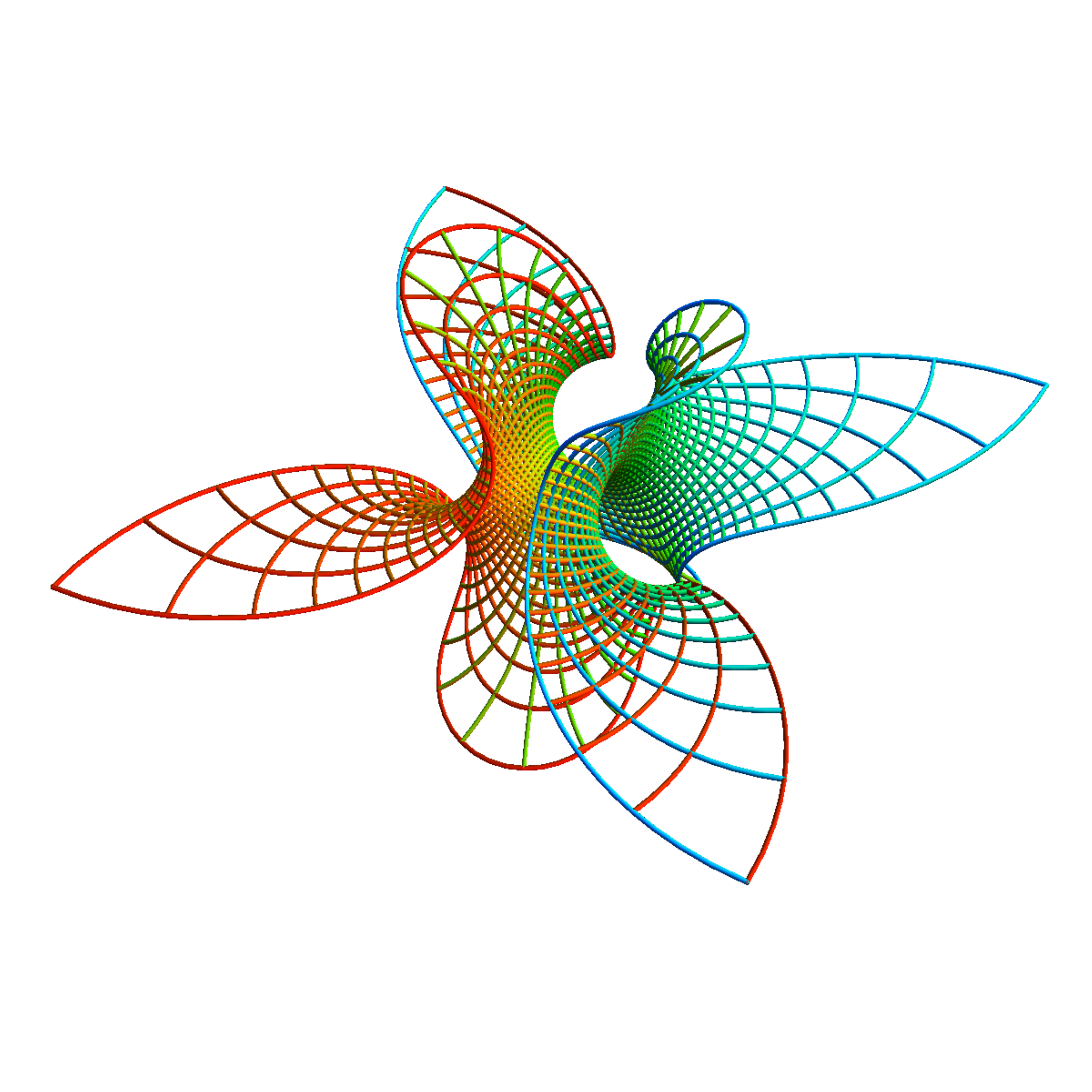}}
} \parbox{6cm}{
The left picture shows the coffee cup caustic. It is an icon of 
catastroph theory. 
The right picture
shows the Costa minimal surface using a parametrization 
found by Gray \cite{GrayMathematica}.
}}

\parbox{16.8cm}{ \parbox{10cm}{
\scalebox{0.121}{\includegraphics{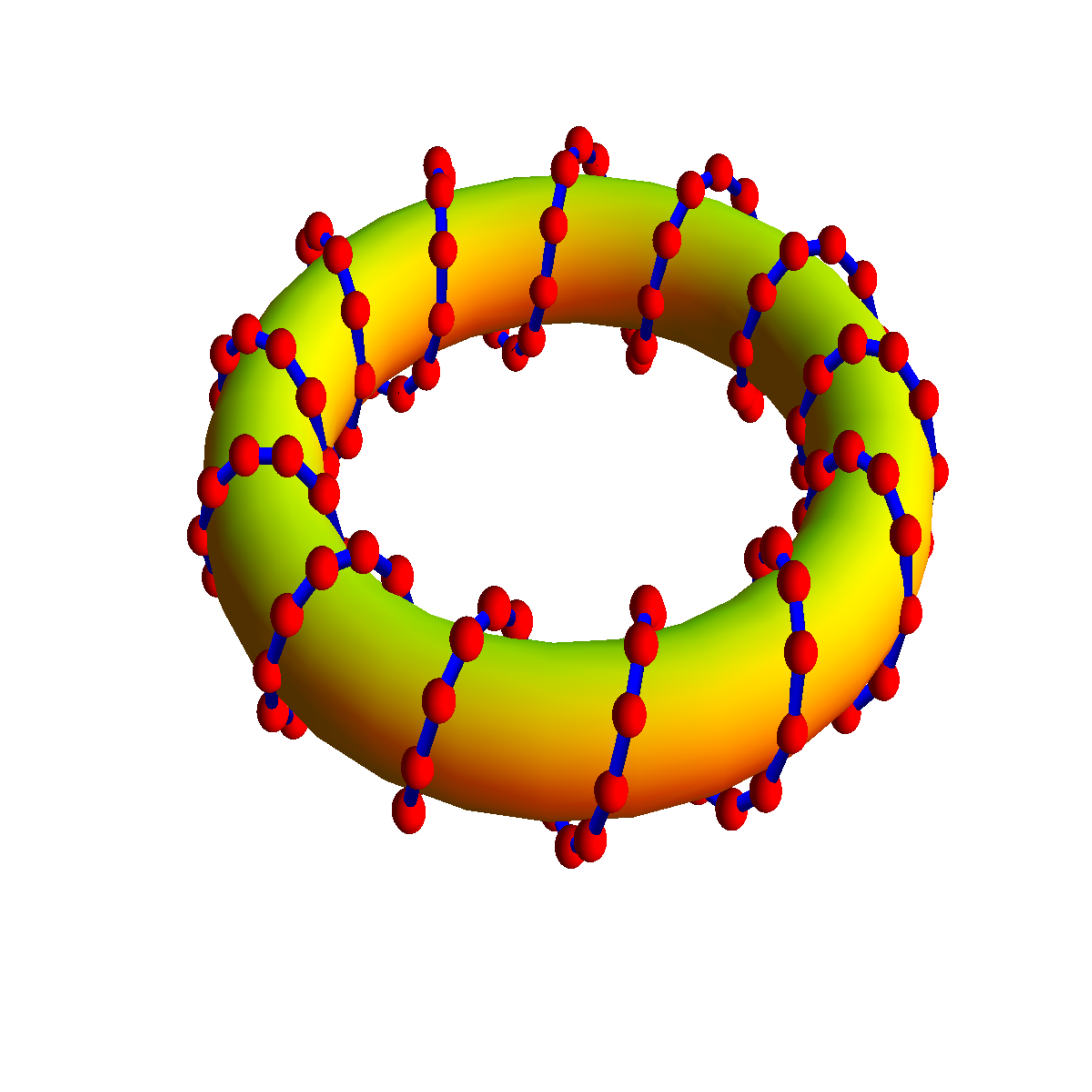}}
\scalebox{0.121}{\includegraphics{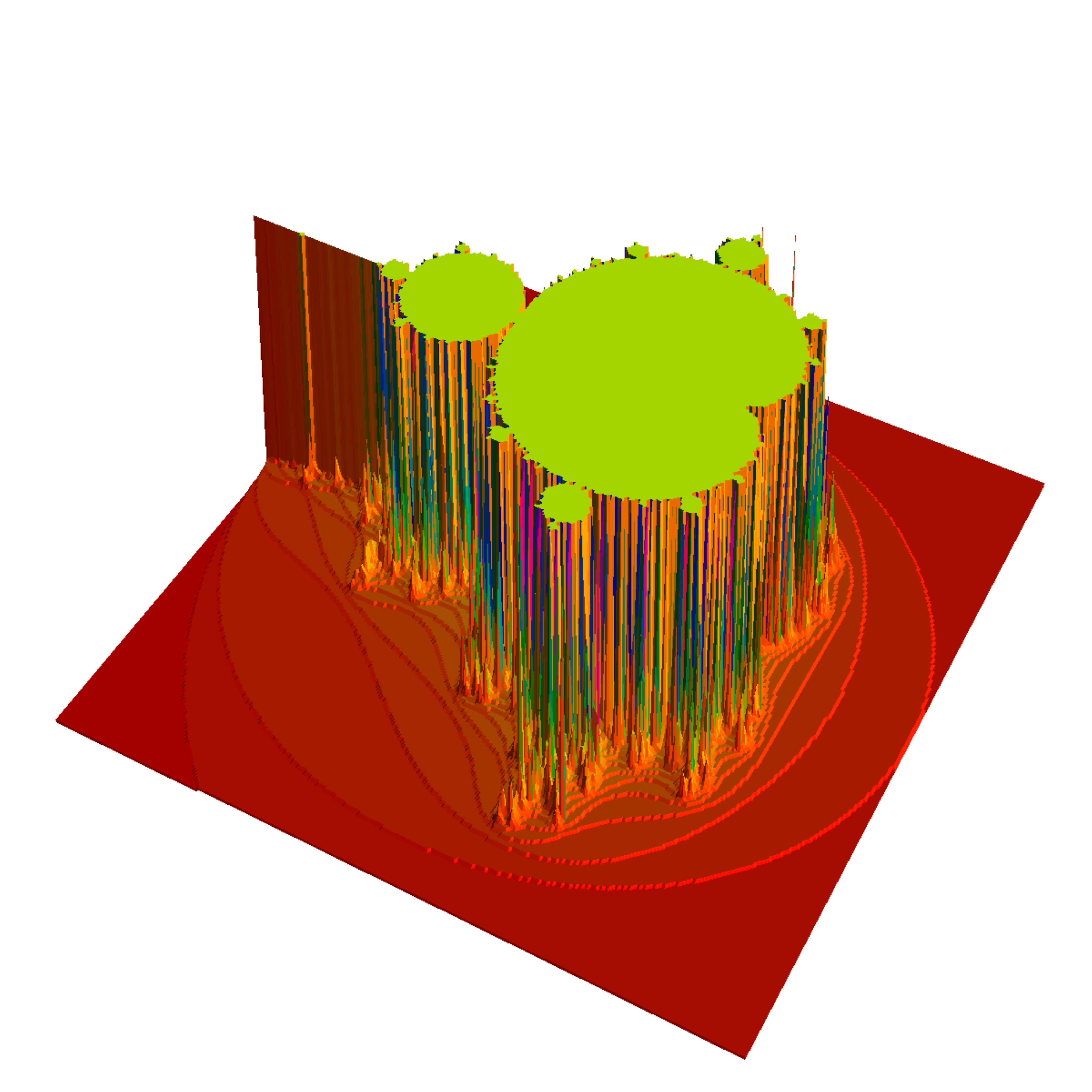}}
} \parbox{6cm}{
The left picture illustrates a torus graph
the right picture shows the Mandelbrot set in 3D.
Fantastic computer generated pictures of the fractal landscape
have been produced already 25 years ago \cite{PeitgenSaupe}. 
}}

\parbox{16.8cm}{ \parbox{10cm}{
\scalebox{0.121}{\includegraphics{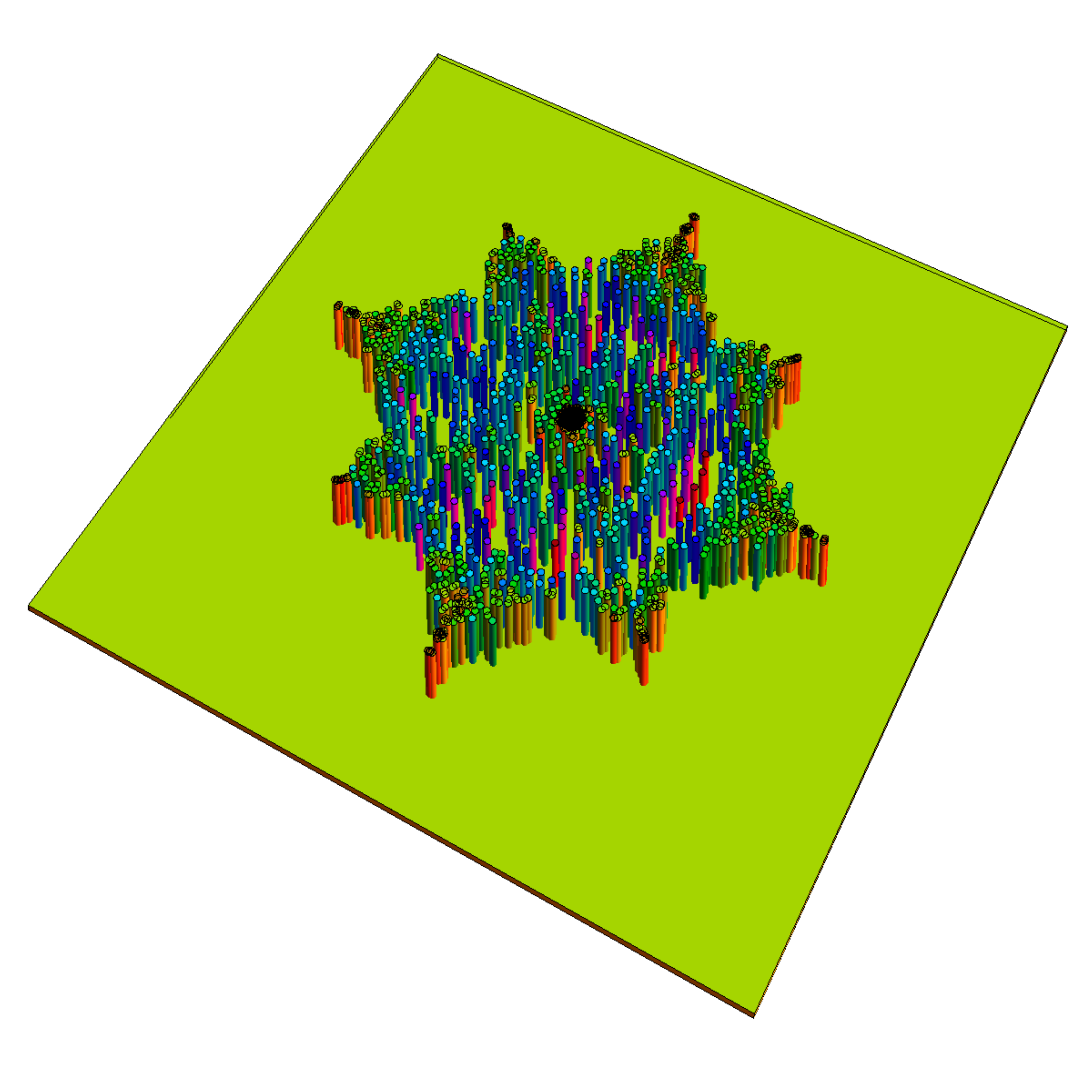}}
\scalebox{0.121}{\includegraphics{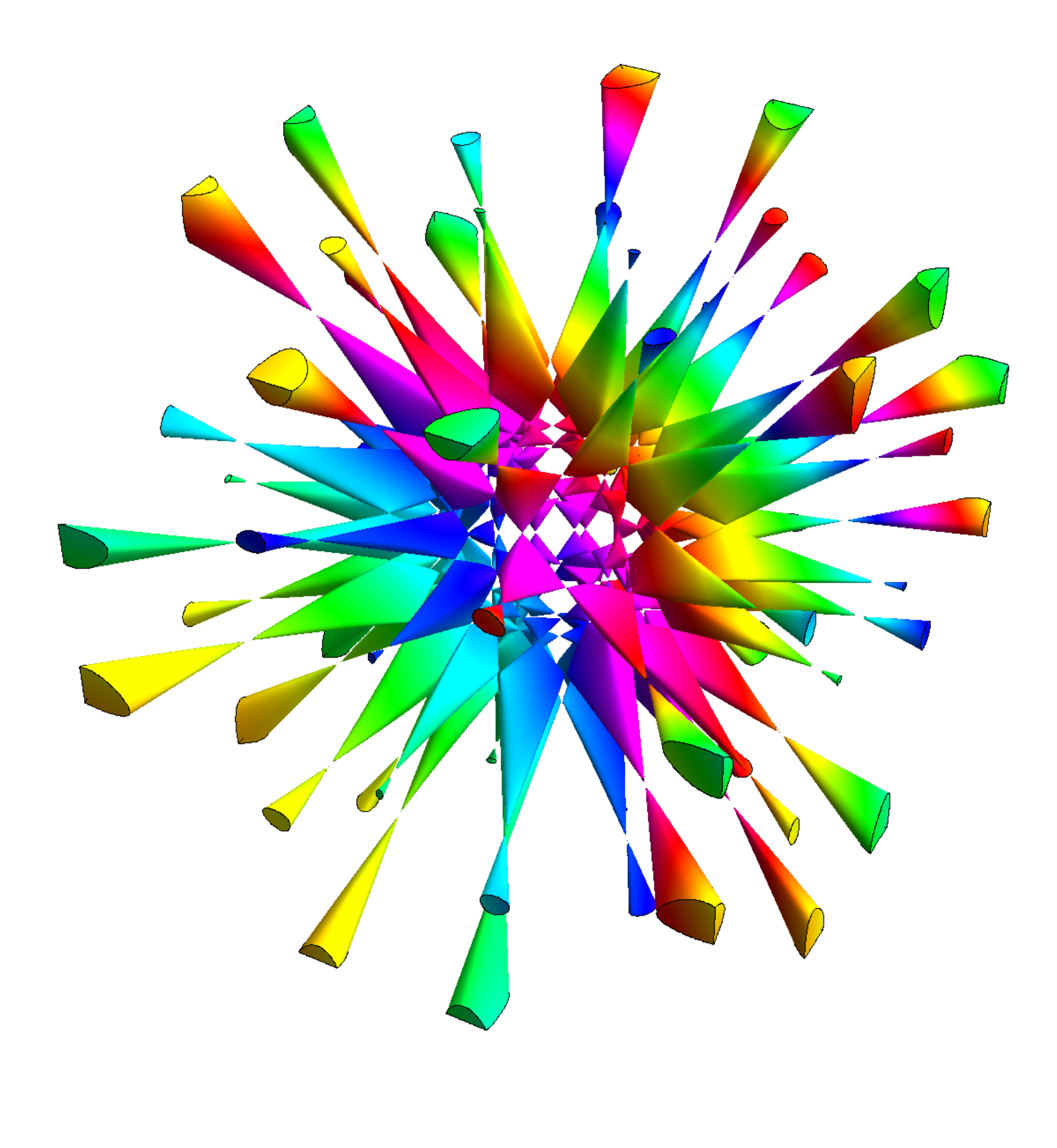}}
} \parbox{6cm}{
The left picture illustrates the spectrum of a matrix, where
the entries are random but correlated. The entries are given 
by the values of an almost periodic function. We have observed
experimentally that the spectrum is of fractal nature
in the complex plane. The picture could
be printed. The right example is a decic surface, the zero locus
$f(x,y,z)=0$ of a polynomial of degree $10$ in three variables .
We show the region $f(x,y,z) \leq 0$. 
}}

\parbox{16.8cm}{ \parbox{10cm}{
\scalebox{0.121}{\includegraphics{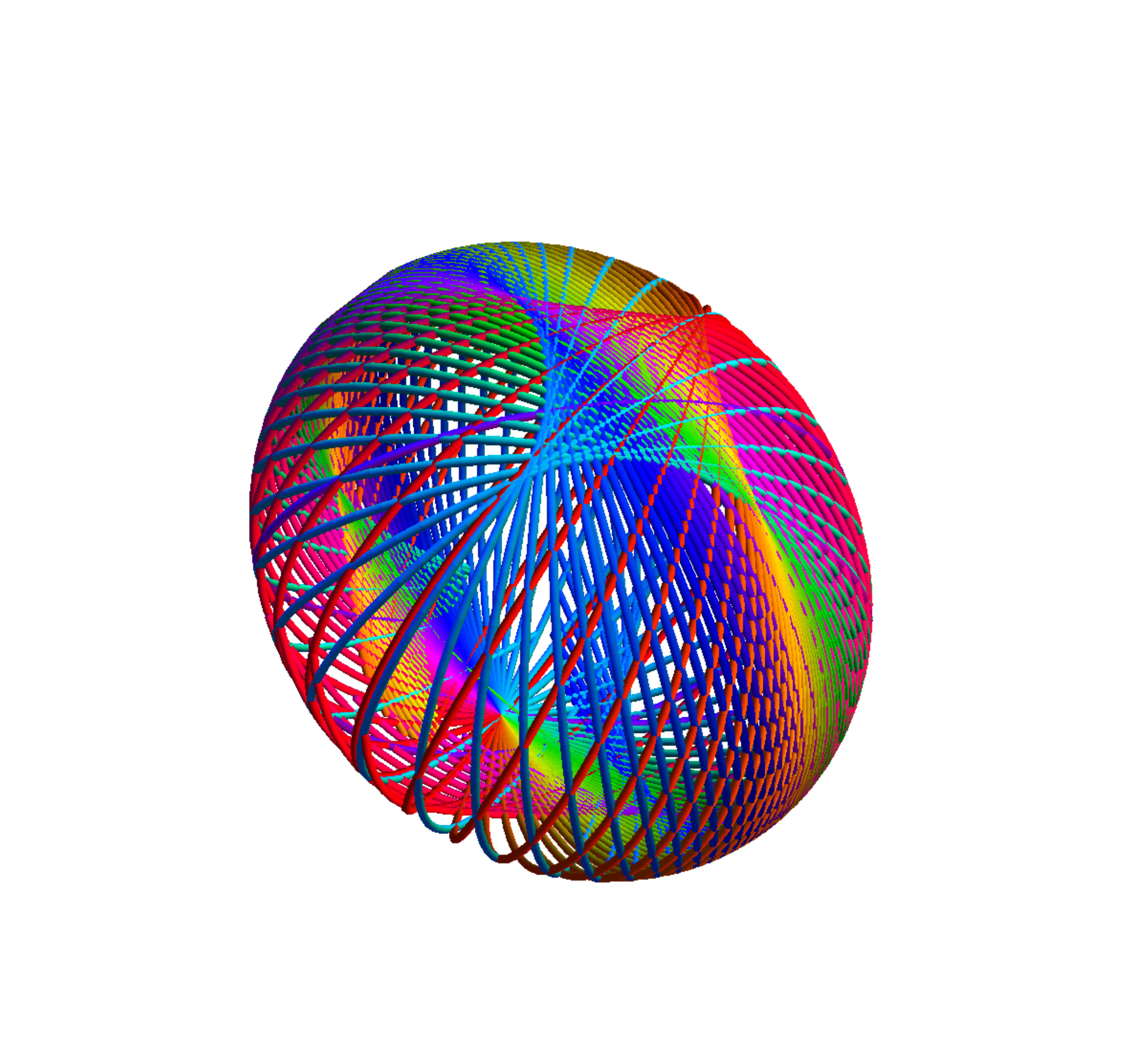}}
\scalebox{0.121}{\includegraphics{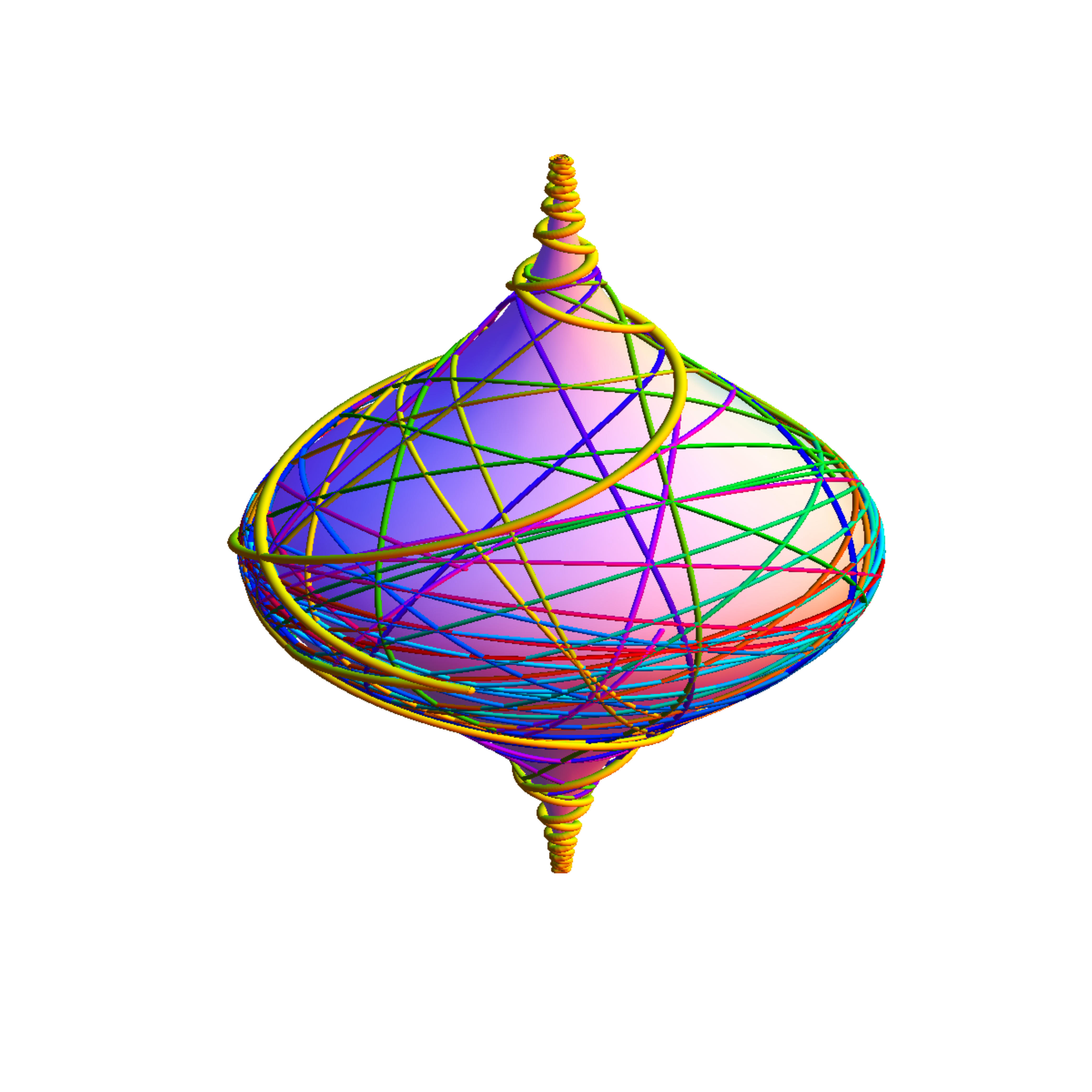}}
} \parbox{6cm}{
The left figure shows the geodesic flow on an ellipsoid 
without rotational symmetry.
Jacobi's last theorem - still an open problem - claims
that all caustics have 4 cusps. The right picture shows some
geodesics starting at a point of a surface of revolution. 
}}

\parbox{16.8cm}{ \parbox{10cm}{
\scalebox{0.121}{\includegraphics{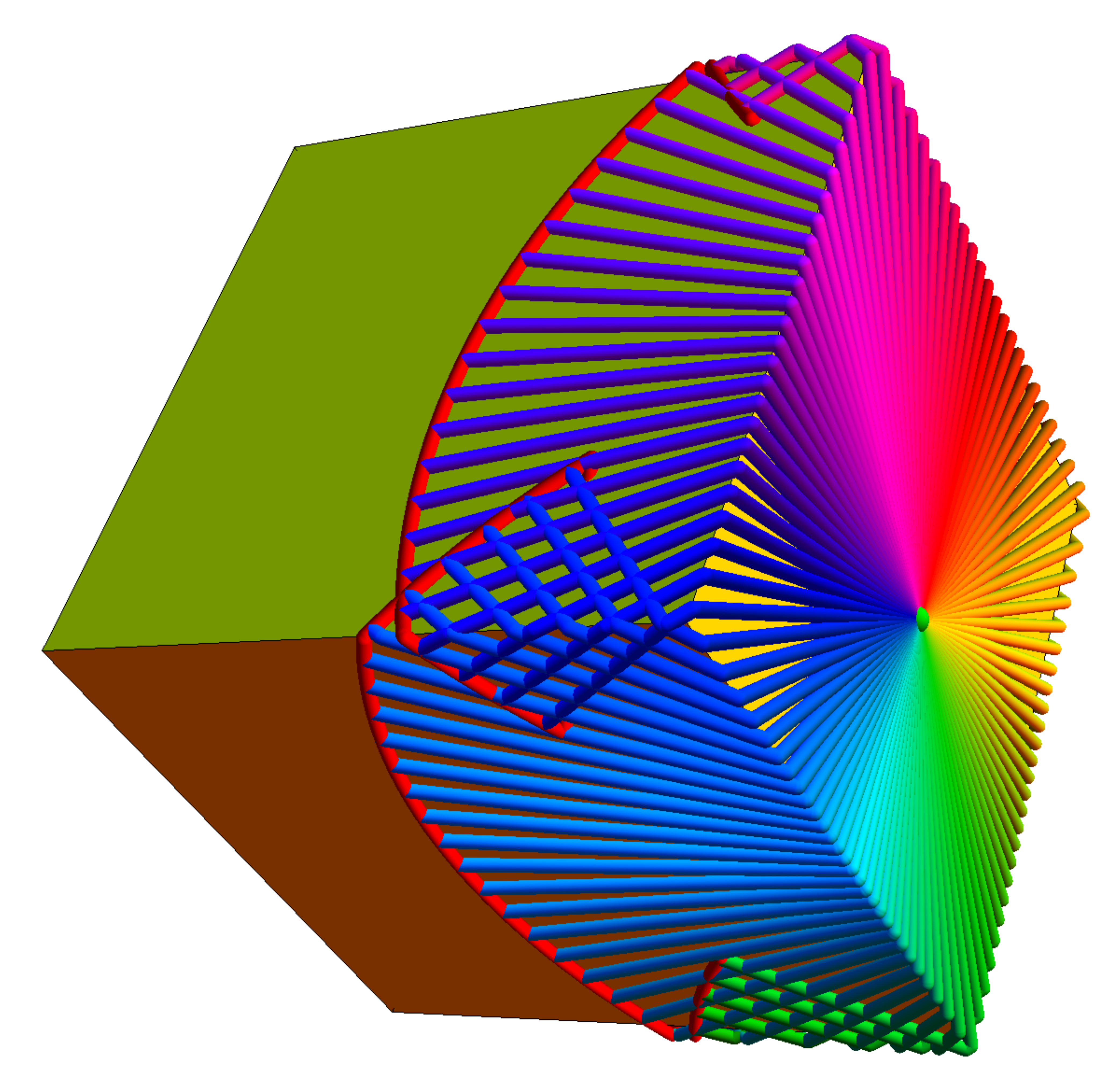}}
\scalebox{0.119}{\includegraphics{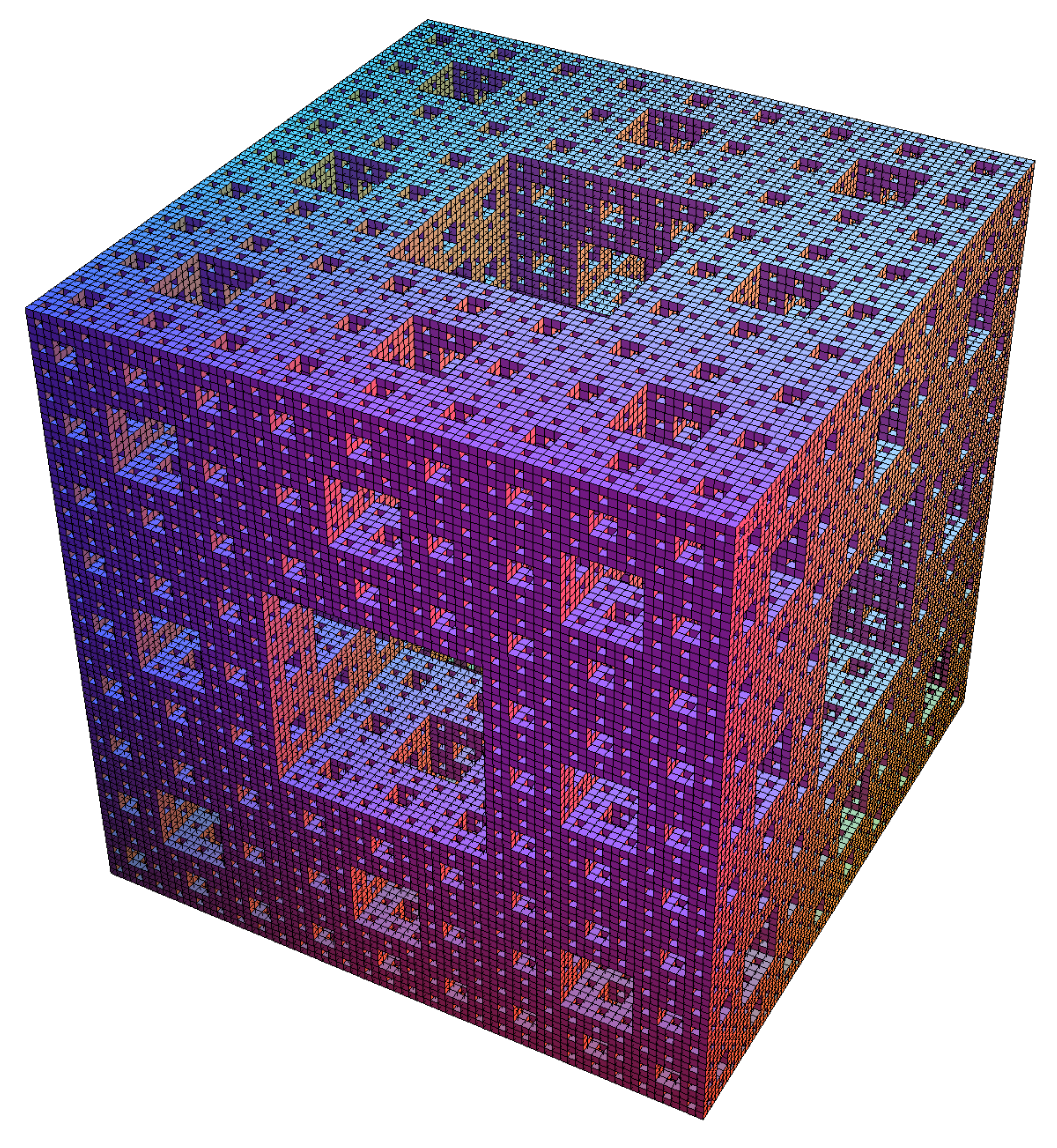}}
} \parbox{6cm}{
A wave front on a cube. Despite the simplicity of the setup, the 
wave fronts become very complicated. 
The right figure shows an approximation to the Menger sponge, 
a fractal in three dimensional space. It is important in topology
because it contains every compact metric space of topological dimension 
$1$.
}}

\parbox{16.8cm}{ \parbox{10cm}{
\scalebox{0.121}{\includegraphics{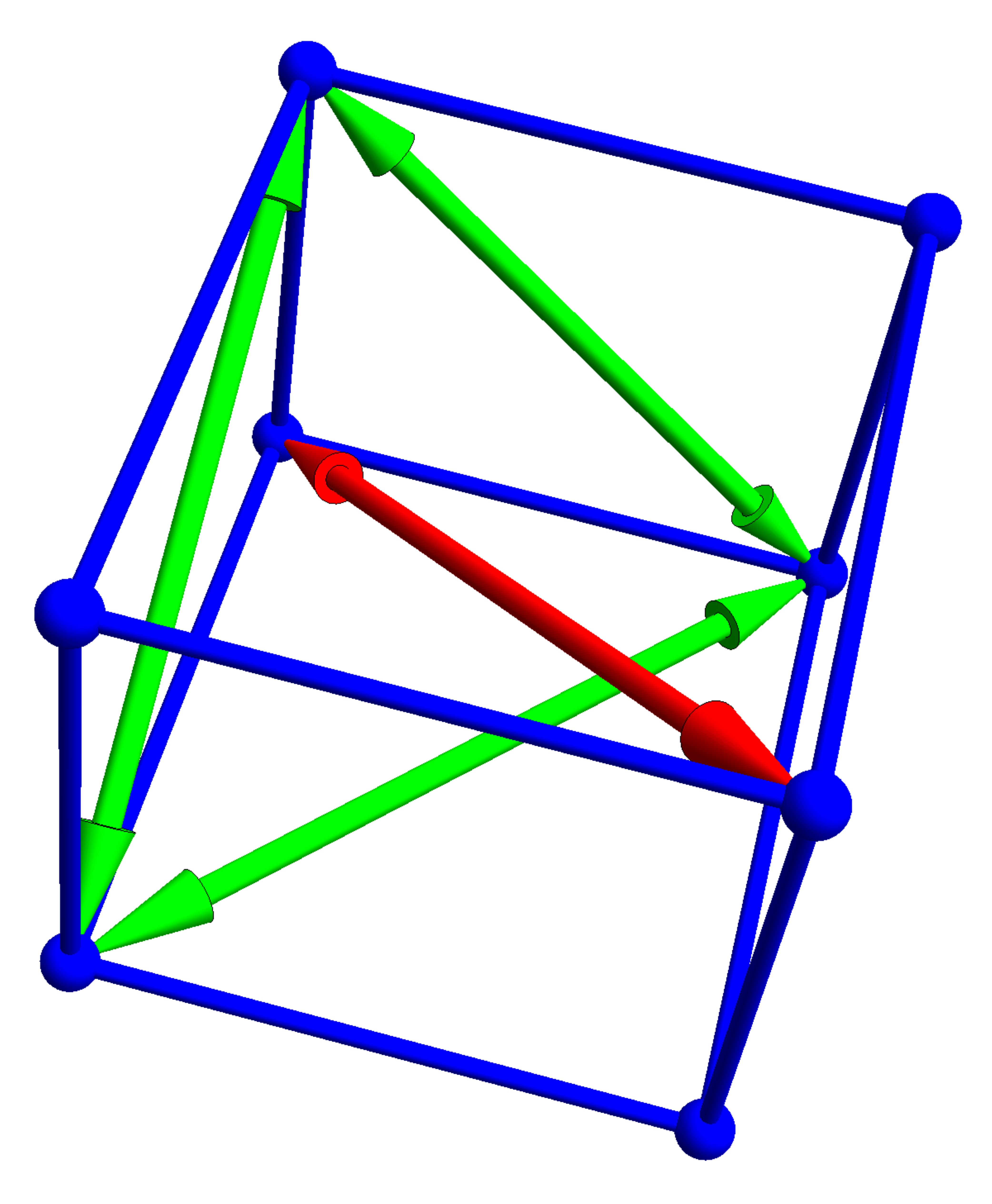}}
\scalebox{0.121}{\includegraphics{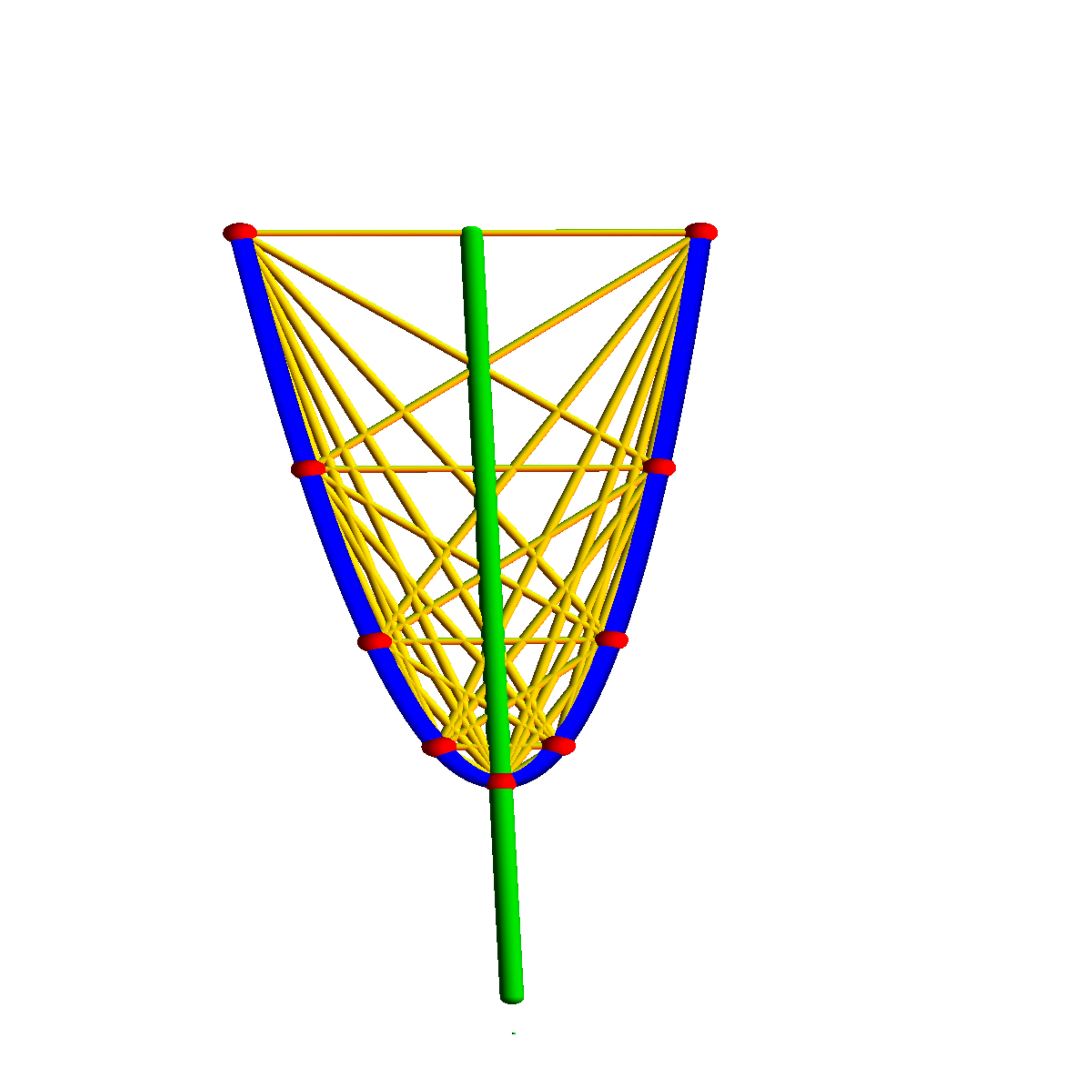}}
} \parbox{6cm}{
The left figure illustrates an Euler brick. It is unknown whether 
a cuboid for which all side lengths are integers and for which also all
face and space diagonals are integers. If all face diagonals have integer
length, it is called an Euler brick. If also the space diagonal is an integer
it is a perfect Euler brick. The right figure shows how one
can realize the multiplication of numbers using parabola.
}}

\parbox{16.8cm}{ \parbox{10cm}{
\scalebox{0.121}{\includegraphics{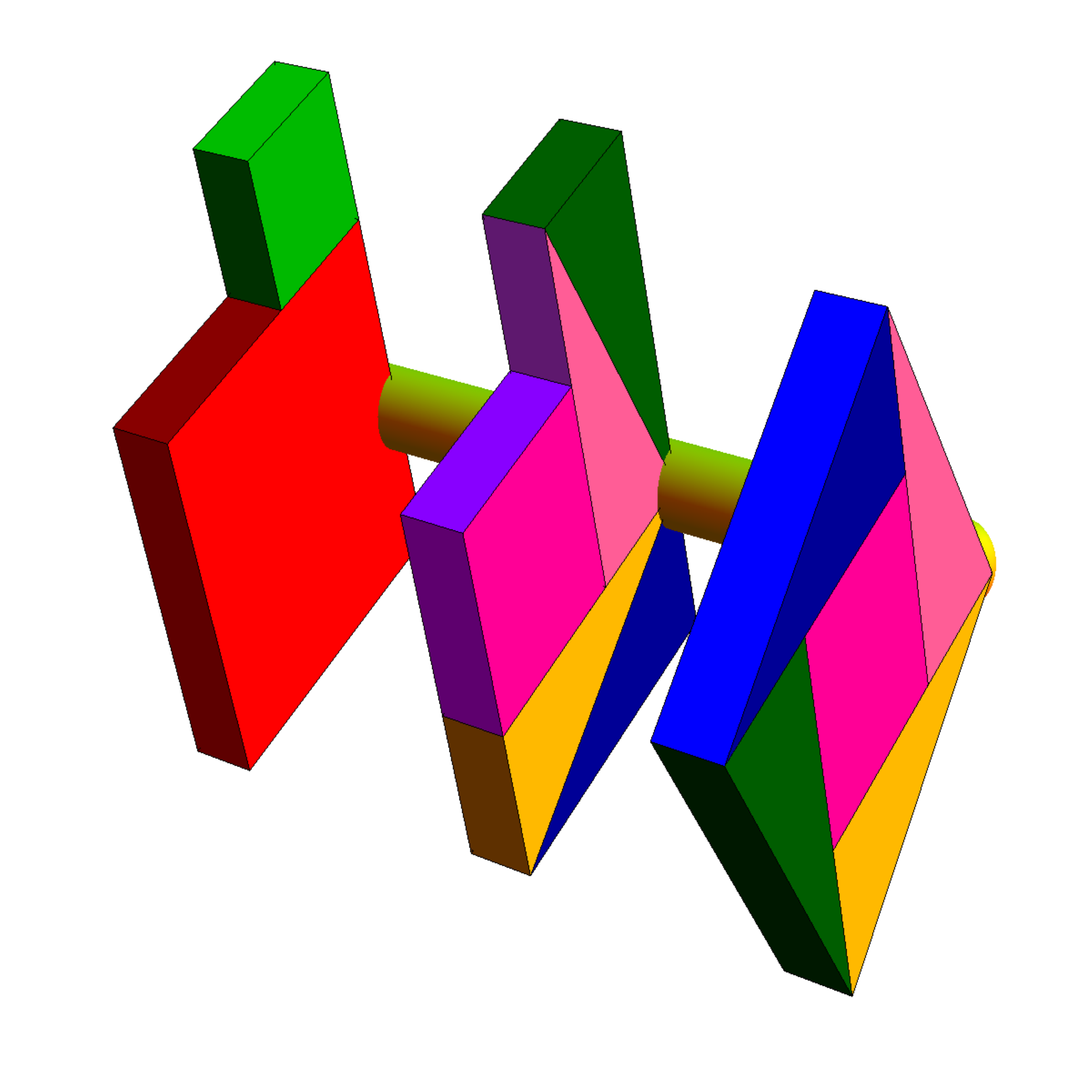}}
\scalebox{0.121}{\includegraphics{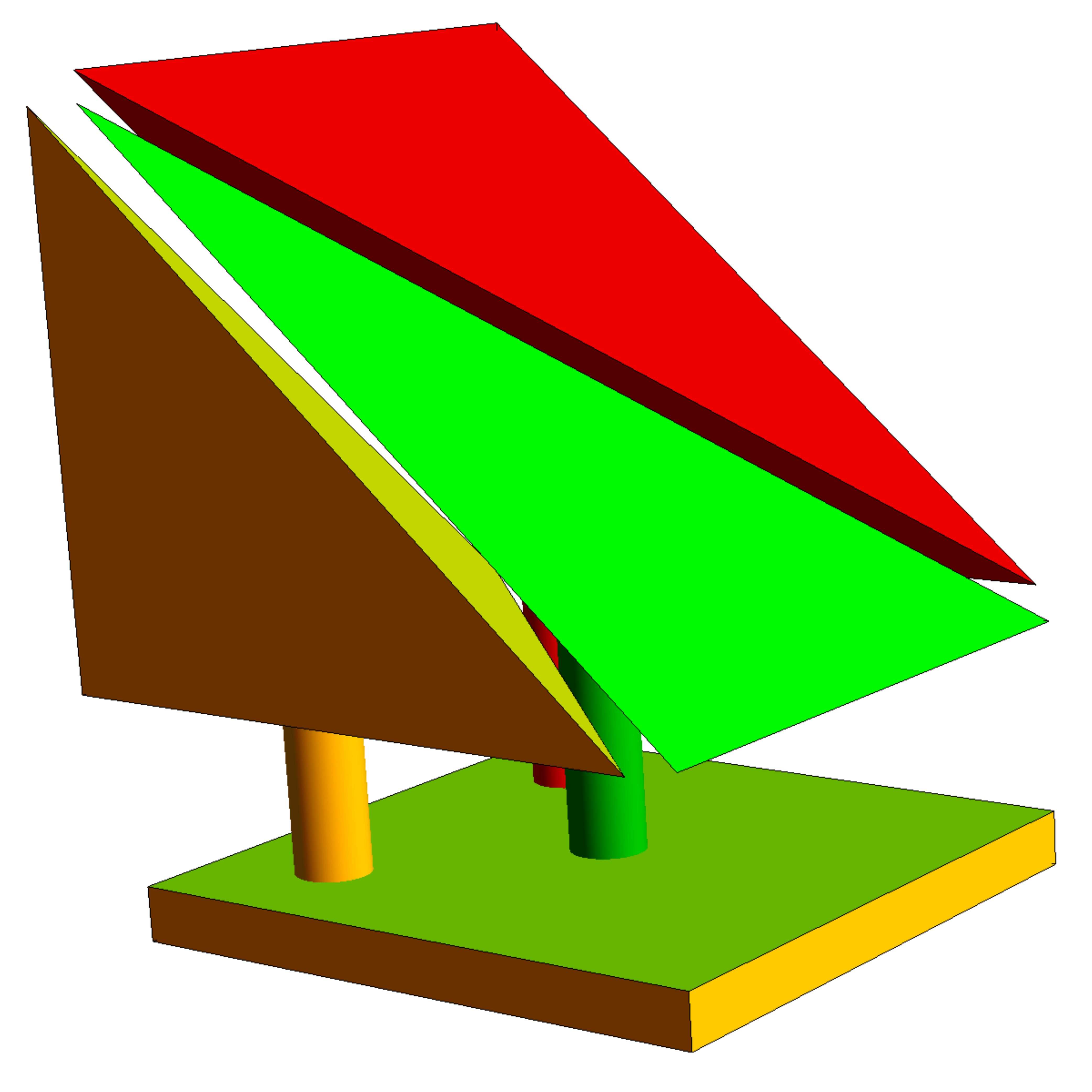}}
} \parbox{6cm}{
The left figure illustrates the proof of the Pythagoras theorem
\cite{Eves}.  The right figure is a proof that a pyramid has volume one third 
times the area of the base times height. 
}}

\parbox{16.8cm}{ \parbox{10cm}{
\scalebox{0.121}{\includegraphics{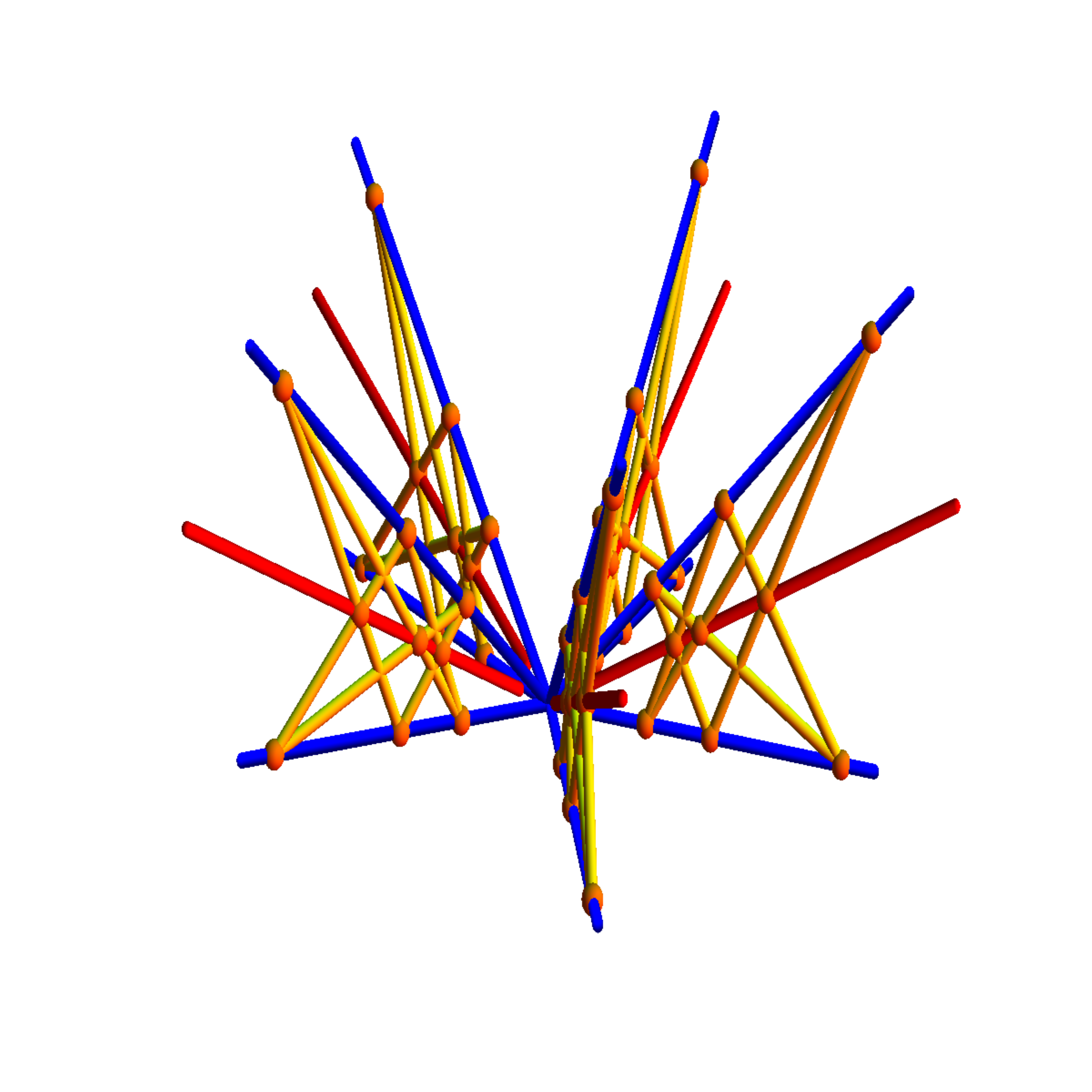}}
\scalebox{0.121}{\includegraphics{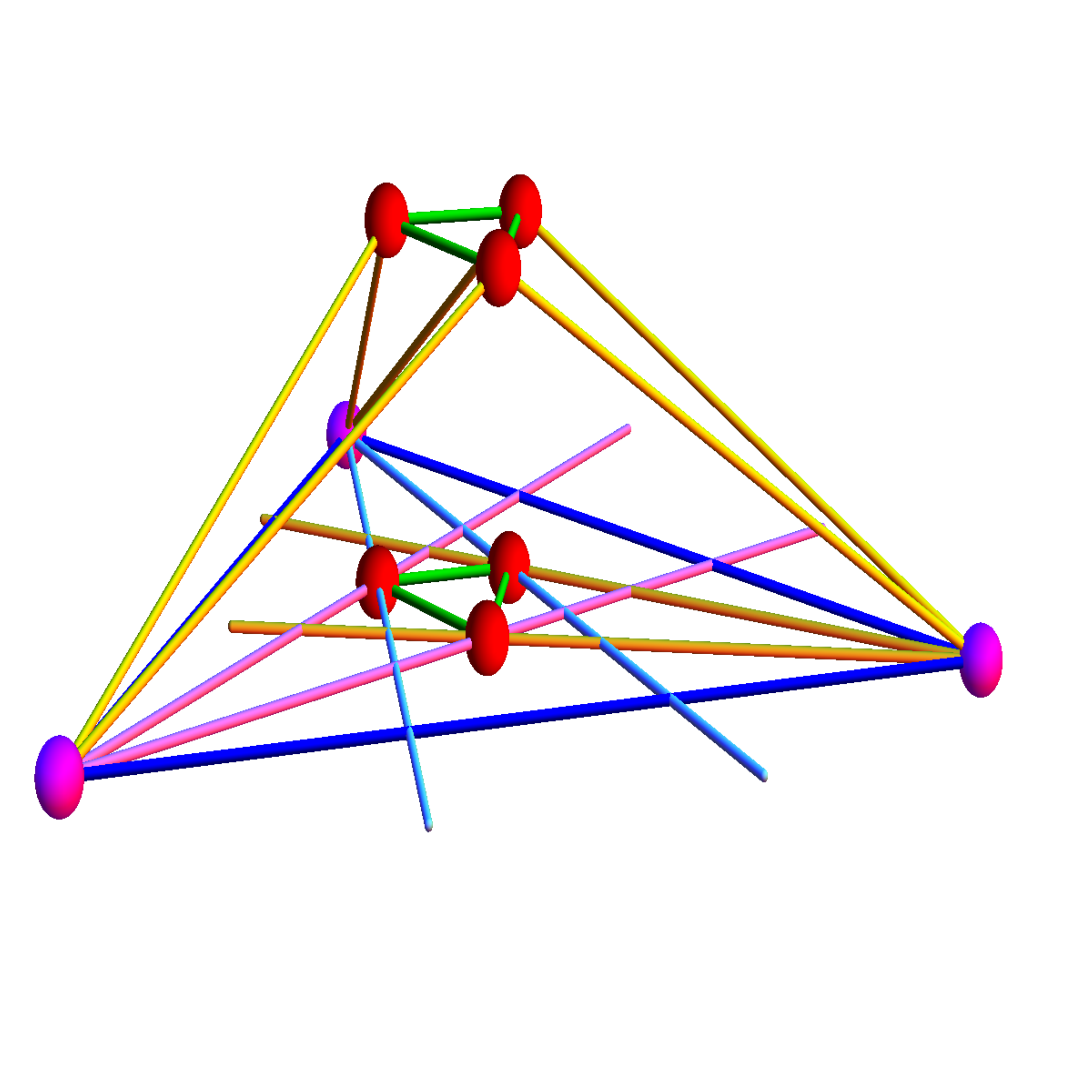}}
} \parbox{6cm}{
A theme on the Pappus theorem to the left and an illustration 
of the Morley theorem which tells that the angle trisectors of
an arbitrary triangle meet in an equilateral triangle. 
}}

\parbox{16.8cm}{ \parbox{10cm}{
\scalebox{0.121}{\includegraphics{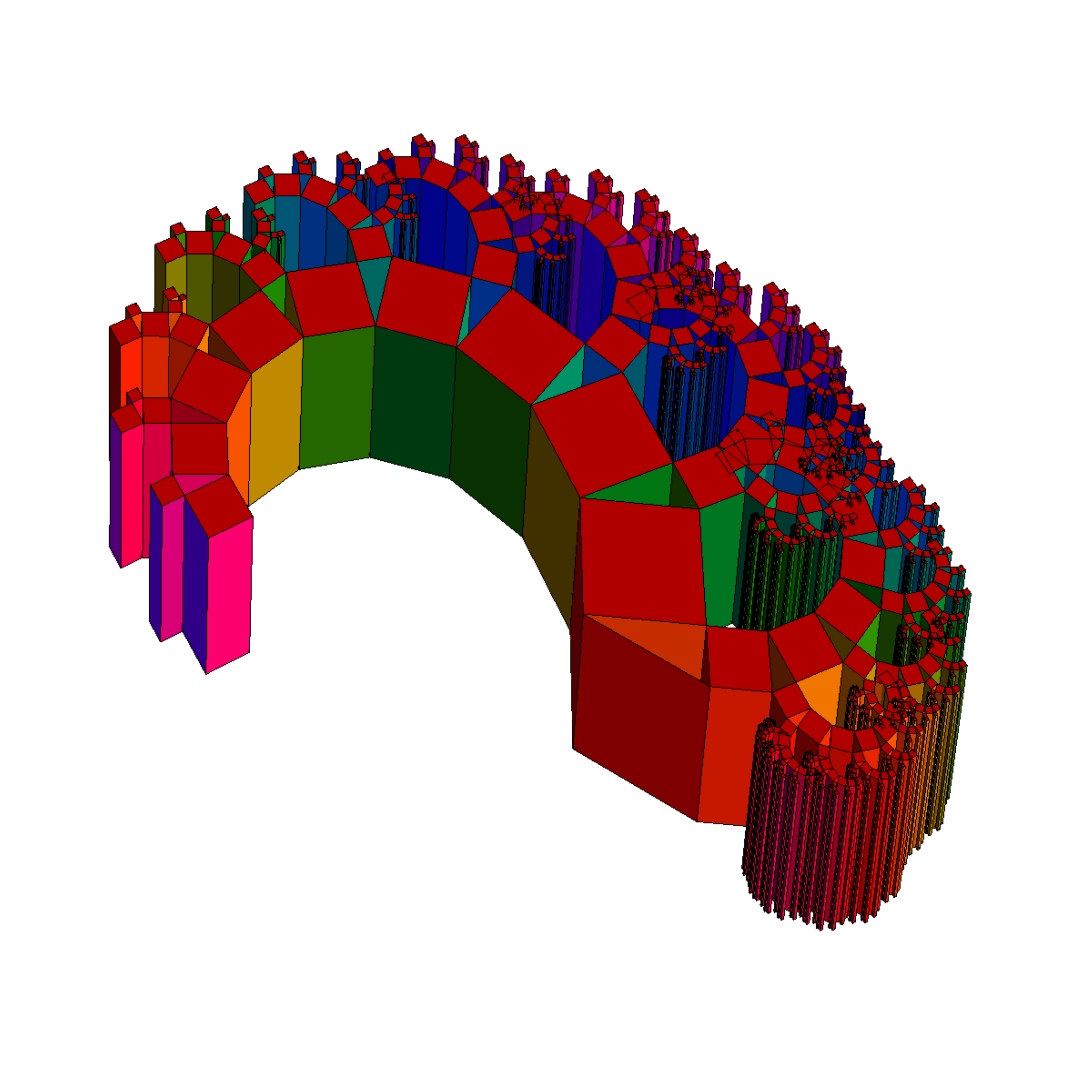}}
\scalebox{0.121}{\includegraphics{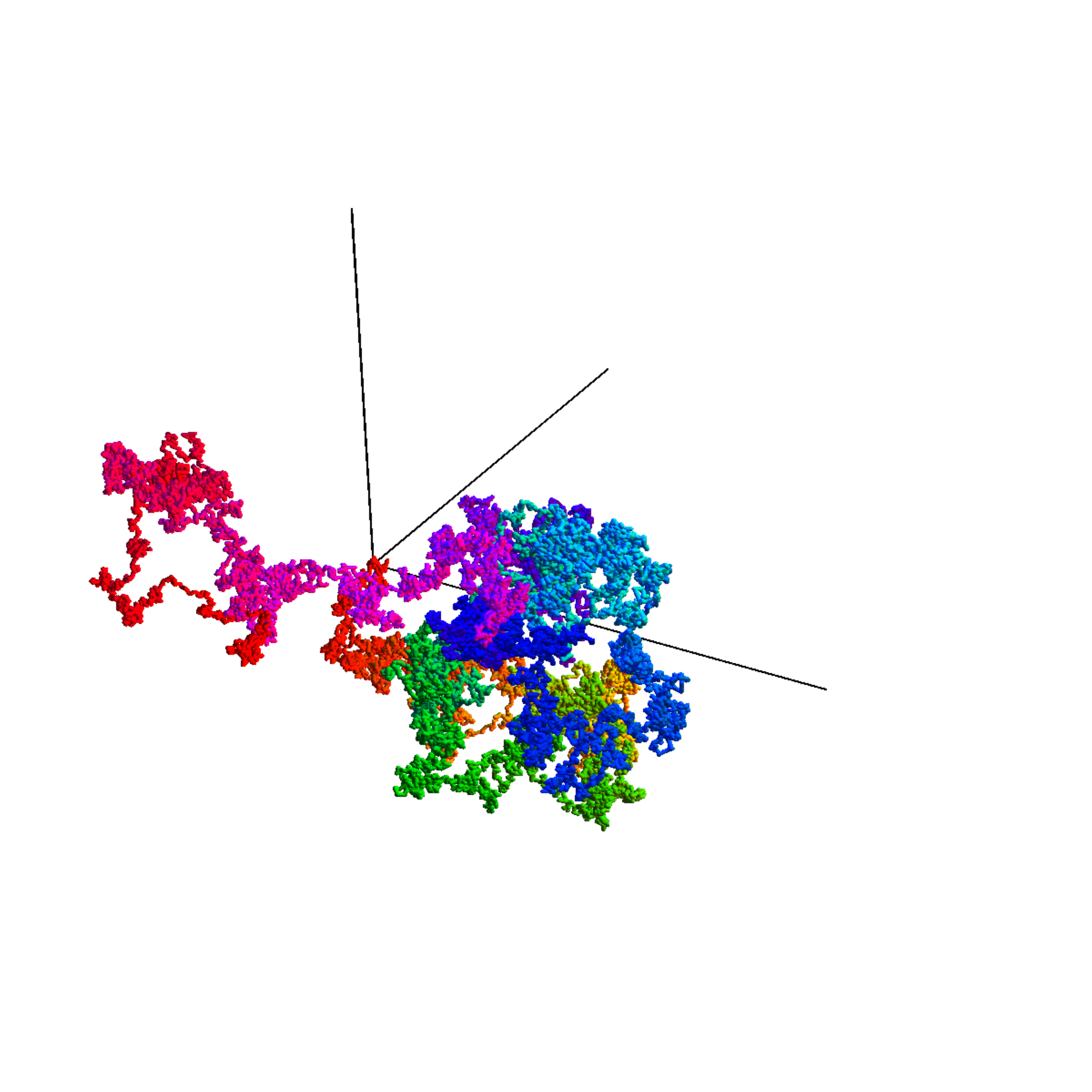}}
} \parbox{6cm}{
The left figure shows a fractal called the ``tree of pythagoras".
The right picture shows the random walk
in three dimensions. Unlike in dimensions one or two, the 
random walker in three dimensions does not return with 
probability $1$ \cite{Feller}.
}}

\parbox{16.8cm}{ \parbox{10cm}{
\scalebox{0.12}{\includegraphics{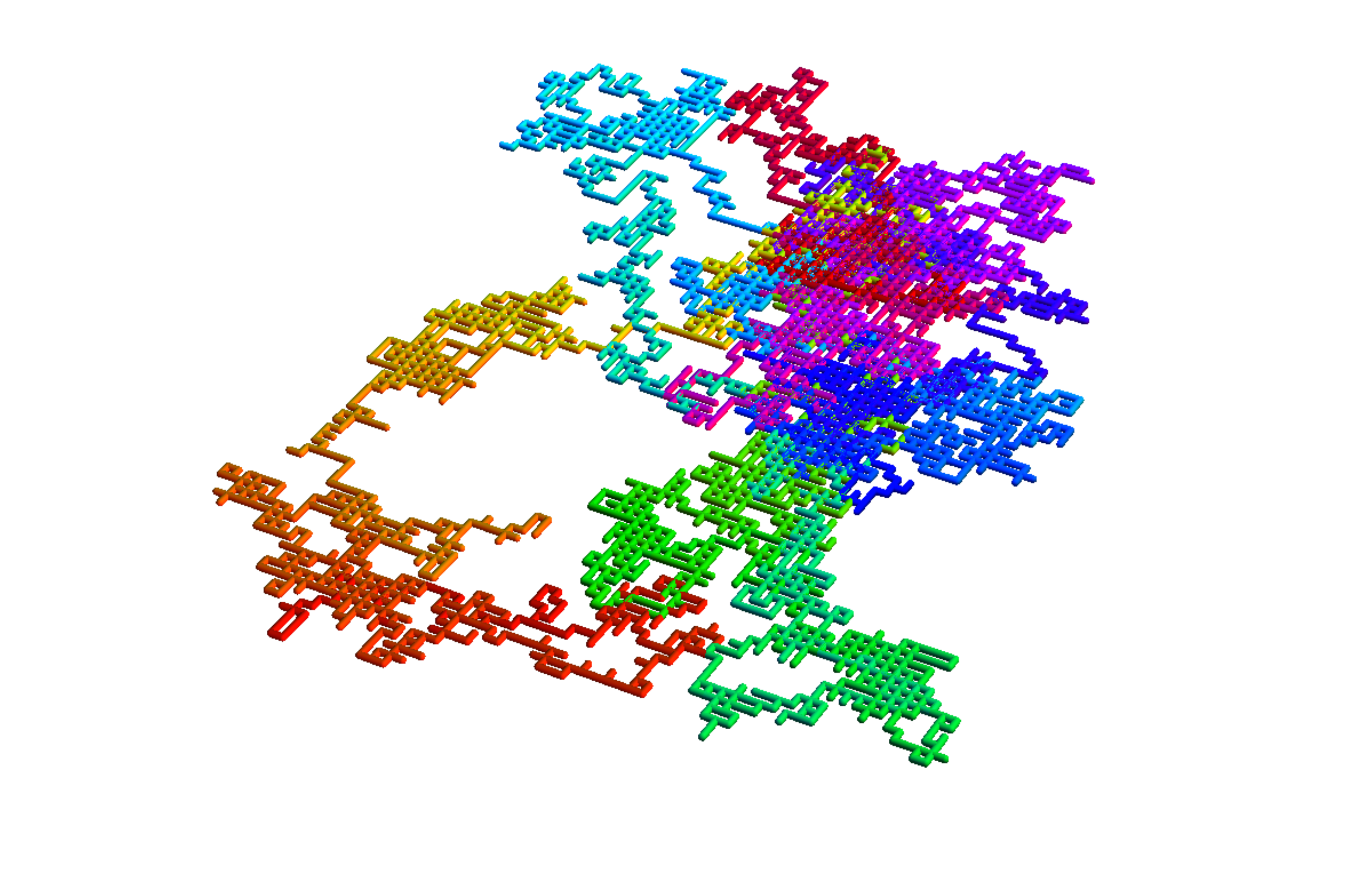}}
\scalebox{0.12}{\includegraphics{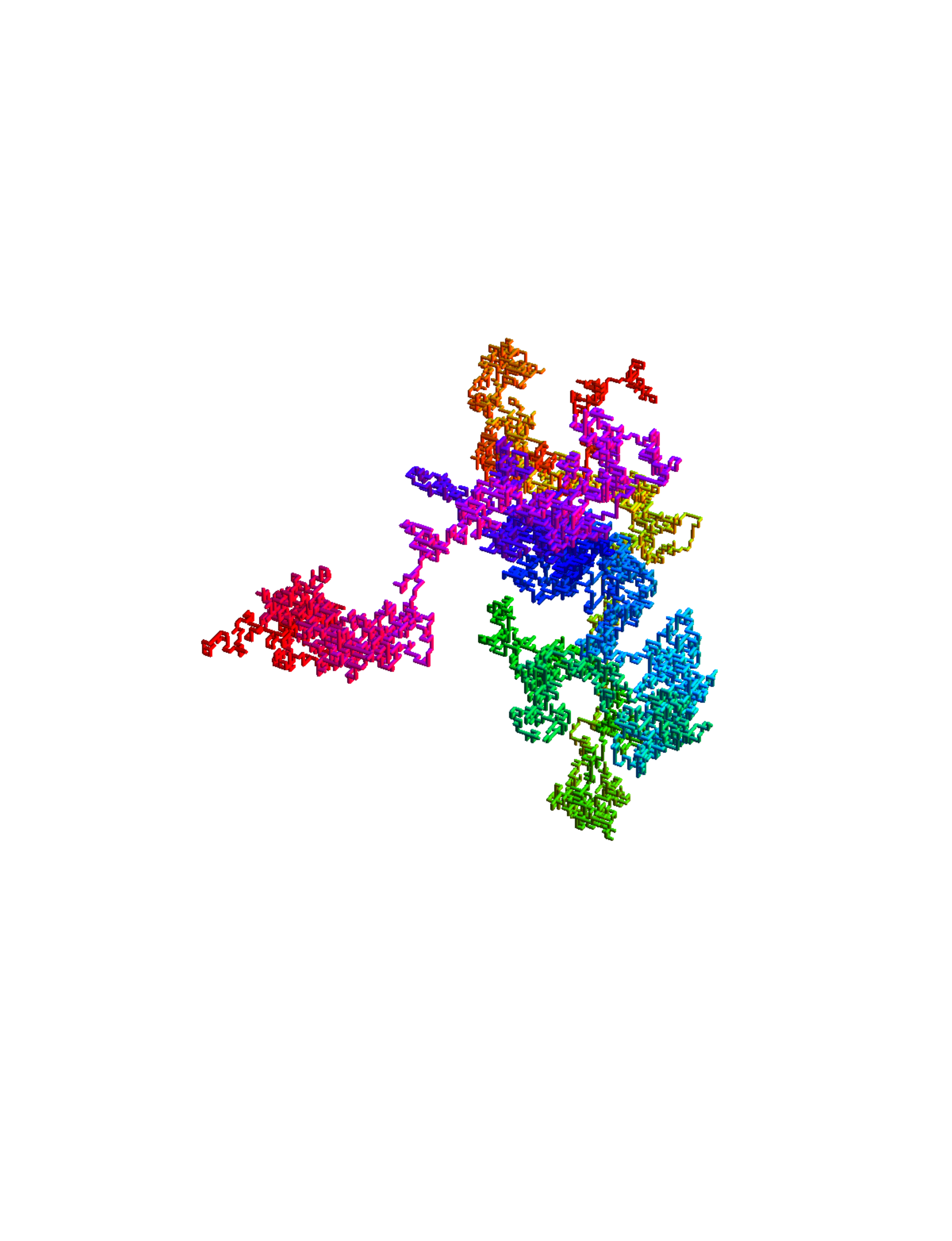}}
} \parbox{6cm}{
The number $\pi$ is believed to be normal in any number system.
An expansion in the number system 4 produces a random walk in the
plane \cite{BBCDDY}. In the number system 6, one sees a random walk in space.
One can print these random walks.
}}

\parbox{16.8cm}{ \parbox{10cm}{
\scalebox{0.12}{\includegraphics{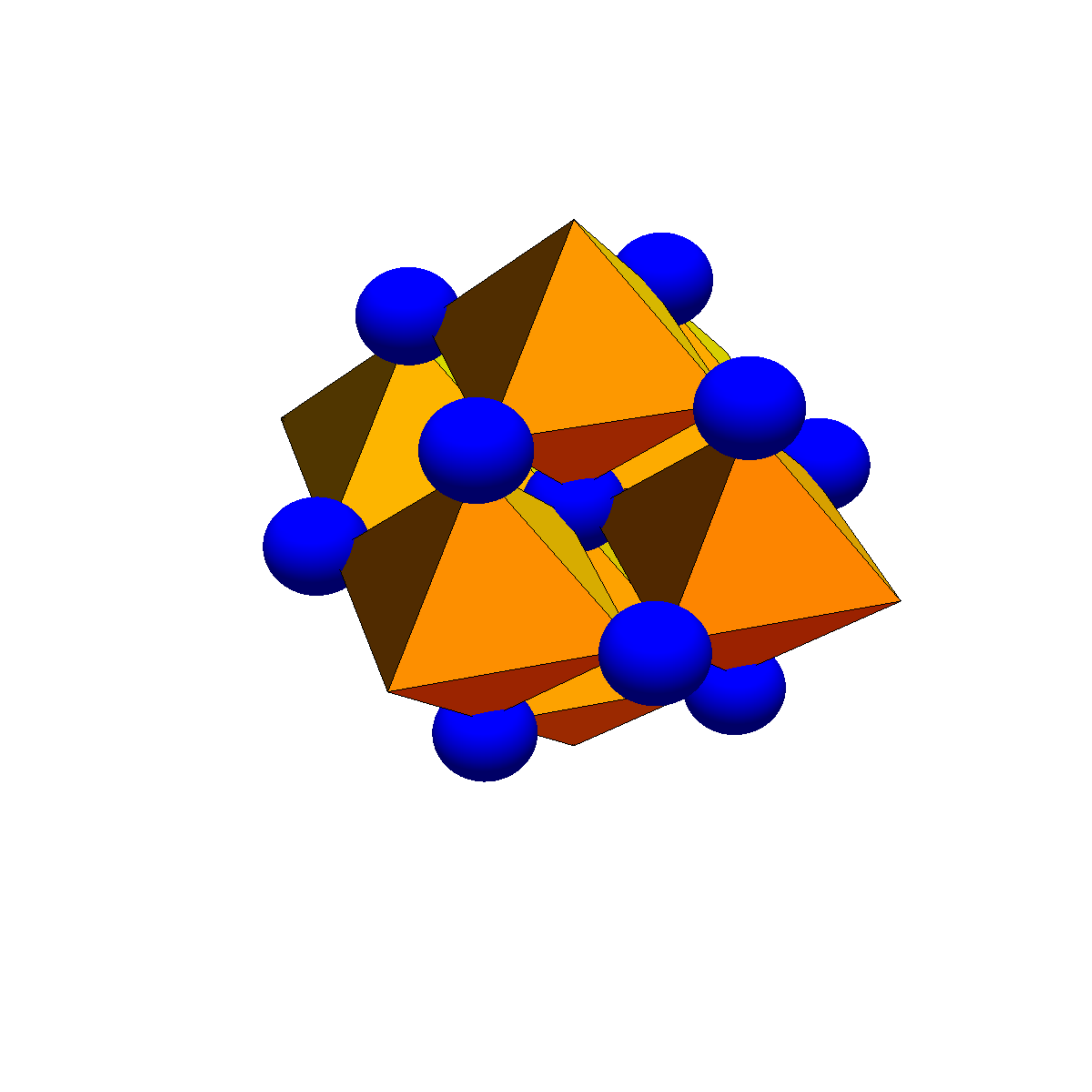}}
\scalebox{0.12}{\includegraphics{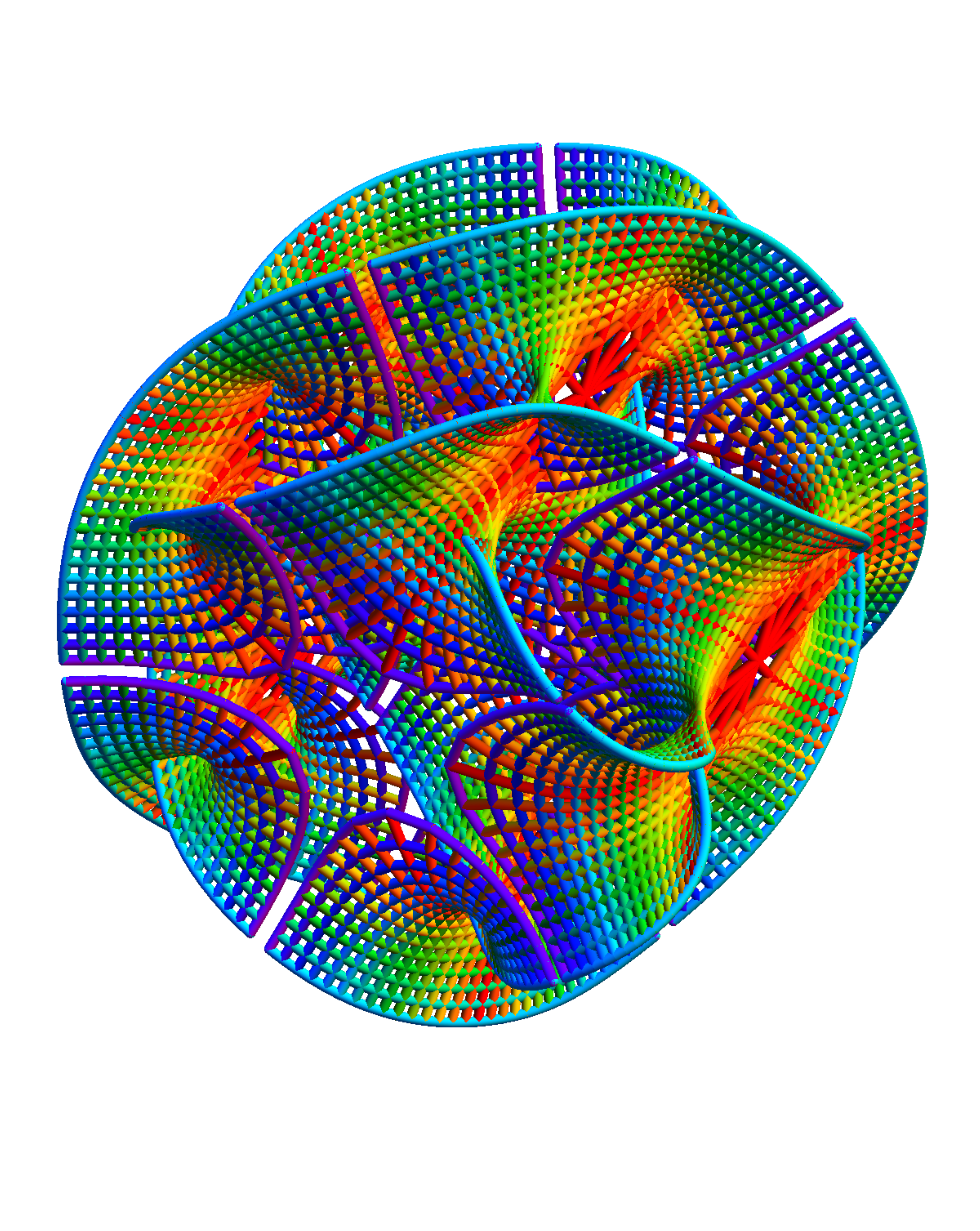}}
} \parbox{6cm}{
A tessellation of space by an alternated cubic honeycomb.
The tessellation consists of alternating octahedra and tetrahedra.
The second picture shows a Calabi-Yau surface. It is an icon
of string theory. In an appendix, we have the source code which 
produced this graphics. 
}}

\parbox{16.8cm}{ \parbox{10cm}{
\scalebox{0.12}{\includegraphics{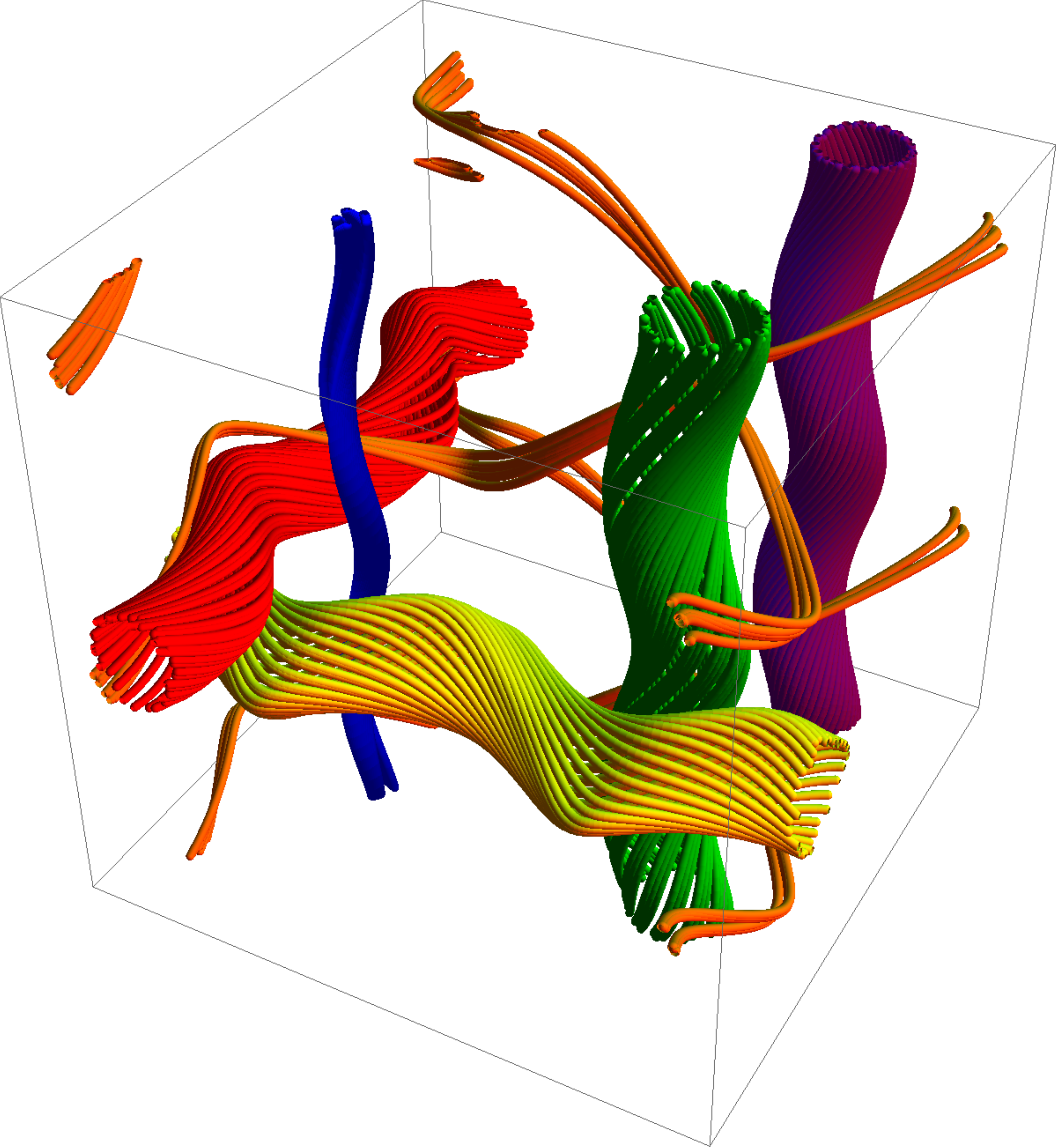}}
\scalebox{0.12}{\includegraphics{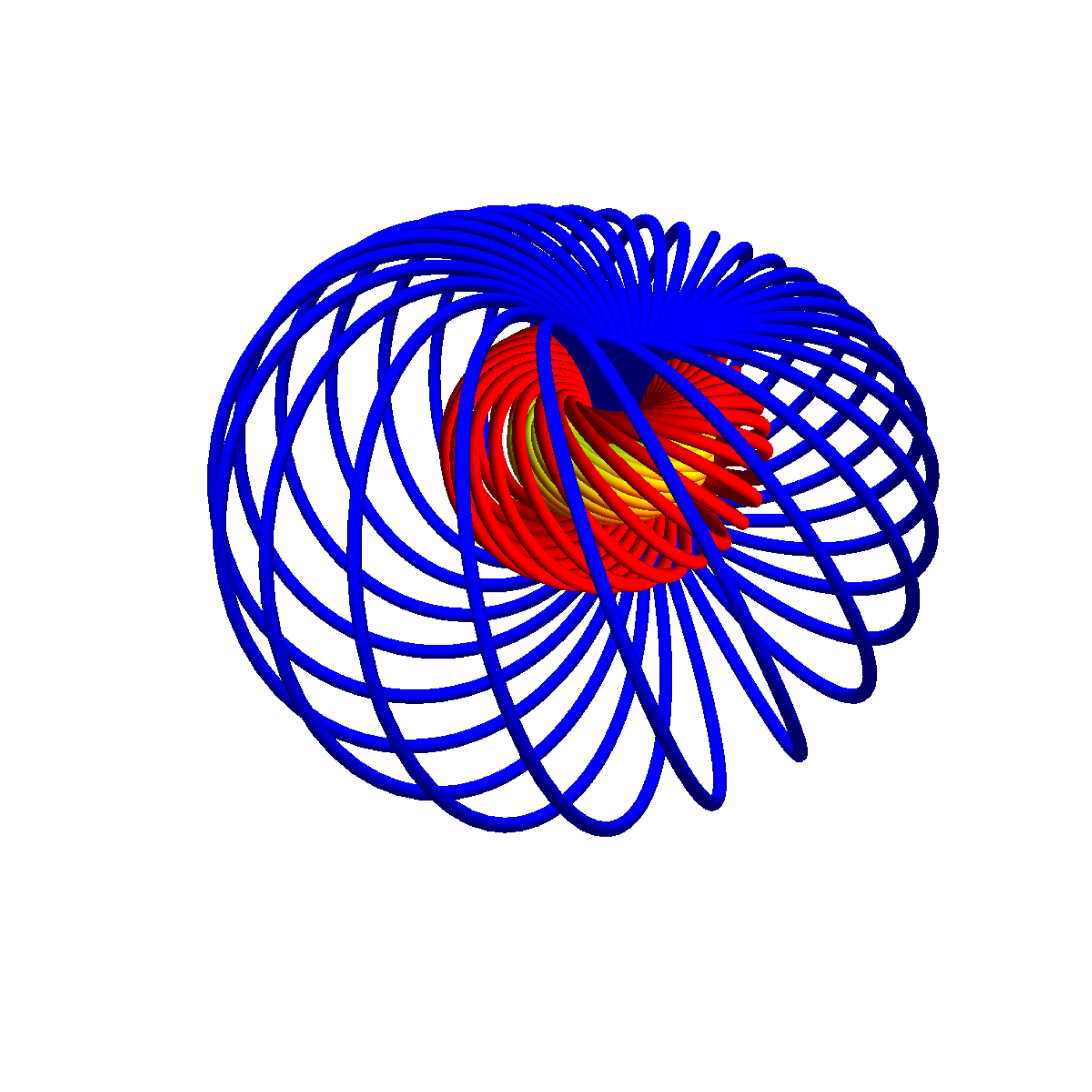}}
} \parbox{6cm}{
The ABC flow is a differential equation on the three dimensional
torus which is believed to show chaos. It is the simplest
volume preserving system on the torus.  We see 6 orbits.
The right picture shows the famous Hopf fibration. It visualizes
the three dimensional sphere. The picture is featured on the cover
of \cite{Mackenzie}.
}}

\parbox{16.8cm}{ \parbox{10cm}{
\scalebox{0.12}{\includegraphics{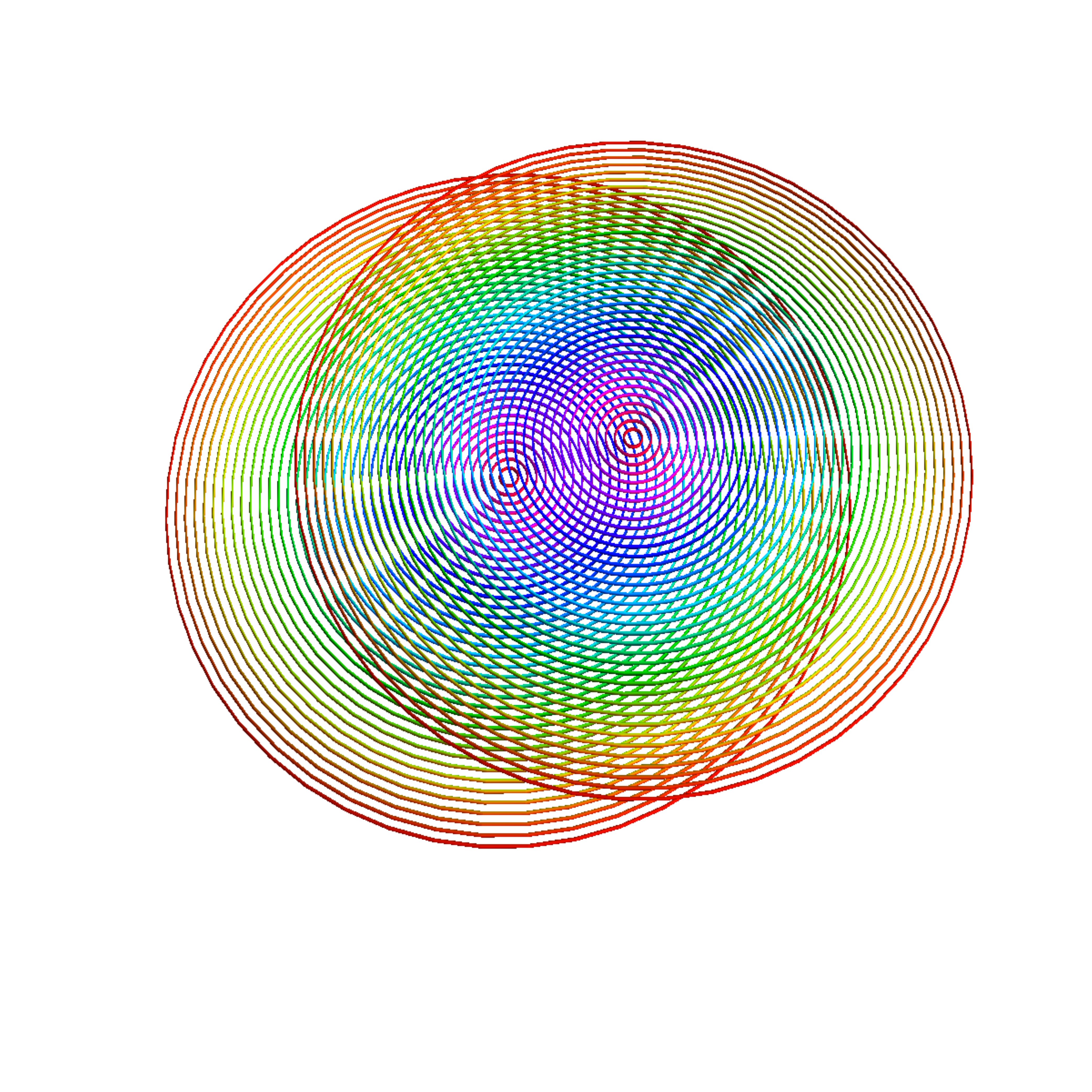}}
\scalebox{0.12}{\includegraphics{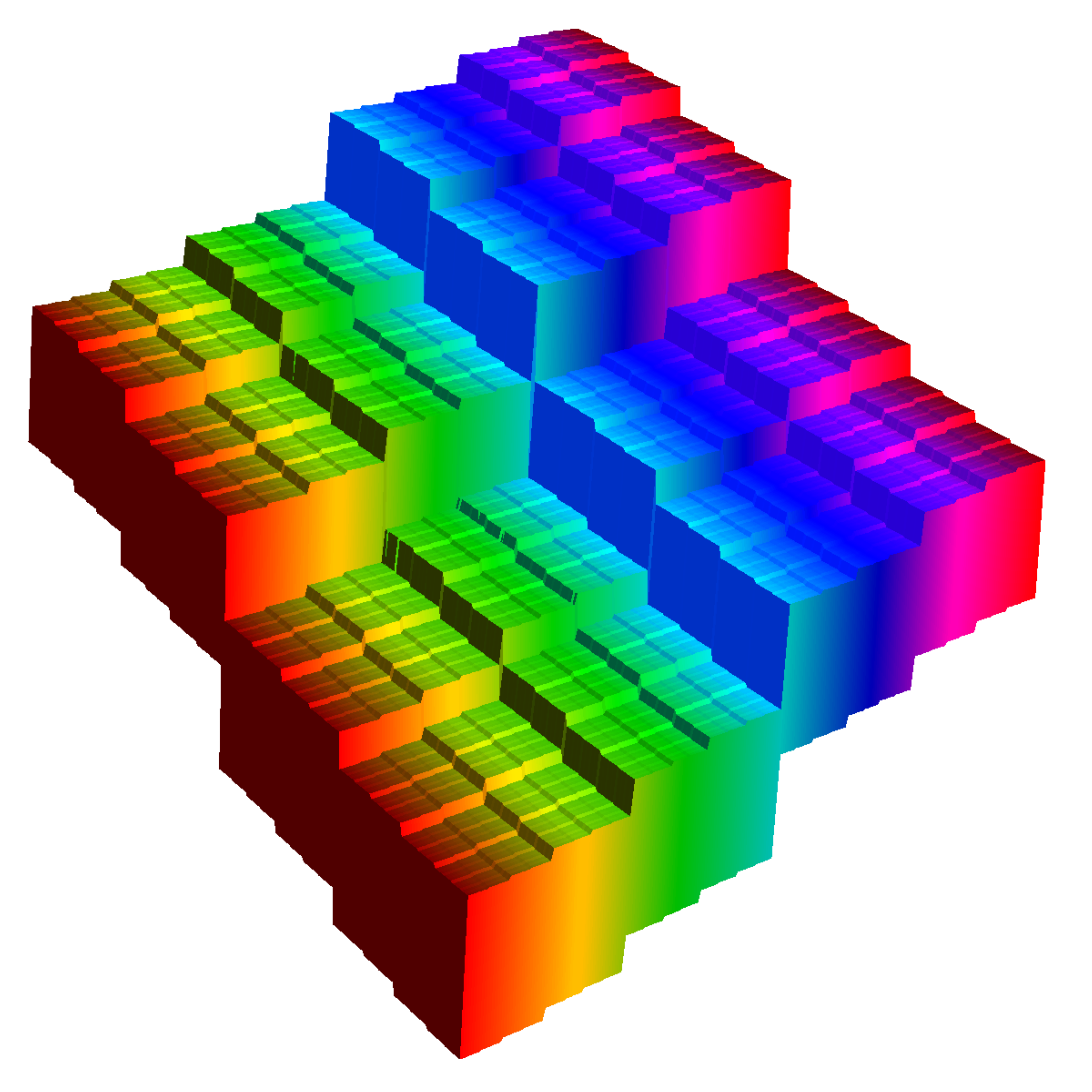}}
} \parbox{6cm}{
The left picture shows two circular overlaying patterns
producing the Moir\"e effect. The right picture shows
a double devil staircase. Take a devil staircase $g(x)$
which has jump discontinuity at all rational points, then
define $f(x,y) = g(x) + g(y)$. It is a function of two variables
which is discontinuous on a dense set of points in the domain.
}}

\parbox{16.8cm}{ \parbox{10cm}{
\scalebox{0.12}{\includegraphics{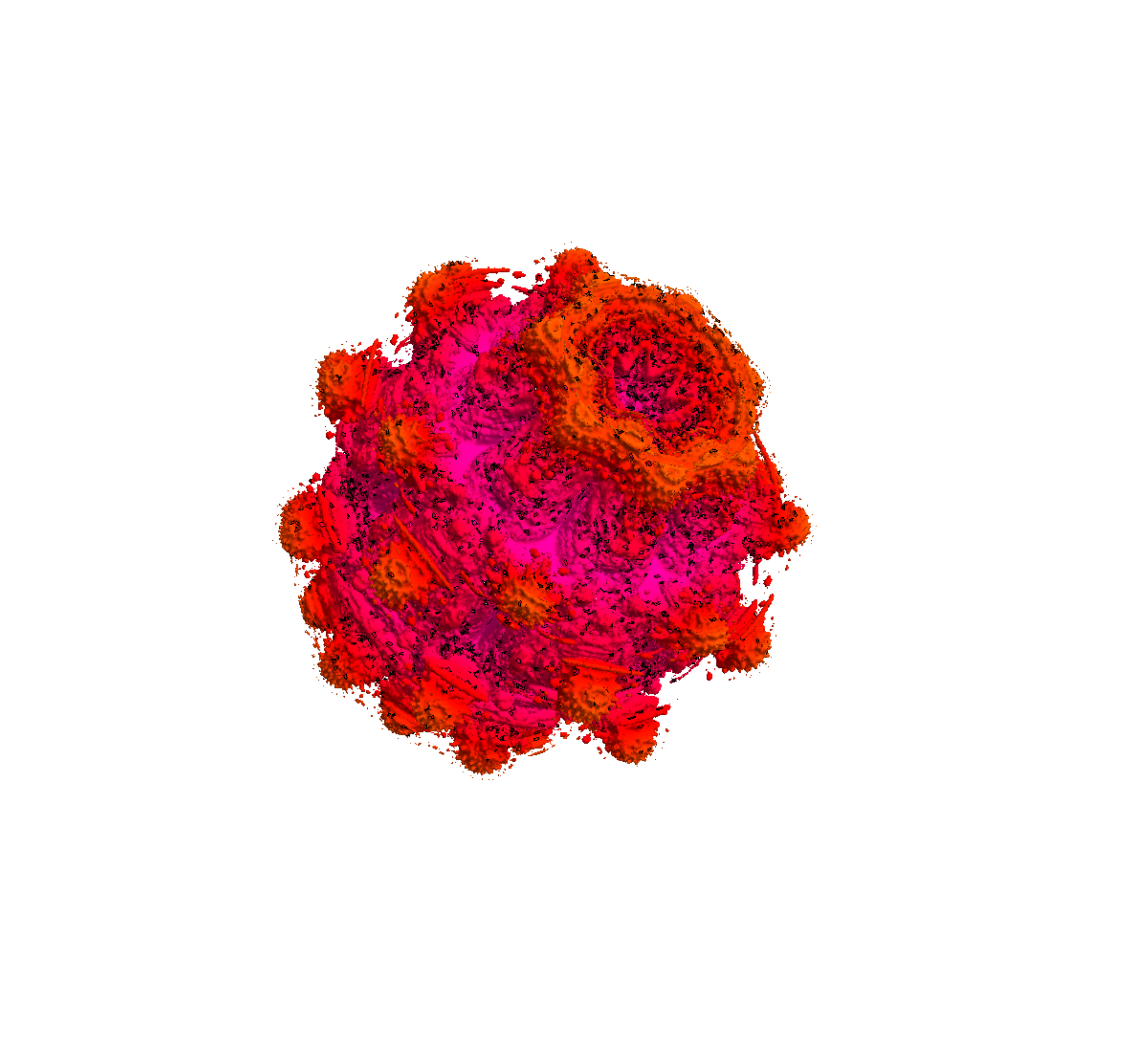}}
\scalebox{0.11}{\includegraphics{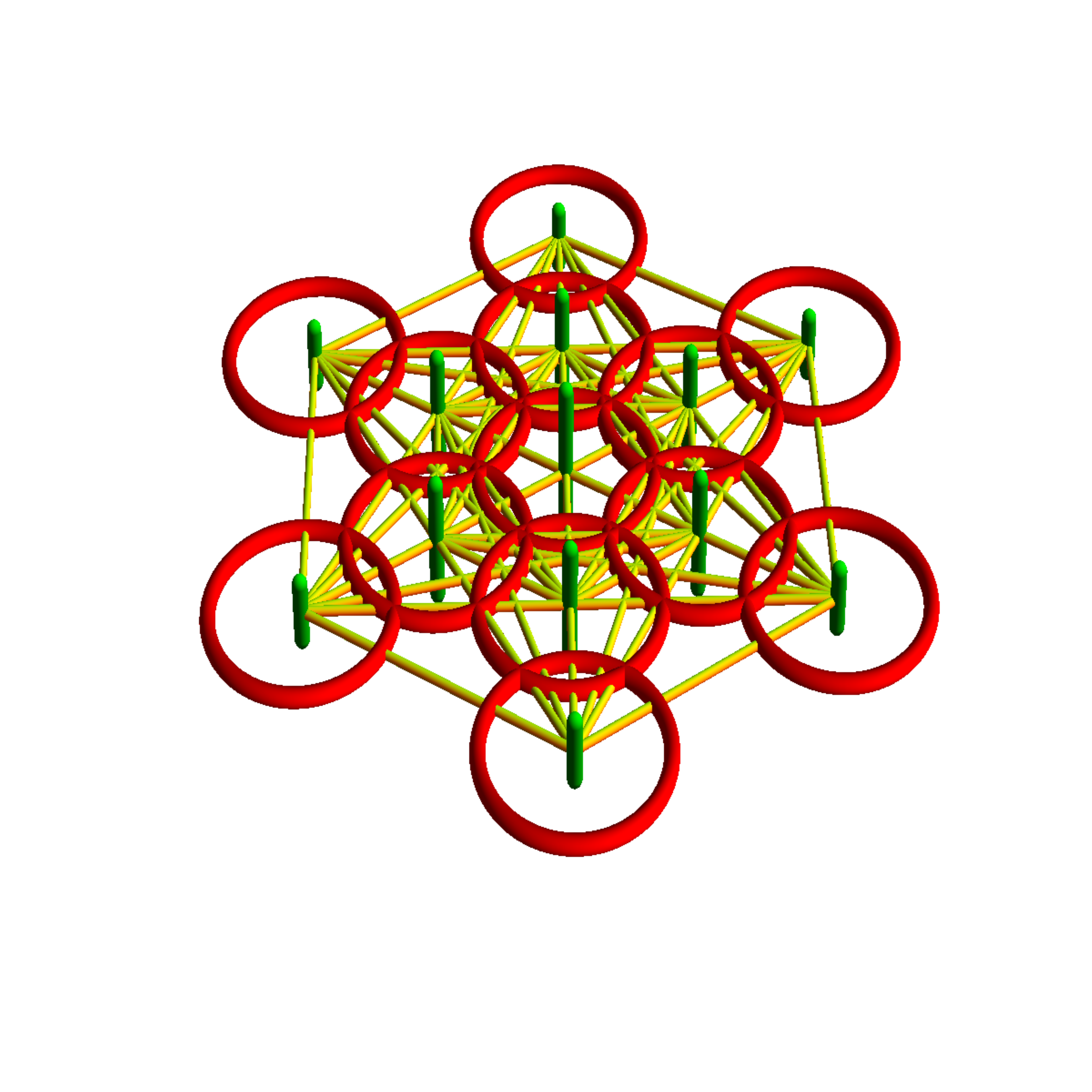}}
} \parbox{6cm}{
The Mandelbulb $M_8$ is defined as the set of vectors $c \in R^3$
for which $T(x)=x^8+c$ leads to a bounded orbit of $0=(0,0,0)$, where $T(x)$
has the spherical coordinates $\rho^8,8 \phi, 8\theta)$ if $x$ has the
spherical coordinates $rho,\phi,\theta)$. The topology of $M_8$ is unexplored.
The right figure shows the ``Metatron", sacred-mathematical inspiration for 
artists (see also chapter 14 in \cite{Grattan-Guinness}). 
For less mystical folks, it makes a frame for a remote controlled drone. }}

\parbox{16.8cm}{ \parbox{10cm}{
\scalebox{0.12}{\includegraphics{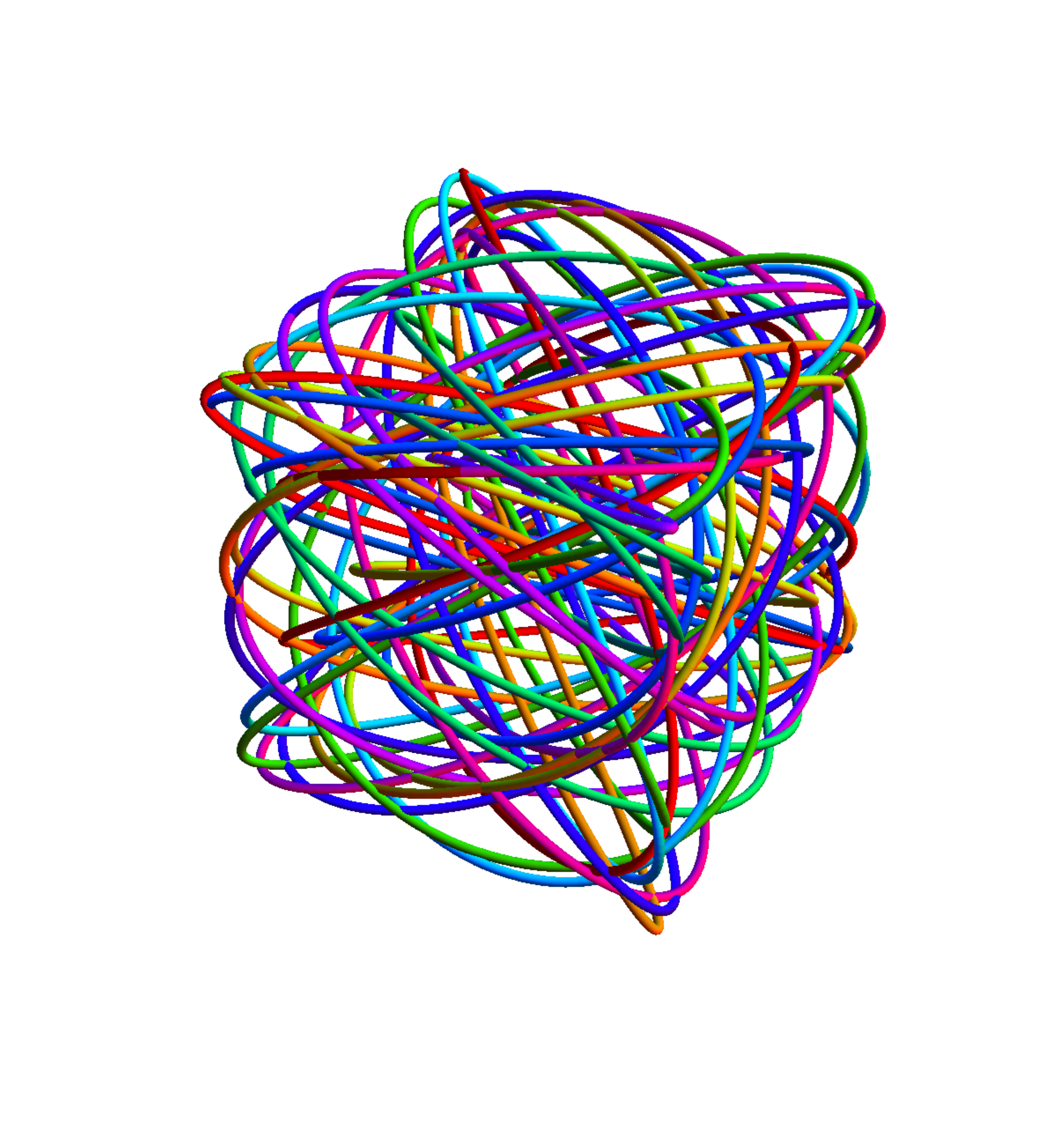}}
\scalebox{0.12}{\includegraphics{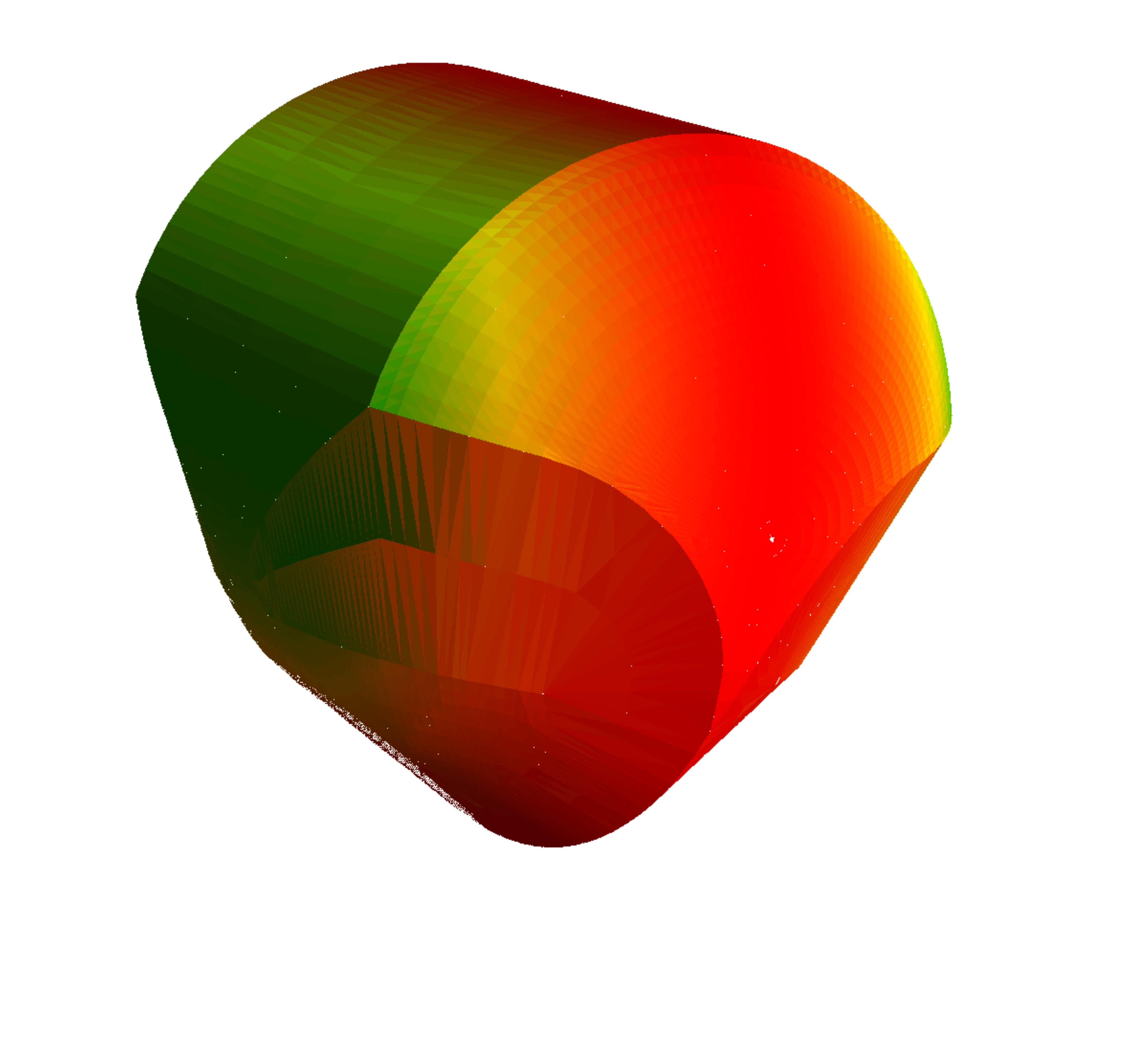}}
} \parbox{6cm}{
A Lissajoux figure in 3D given by the biorthm curve
$\vec{r}(t) = ( \sin(2\pi x/23),\sin(2\pi x/28),\sin(2\pi x/33) )$
measuring physical, emotional and intellectual strength $x$ 
days after birth. While hardly taken seriously by anybody, it can be 
a first exposure to trigonometry. The right picture shows the G\"omb\"oc \cite{Gomboc}, 
a convex solid which has one stable and one unstable equilibrium point.
It was constructed in 2006 \cite{VD}.  }} 

\parbox{16.8cm}{ \parbox{10cm}{
\scalebox{0.10}{\includegraphics{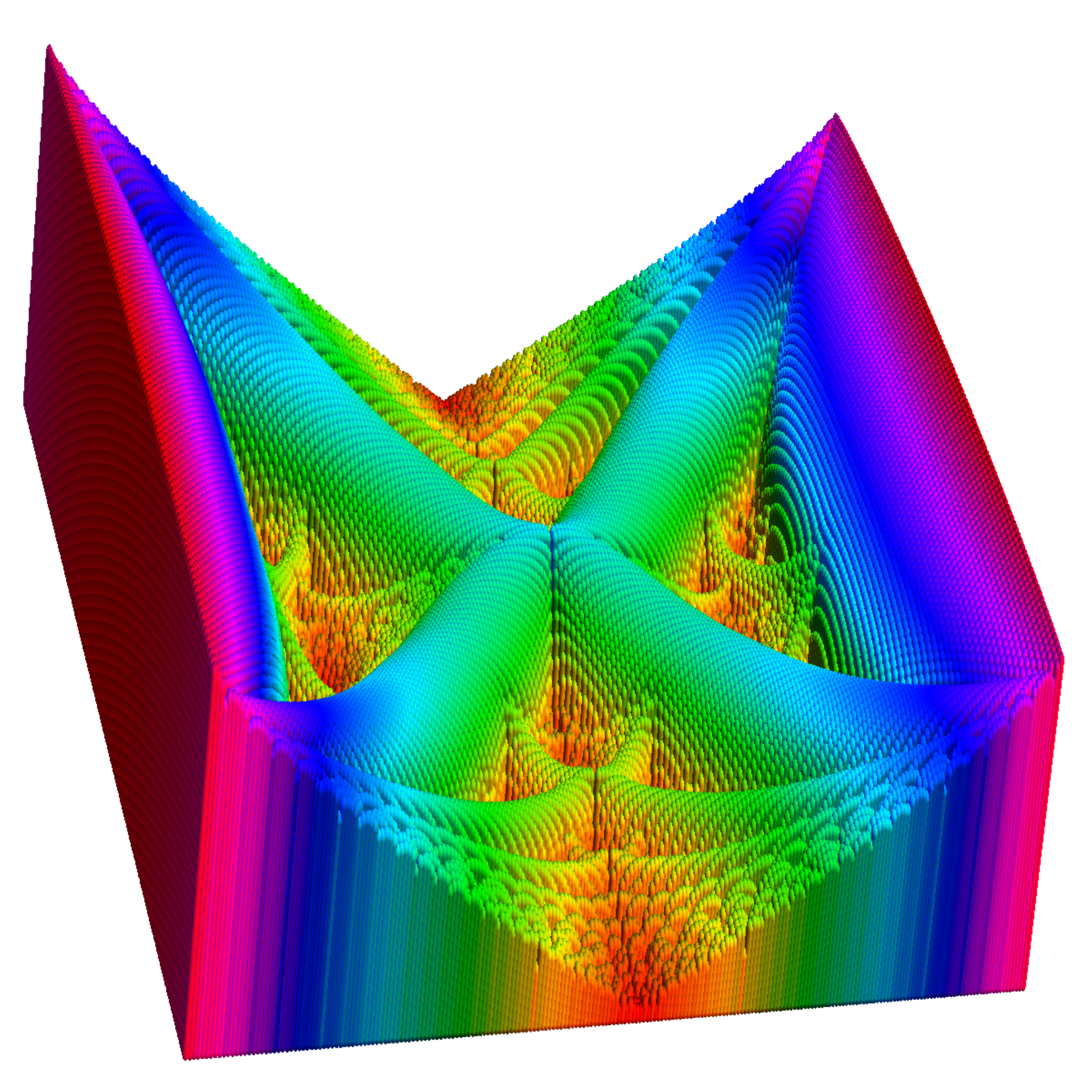}}
\scalebox{0.13}{\includegraphics{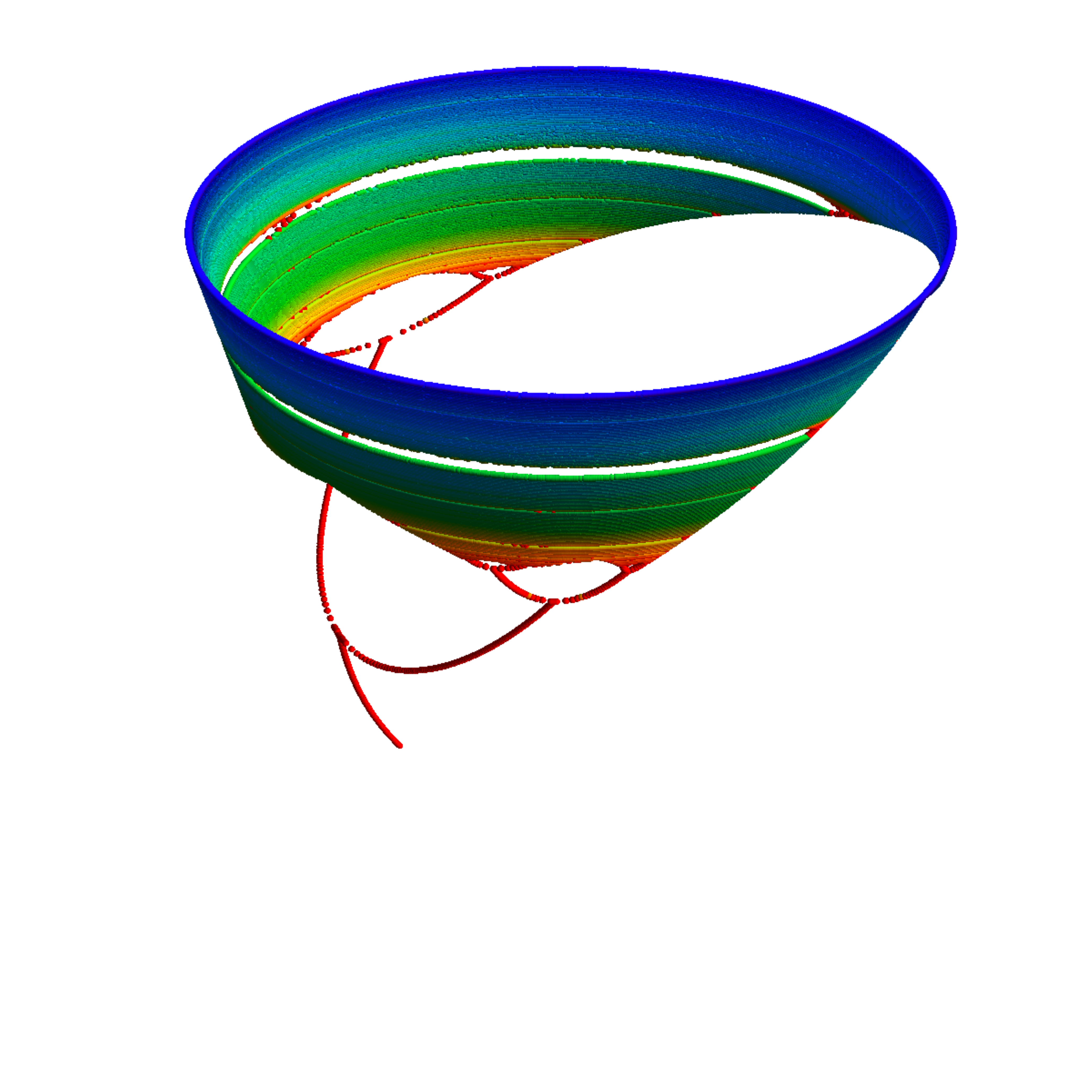}}
} \parbox{6cm}{
The Hofstadter butterfly popularized in \cite{Hofstadter} to the left
encodes the spectrum of the discrete Mathieu operator. The
butterfly is a Cantor set \cite{Tenmartini} of zero Lebesgue measure
\cite{ZeroMeasure}, its wings are still not understood \cite{LastReview,ExoticSpectra}.
To the right, the bifurcation diagram of the one parameter
family of circle maps $T_c(x) = c \sin(x)$, which features universality. The
bifurcations happen at parameter values such that the quotient of 
adjacent distances converges \cite{Feigenbaum1978}.
}}

\parbox{16.8cm}{ \parbox{10cm}{
\scalebox{0.10}{\includegraphics{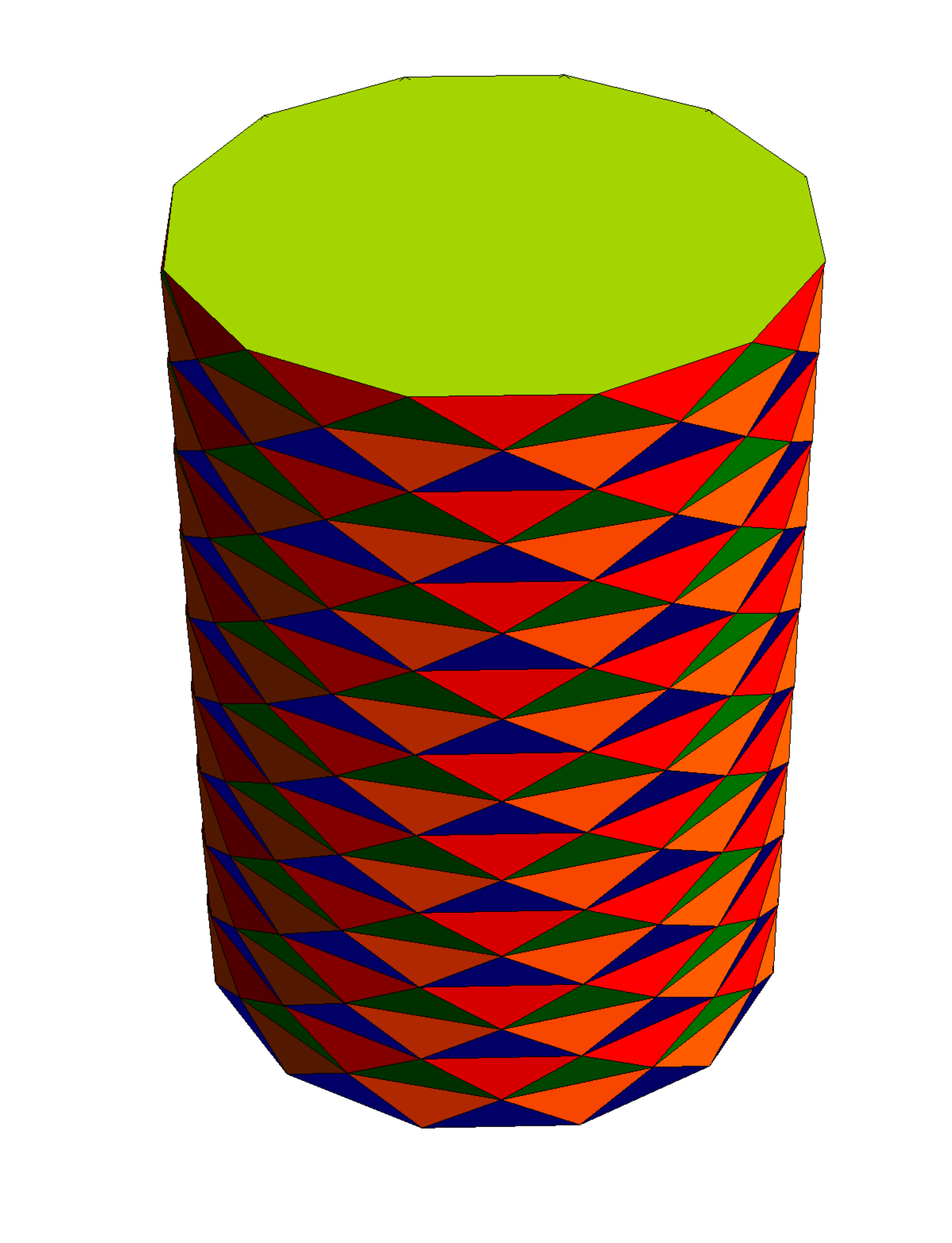}}
\scalebox{0.14}{\includegraphics{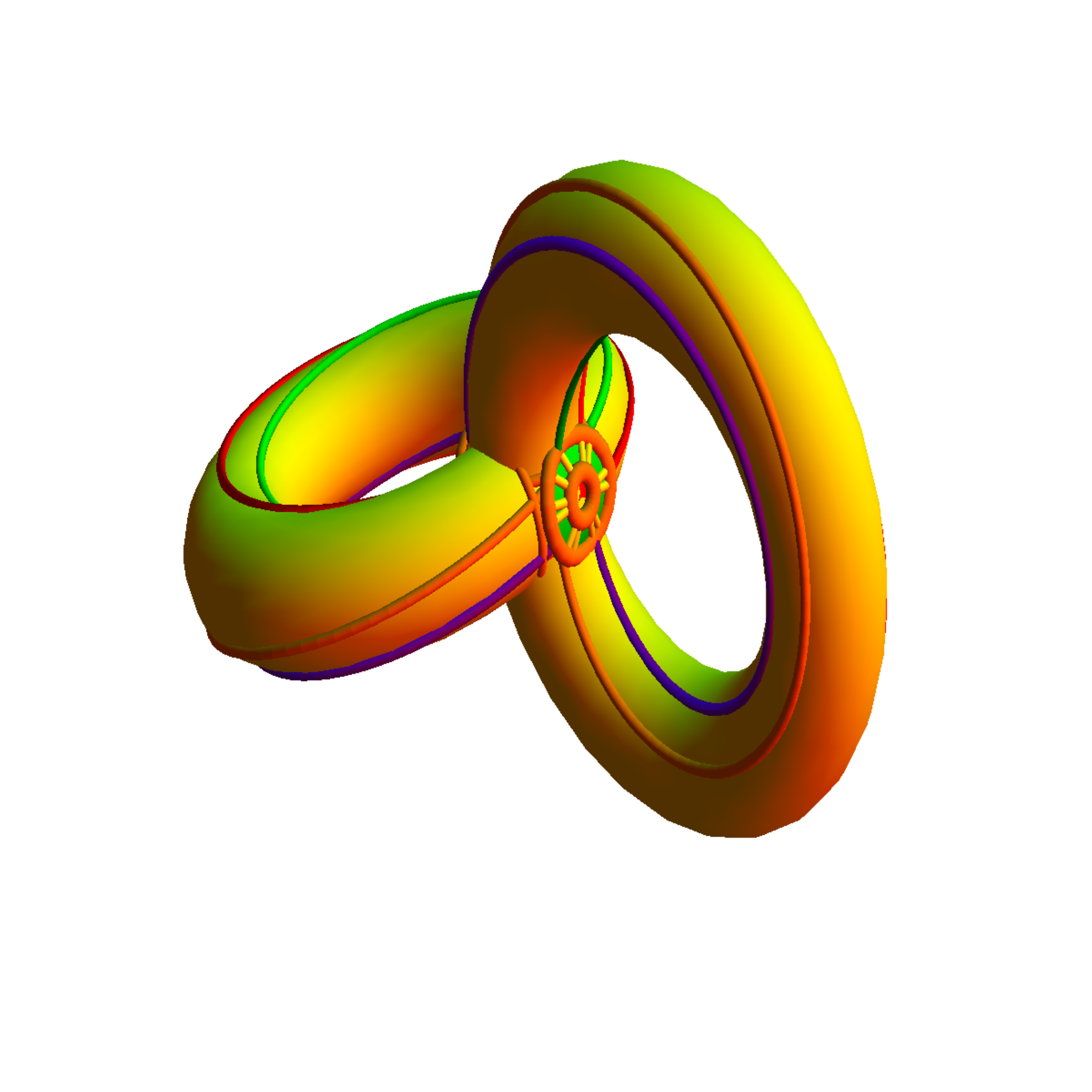}}
} \parbox{6cm}{
To the left the H.A. Schwarz example from 1880 
\cite{Blatter,Plateau}. It illustrates that the surface area 
of a triangularization does not need to converge to the surface 
area of the surface. The area of Chinese lantern triangularizations
can be arbitrary large, even it approximates the surface arbitrary well. 
The right figure shows the Tucker group, the only group of genus 2
\cite{Tucker,Ferguson}. The group had been realized as a sculpture 
by De Witt Godfrey and Duane Martinez. }}

\pagebreak

\section{Tables and code snippets}
\appendix{A) {\bf Revolutions.}}  The first table summarizes information and industrial revolutions. \\

\begin{tabular}{|l|l|l|l|}   \hline
{\bf Information revolutions}  &         &   {\bf Industrial revolutions}     &          \\  \hline
Gutenberg Press                & 1439    &   Steam Engine, Steel and Textile  &  1780\\
Mechanical Computer            & 1642    &   Automotive, Chemistry            &  1850\\
Personal Computer, Cell phone  & 1973    &   Personal Computer, Rapid Prototyping &  1969\\ \hline
\end{tabular}
\vspace{2mm}

For industrial revolutions, see \cite{Deane} page 3, for the second industrial 
revolution \cite{Levin} page 2, for the third one \cite{GRS}  page 34,\cite{worldfinancialreview}.

\vspace{2mm}

\appendix{B) {\bf Changes in communication, perception and classroom.}} This table gives examples of breakthroughs in 
communication and in the classroom. The middle number indicates how many years ago, the event happened. \\

\begin{tabular}{|c|c|c|}   \hline
{\bf Communication}  &
{\bf Perception}    & 
{\bf Classroom}   \\  \hline
\begin{tabular}{lrl}
Alphabet    &  20K&   Ishango bone    \\
Figures     &   6K&   Clay tablet \\
Models      &   2K&   Apollonius \\
Books       &  560&   Gutenberg  \\
Photo       &  170&   Photographia \\
Film        &  130&   Kineograph \\
3D-Print    &   30&   Stereolithog.
\end{tabular}
&
\begin{tabular}{lrl}
Projection  &     1K&  Camera obscura \\
Eye glass   &    730&  Aalvina Armate \\
Microscope  &    420&  Janssen \\
Telescope   &    400&  Kepler \\
Xrays       &    110&  Roentgen  \\
MRI         &     60&  Catscan \\
3D scan     &     25&  Cyberware 
\end{tabular}
&
\begin{tabular}{lrl}
Models      &2K&  Greeks \\
Abacus      &1.5K&  Abacus  \\
Blackboard  &1K&  Tarikh Al-Hind \\
CAS         &    50&  Schoonship \\
Calculator  &    40&  Busicom  \\
Powerpoint  &    30&  Presenter  \\
3D-Models   &    15&  Makerbot
\end{tabular}  \\ \hline
\end{tabular}
\vspace{2mm}

For classroom technology and teaching see \cite{Toolsteaching}.
Early computer algebra systems (CAS) in the 1960ies were Mathlab, Cayley, Schoonship, 
Reduce, Axiom and Macsyma \cite{Weinzierl}. For Multimedia development in mathematics,
see \cite{BorweinMoraelsPOlthierRodrigues}.
The first author was exposed as a student to Macsyma, Cayley (which later became Magma) 
and Reduce. We live in a time when even the three categories start to blur: cellphones with visual 
and audio sensors, possibly worn as glasses connect to the web. In the classroom, teachers
already today capture student papers by cellphones and have it automatically graded. Students
write on intelligent paper and software links the recorded audio with the written text. A time
will come when students can print out a physics experiment and work with it.\\

\appendix{C) {\bf Source code for exporting an STL or X3D file.}} The following
Mathematica lines generate an object with 13 kissing spheres. \\

\hspace{-1mm}
\begin{small}
\begin{lstlisting}[language=Mathematica]
s=Table[{0,n,m*GoldenRatio},{n,-1,1,2},{m,-1,1,2}];
s=Partition[Flatten[s],3]; s=Append[s,{0,0,0}];
s=Union[s,Map[RotateRight,s],Map[RotateLeft,s]];
S=Graphics3D[Table[{Hue[Random[]],Sphere[s[[k]]]},{k,12}]]
Export["kissing.stl",S,"STL"]; Export["kissing.x3d",S];
\end{lstlisting}
\end{small}

\vspace{5mm}

\appendix{D) {\bf The STL format }}
Here is top of the file kissing.stl converted using ``admesh" to the human readable 
ASCII format. The entire file has 104'000 lines and contains 14'640 facets. 
The line with ``normal" contains a vector indicating the orientation of the 
triangular facet.  \\

\hspace{-1mm}
\begin{small}
\begin{lstlisting}[language=Mathematica]
solid  Processed by ADMesh version 0.95
  facet normal  2.45300293E-01 -3.88517678E-02  9.68668342E-01
    outer loop
      vertex  1.64594591E-01  0.00000000E+00  9.86361325E-01
      vertex  1.56538755E-01 -5.08625247E-02  9.86361325E-01
      vertex  3.08807552E-01 -1.00337654E-01  9.45817232E-01
    endloop
  endfacet
\end{lstlisting}
\end{small}

\appendix{E) {\bf Mathematica Examples}}

Here are examples of basic ``miniature programs" which can be used to 
produce shapes: 

\parbox{16.8cm}{
\parbox{4.1cm}{\scalebox{0.113}{\includegraphics{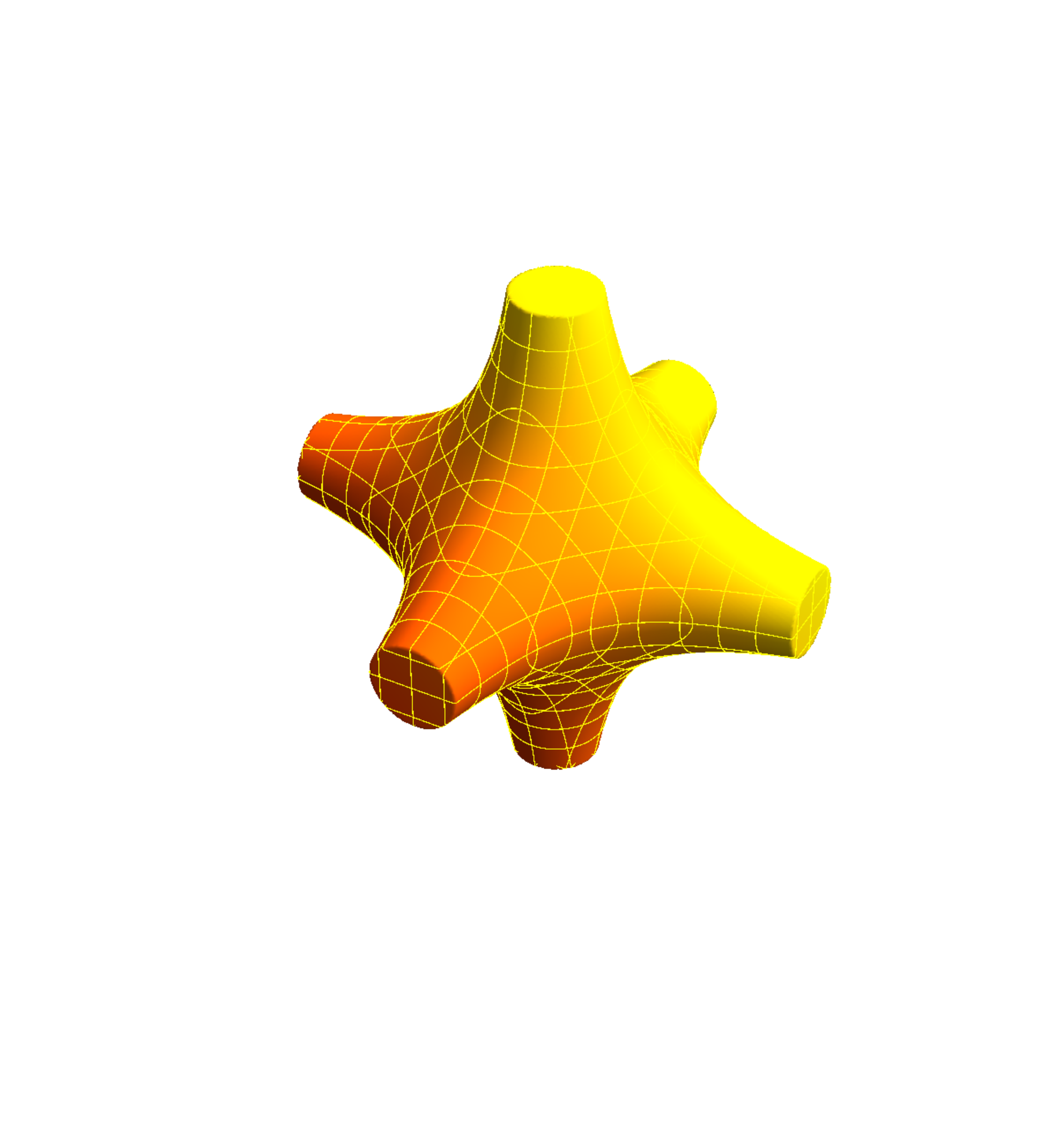}}}
\parbox{4.1cm}{\scalebox{0.113}{\includegraphics{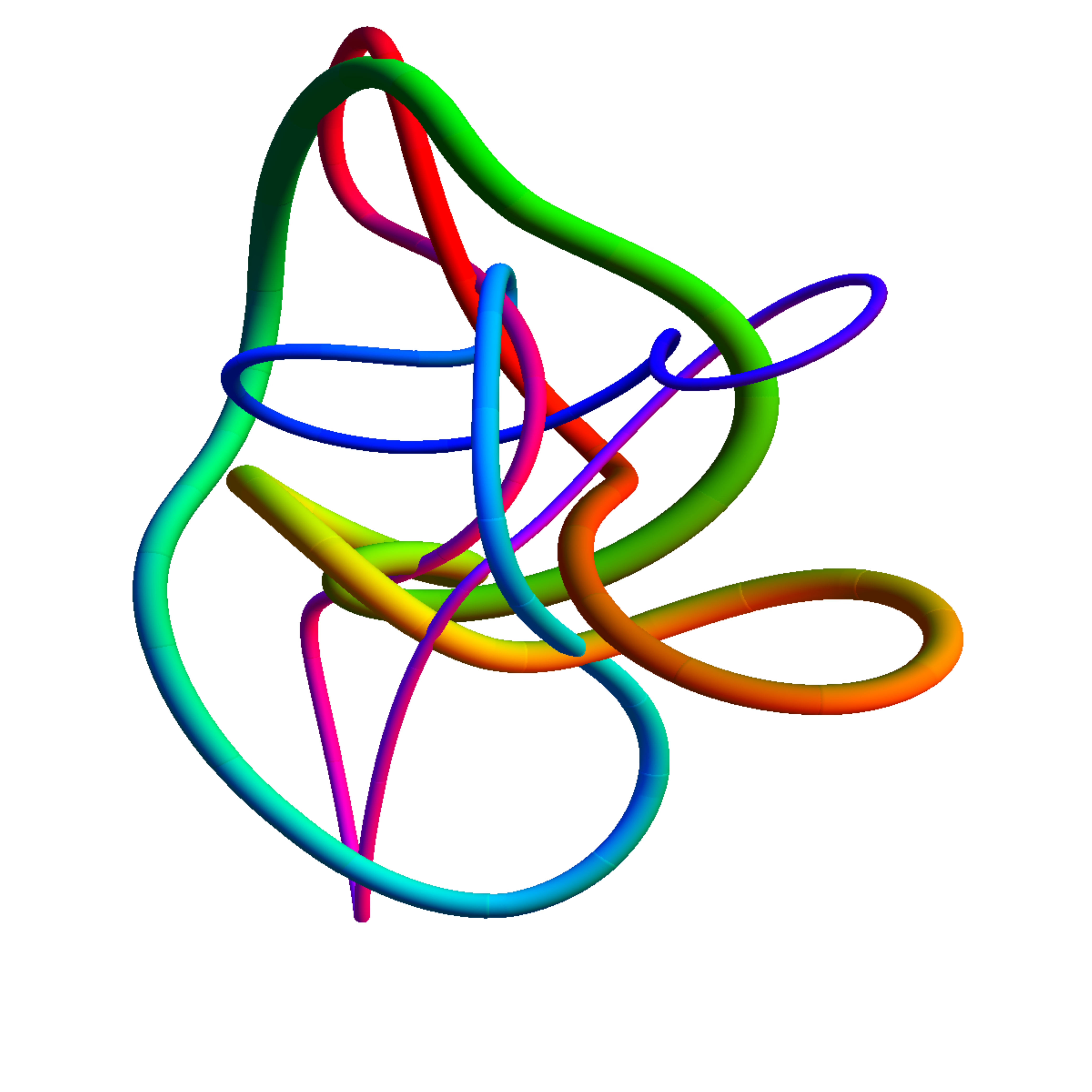}}}
\parbox{4.1cm}{\scalebox{0.113}{\includegraphics{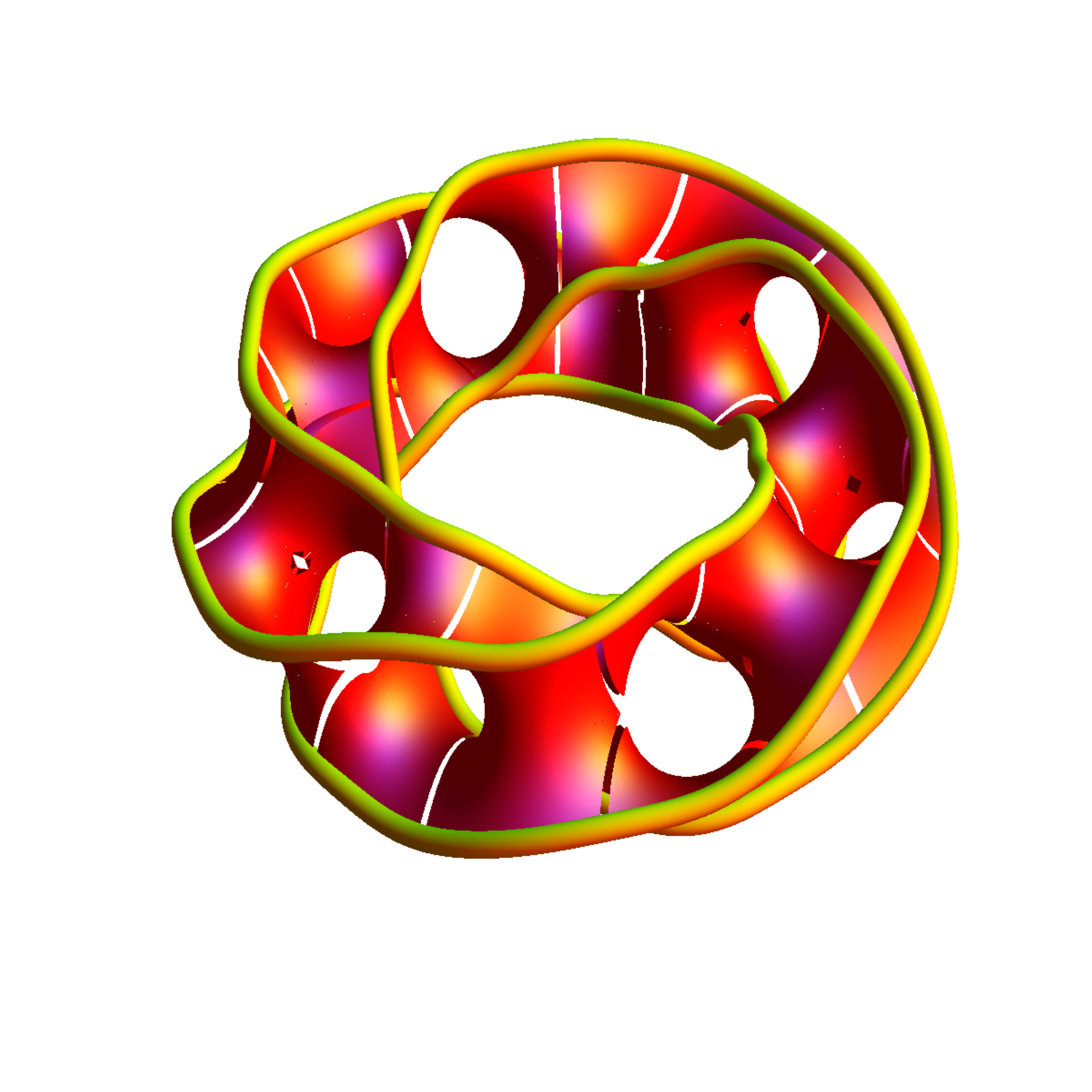}}}
\parbox{4.1cm}{\scalebox{0.113}{\includegraphics{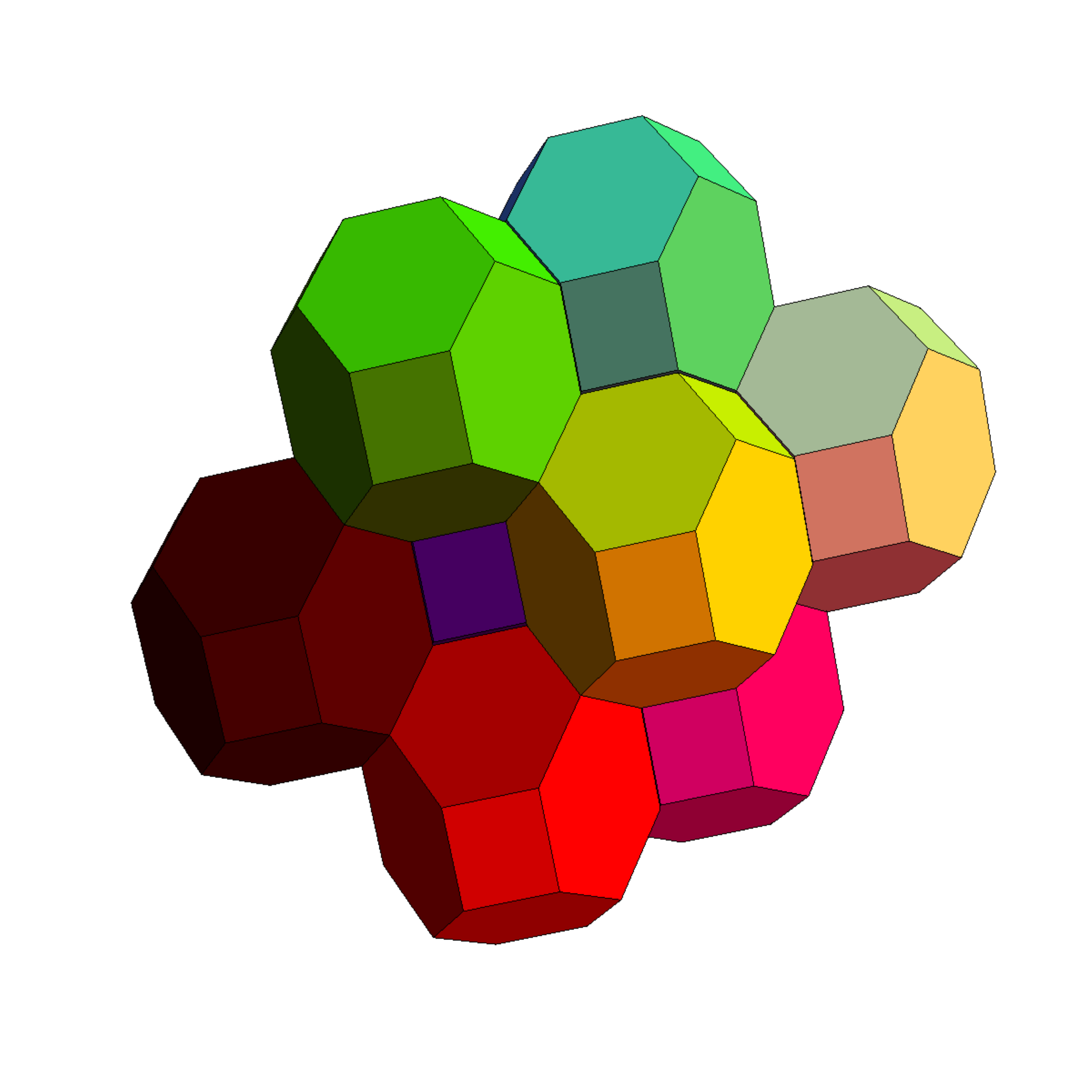}}}
}

\absatz{E1) A region plot:}

\begin{small}
\begin{lstlisting}[language=Mathematica]
RegionPlot3D[
 x^2 y^2+z^2 y^2+x^2 z^2<1 && x^2<6 && y^2<6 && z^2<6, 
 {x,-3,3},{y,-3,3},{z,-3, 3}]
\end{lstlisting}
\end{small}

\absatz{E2) Addition of some "knots"}

\begin{small}
\begin{lstlisting}[language=Mathematica]
u = KnotData[{"PretzelKnot",{3, 5, 2}}, "SpaceCurve"];
v = KnotData[{"TorusKnot", {5, 11}}, "SpaceCurve"];
Graphics3D[Tube[Table[3u[t]-4v[t],{t,0,2 Pi,0.001}],0.3]]
\end{lstlisting}
\end{small}

\absatz{ E3) A Scherk-Collins surface \cite{Collins}} which is close to a minimal surface:

\begin{small}
\begin{lstlisting}[language=Mathematica]
T[{x_,y_,z_},t_]:={x Cos[t]-y*Sin[t],x Sin[t]+y*Cos[t],z};
W[{x_,y_,z_},t_]:={(x+3)Cos[t],(x+3)Sin[t],y}; A=ArcTan; L=Log;
f[z_]:=Module[{p=Sqrt[2Cot[z]],q=Cot[z]+1,r=Re[z]/3},
W[T[Re[{u(L[p-q]-L[p+q])/Sqrt[8],v*I(A[1-p]-A[1+p])/Sqrt[2],z}],r],r]];
Show[Table[ParametricPlot3D[f[x+I*y],{x,0,6Pi},{y,0,6}],{u,-1,1,2},{v,-1,1,2}]]
\end{lstlisting}
\end{small}

\absatz{E4) A polyhedral tessellation} At present, Mathematica commands ``Translate" and ``Rotate" or ``Scale" produce STL files which
are not printable. This requires to take objects apart and put them together again. Here is an example,
which is a visual proof that one can tesselate space with truncated octahedra.

\hspace{-1cm}
\begin{small}
\begin{lstlisting}[language=Mathematica]
T[s_, scale_, trans_] := Module[{P,E},
   P = Table[scale*s[[1,k]] + trans,{k,Length[s[[1]]]}];
   E = s[[2,1]]; Graphics3D[Table[{Polygon[Table[P[[E[[k,l]]]],
   {l,Length[E[[k]]]}]]},{k,Length[E]}]]];
P=PolyhedronData["TruncatedOctahedron","Faces"]; 
Show[Table[T[P,1,{k+l+2m,k-l,3m/2}],{k,-1,1,2},{l,-1,1,2},{m,0,1}]]
\end{lstlisting}
\end{small}

\appendix{F) {\bf Conversion to OpenScad}}

The book \cite{CFZ} as well as a ICTP workshop in Trieste made us aware of OpenScad, a 3D compiler 
which is in spirit very close to a computer algebra system and which also has elements of the 
open source raytracing programming language Povray. An other scripting language for solid modelling is 
``Plasm". The computer algebra system ``Mathematica" has only limited
support for intersecting objects, especially does not yet intersect object given as a mesh. 
Both Povray and OpenScad can do that. Here is a quick 
hack, which allows Mathematica to export graphics to OpenScad.  \\

\hspace{-1cm}
\begin{small}
\begin{lstlisting}[language=Mathematica]
A=PolyhedronData["GreatRhombicosidodecahedron","Faces"]; name="out.scad";
p=10*A[[1]]; t= A[[2,1]]-1; WS=WriteString;
WS[name,"polyhedron("]; WS[name,"points = "]; Write["out.scad",N[p]]; WS[name,","];
WS[name,"triangles = "]; Write["out.scad",t]; WS[name,");"]
Run["cat out.scad |tr \"{\" \"[\"|tr \"}\" \"]\" >tmp; mv tmp out.scad"];
\end{lstlisting}
\end{small}

Once these lines are run through Mathematica the file ``out.scad" can be opened
in OpenScad.   \\

\appendix{G) {\bf OpenScad - CAS comparison  }}

To illustrate intersection of regions we look at the intersection of
three cylinders. 

\absatz{G1) In Mathematica with Region Plot}

\hspace{-1cm}
\begin{small}
\begin{lstlisting}[language=Mathematica]
RegionPlot3D[x^2+y^2<1 && y^2+z^2<1 && x^2+z^2 <1 && x^2 <1 && y^2 <1 
 && z^2<1,{x,-2,2},{y,-2,2},{z,-2,2},PlotPoints->100,Mesh->False]
\end{lstlisting}
\end{small}

The code to parametrize the boundary is much longer and not included here.

\absatz{G2) In OpenScad}

\hspace{-1cm}
\begin{small}
\begin{lstlisting}[language=Mathematica]
module cyl(a){rotate(90,a) cylinder(r=10,h=50,center=true,$fn=200);}  
intersection() { cyl([0,0,0]); cyl([1,0,0]); cyl([0,1,0]); }
\end{lstlisting}
\end{small}

\absatz{G3) In Povray, the open source ray tracer}

\hspace{-1cm}
\begin{small}
\begin{lstlisting}[language=Mathematica]
camera{up y right x location <1,3,-2> look_at <0,0,0>}
light_source { <0,300,-100> color rgb <1,1,1> }
background { rgb<1,1,1> }
#macro r(c) pigment{rgb c} finish {phong 1 ambient 0.5} #end
intersection {
  cylinder{<-1,0,0>,<1,0,0>, 1 texture {r(<1,0,0>)}}
  cylinder{< 0,-1,0>,<0,1,0>,1 texture {r(<1,1,0>)}}
  cylinder{<0,0,-1>,<0,0,1>, 1 texture {r(<0,0,1>)}}
}
\end{lstlisting}
\end{small}

\parbox{16.8cm}{
\parbox{4.1cm}{\scalebox{0.113}{\includegraphics{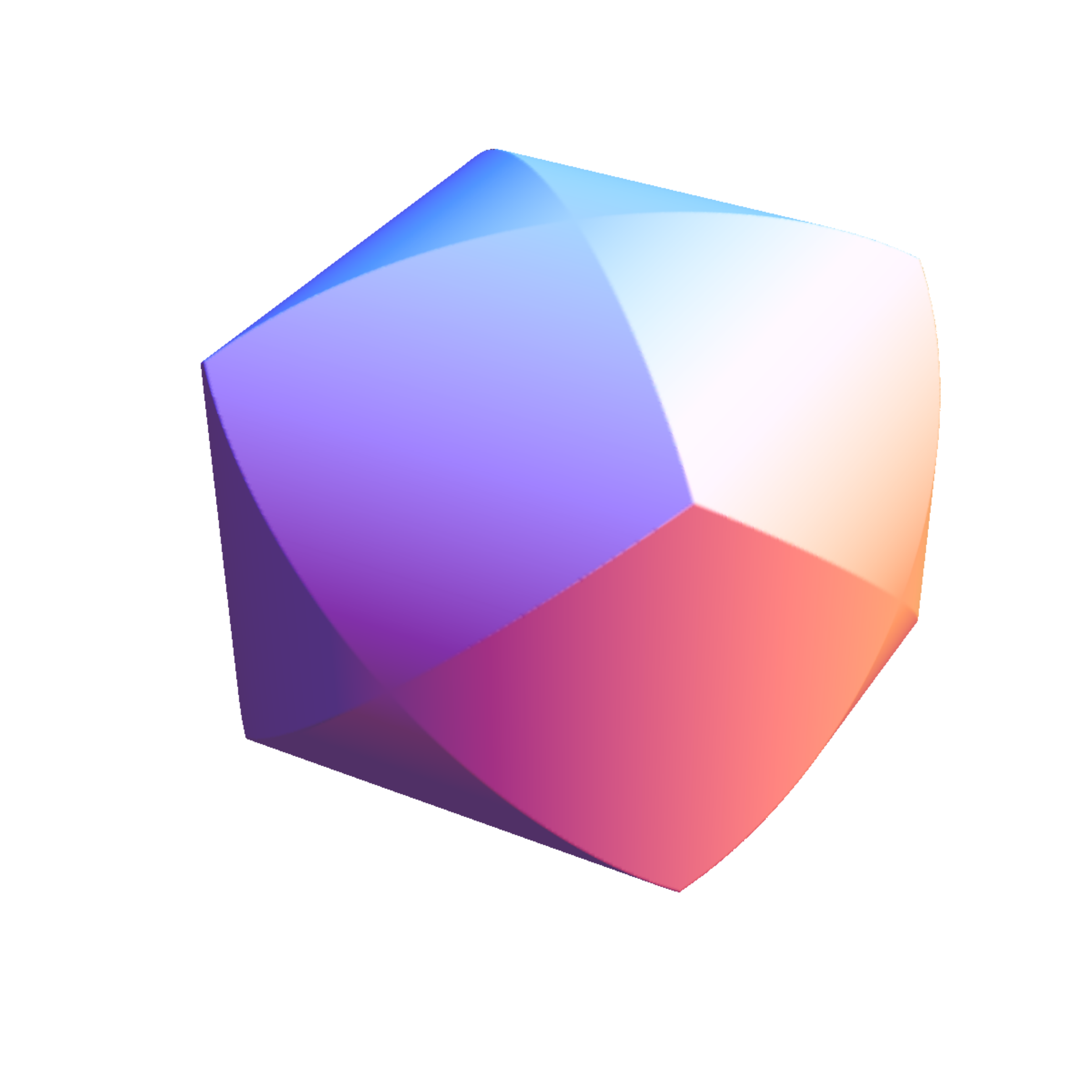}}}
\parbox{4.1cm}{\scalebox{0.113}{\includegraphics{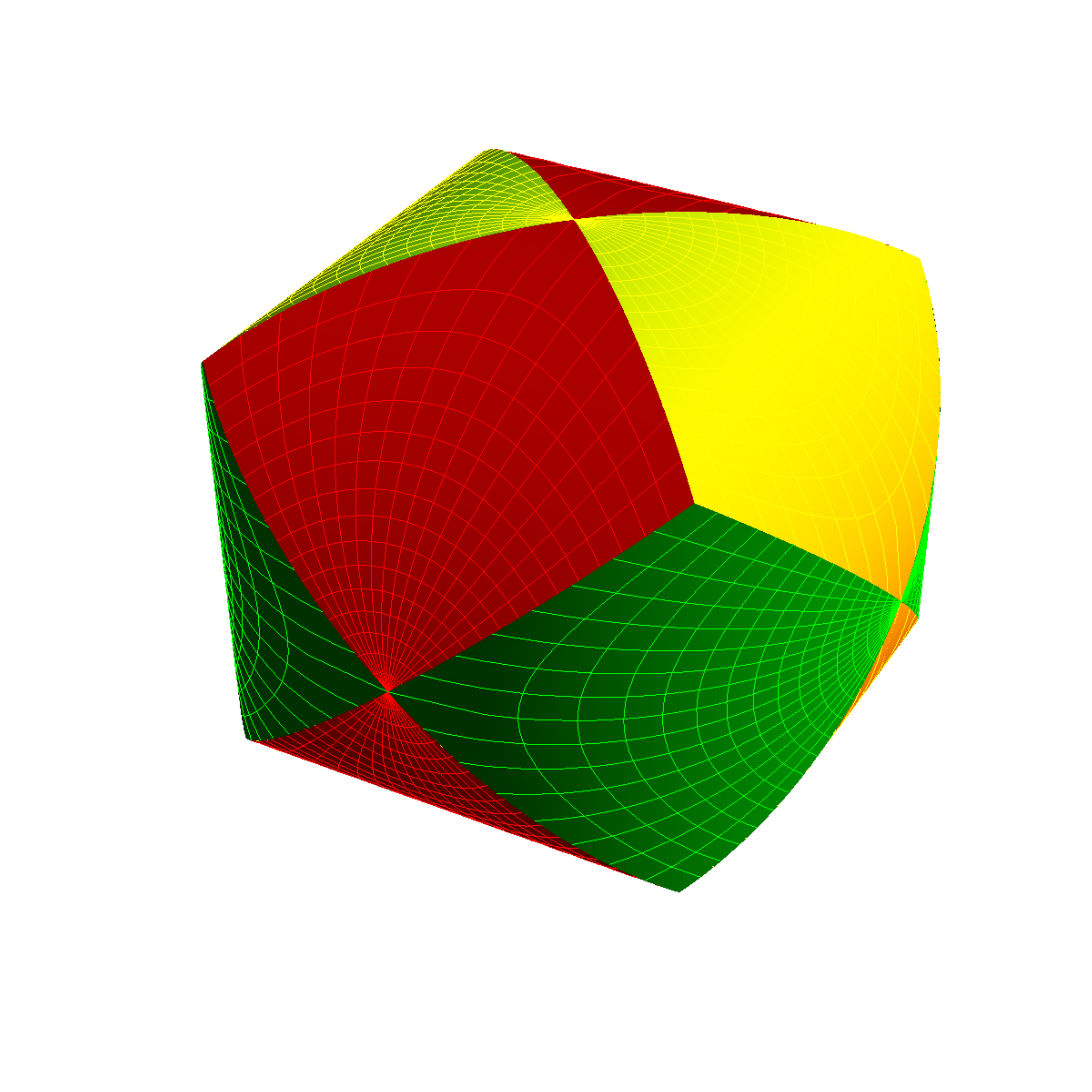}}}
\parbox{4.1cm}{\scalebox{0.113}{\includegraphics{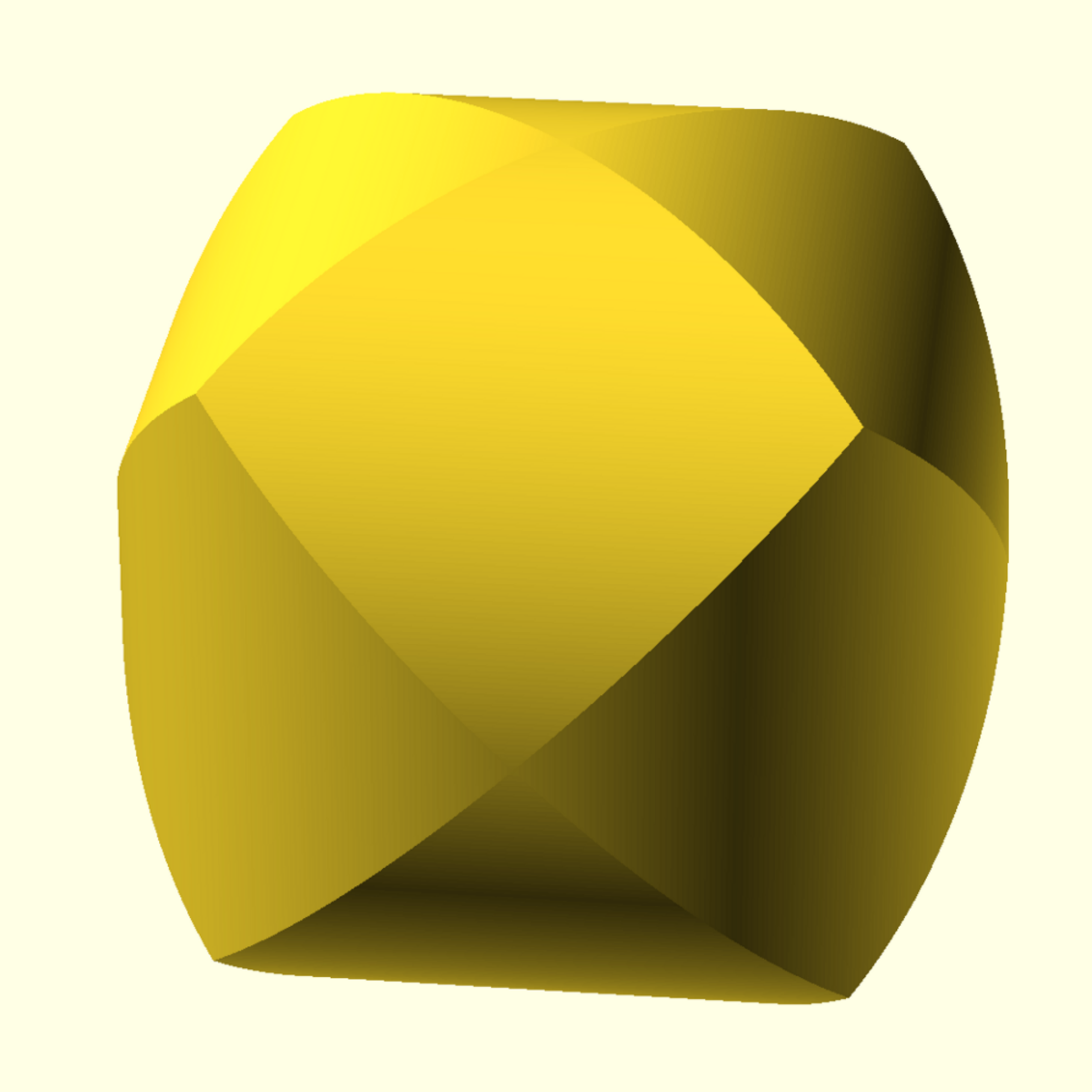}}}
\parbox{4.1cm}{\scalebox{0.113}{\includegraphics{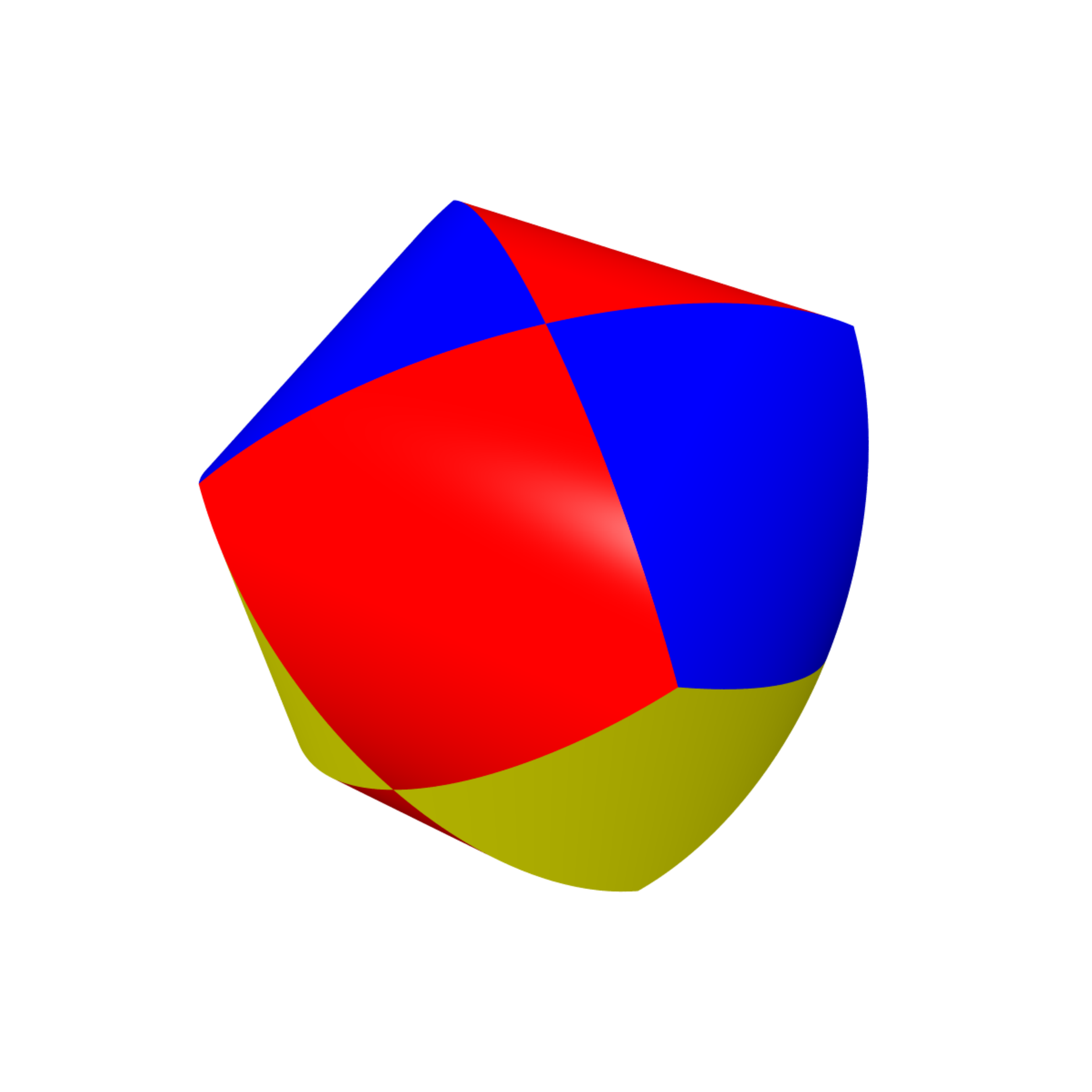}}}
}

\parbox{16.8cm}{
\parbox{4.1cm}{Regionplot Mathematica}
\parbox{4.1cm}{Surface Mathematica}
\parbox{4.1cm}{OpenScad}
\parbox{4.1cm}{Povray}
}

\vspace{3mm}

\appendix{H) {\bf 3D File formats}}
Here are a couple of file formats which can be used for 3D objects. Some can be used for 
3D printing. STL is a format specifically for 3D printing. It was made by ``3D systems",
X3D is the sucessor of WRL, a virtual reality language file format. The format 3DS was introduced
by ``3D Studio", DXF is a CAD data file format by ``Autodesk", 
OBJ used by ``Wavefront technologies", PDF3D by ``Adobe", SCAD by the OpenSCAD, a programmers Modeller, 
PLY by Stanforod graphics lab, POV by the open source ray tracer Povray. Blender, Processing and Meshlab
are all open source projects. 

\begin{center}
\begin{tabular}{|l|l|l|l|l|l|l|}   \hline
{\bf File format}  & Color support & 3D printing    & Mathematica  & Blender  & Processing & Meshlab \\  \hline
STL                & no            & yes            & yes          & yes      & library    &  yes    \\  \hline
WRL                & yes           & yes            & yes          & no       & no         &  yes    \\  \hline
X3D                & yes           & yes            & yes          & no       & no         &  yes    \\  \hline
3DS                & yes           & convert        & yes          & yes      & library    &  yes    \\  \hline
OBJ                & no            & yes            & yes          & yes      & library    &  yes    \\  \hline
PDF3D              & yes           & no             & no           & no       & no         &  no     \\  \hline
SCAD               & no            & convert        & no           & no       & no         &  no     \\  \hline
DXF                & yes           & convert        & yes          & yes      & yes        &  yes    \\  \hline
POV                & yes           & no             & no           & no       & no         &  no     \\  \hline
PLY                & yes           & convert        & yes          & yes      & library    &  yes    \\  \hline
\end{tabular}
\end{center}

\appendix{I) {\bf Tips and Tricks}}
Exporting files from Mathematica:

\hspace{-1cm}
\begin{small}
\begin{lstlisting}[language=Mathematica]
Export["file.stl",Graphics3D[{Red,Sphere[{0,0,0},1]}]];
Export["file.x3d",Graphics3D[{Red,Sphere[{0,0,0},1]}]];
Export["file.wrl",Graphics3D[{Red,Sphere[{0,0,0},1]}]];
Export["file.3ds",Graphics3D[{Red,Sphere[{0,0,0},1]}]];
Export["file.obj",Graphics3D[{Red,Sphere[{0,0,0},1]}]];
\end{lstlisting}
\end{small}

3DS files can be imported in ``Sketchup", ``Blender" or ``Google earth". X3D is the 
successor of WRL and allows to be 3D printed in color. 
The formats OBJ and STL do not support color as for now. 
Here are examples of command-line conversions with admesh. The first 
converts an  STL file into ASCII, the second fills holes. 

\hspace{-1cm}
\begin{small}
\begin{lstlisting}[language=Mathematica]
admesh -a out.stl in.stl
admesh -f  in.stl
\end{lstlisting}
\end{small}

And here is an example to use the tool to convert a height PNG height field
to a printable STL file: 

\hspace{-1cm}
\begin{small}
\begin{lstlisting}[language=Mathematica]
png23d -t x -f surface -o stl -l 40 -d 30 -w 100 h.png h.stl
\end{lstlisting}
\end{small}

The next figure shows a relief of the Rheinfalls region in Switzerland, converted
into a printable STL file. The gray shaded hight relief as well as the color coded
third picture had been colorized by hand in the late 1990ies by the first author 
from scanned 1:25'000 maps of the region. Digital elevation data had been very expensive
at that time. Today, services like google earth brought down the prizes. 
The last picture shows a hand colored file of the area with the rhein river near the town of Neuhausen. \\

\parbox{16.8cm}{
\parbox{5cm}{\scalebox{0.12}{\includegraphics{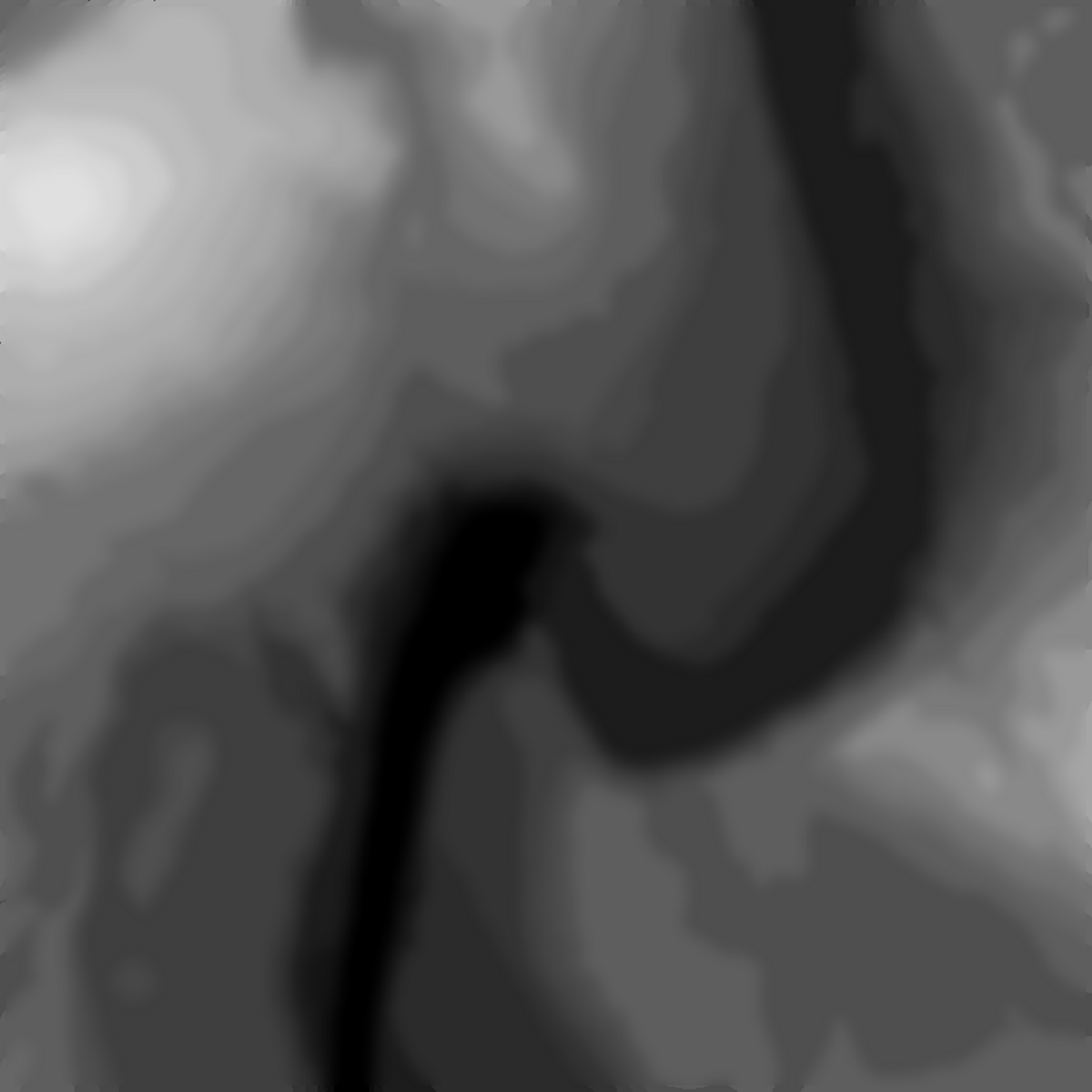}}}
\parbox{5cm}{\scalebox{0.12}{\includegraphics{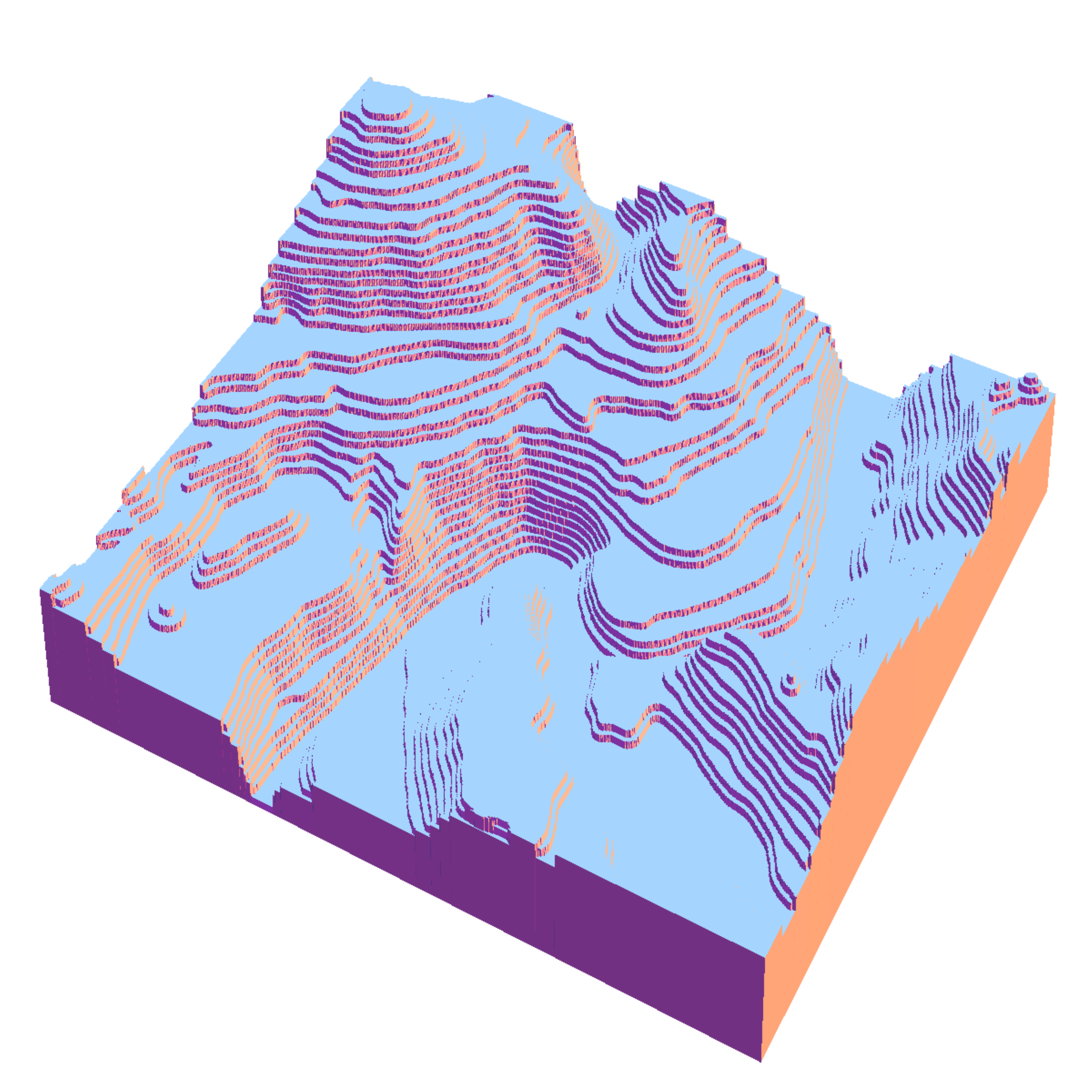}}}
\parbox{5cm}{\scalebox{0.12}{\includegraphics{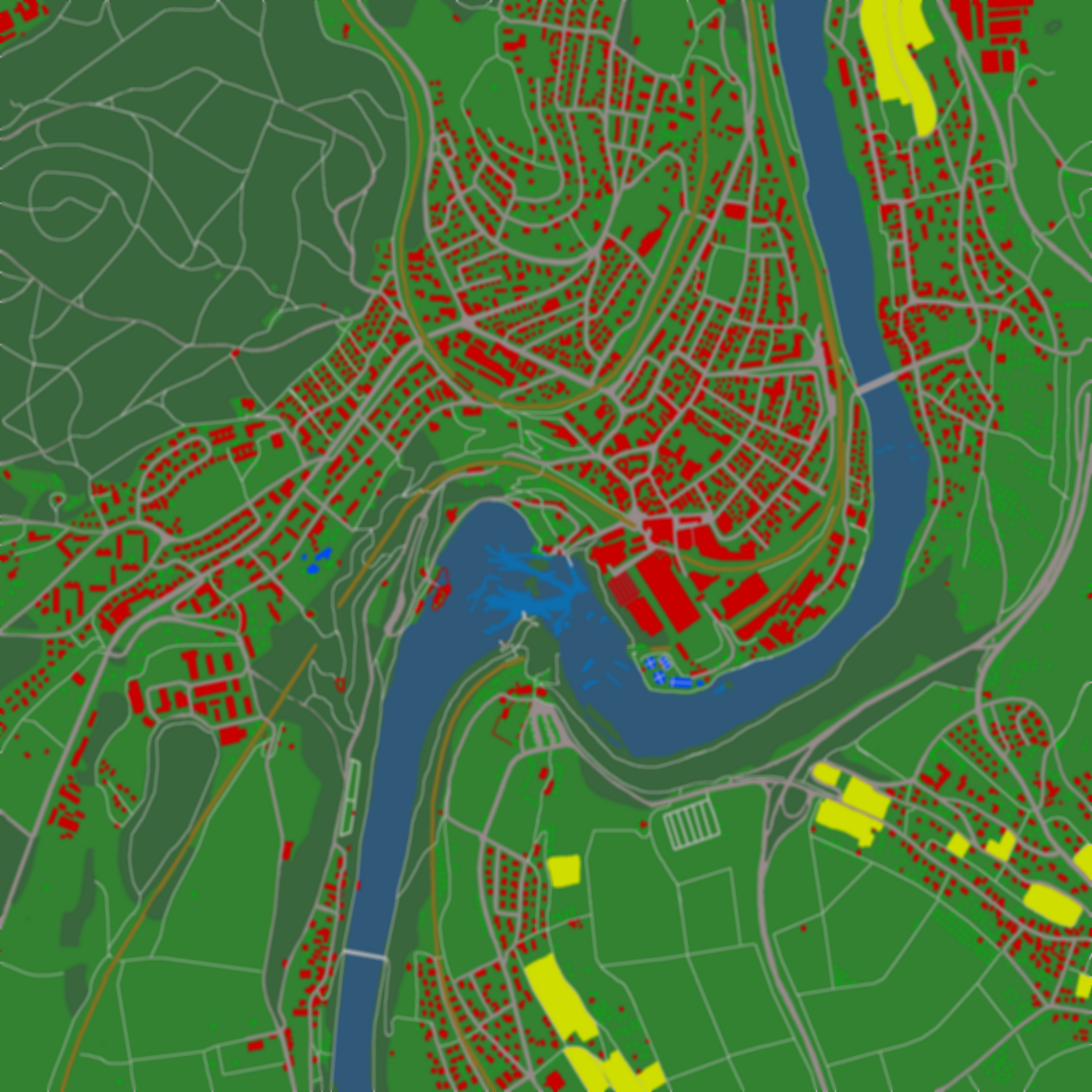}}}
}

\vspace{4mm}

How can one work with STL files in Mathematica?  Here is an example: since the Gomboc specifications are not
public, folks had to reverse engineer it. We grabbed the STL file on thingiverse. Now, one has a graphics
data structure A. With A[[1,2,1]], one can access the polygon tensor which consists of a list of triangles
in space. Now, it is possible to work with this object. As an example, we add at every center of a triangle
a small sphere: 

\begin{small}
\begin{lstlisting}[language=Mathematica]
A = Import["gomboc.stl"];    B=A[[1,2,1]];
S[{a_,b_,c_}]:={Yellow,Sphere[(a+b+c)/3,0.001]};
U={Black,Map[Polygon,B]}; V=Map[S,B]; 
Graphics3D[{U,V},Boxed->False];
\end{lstlisting}
\end{small}

\parbox{16.8cm}{
\parbox{5cm}{\scalebox{0.12}{\includegraphics{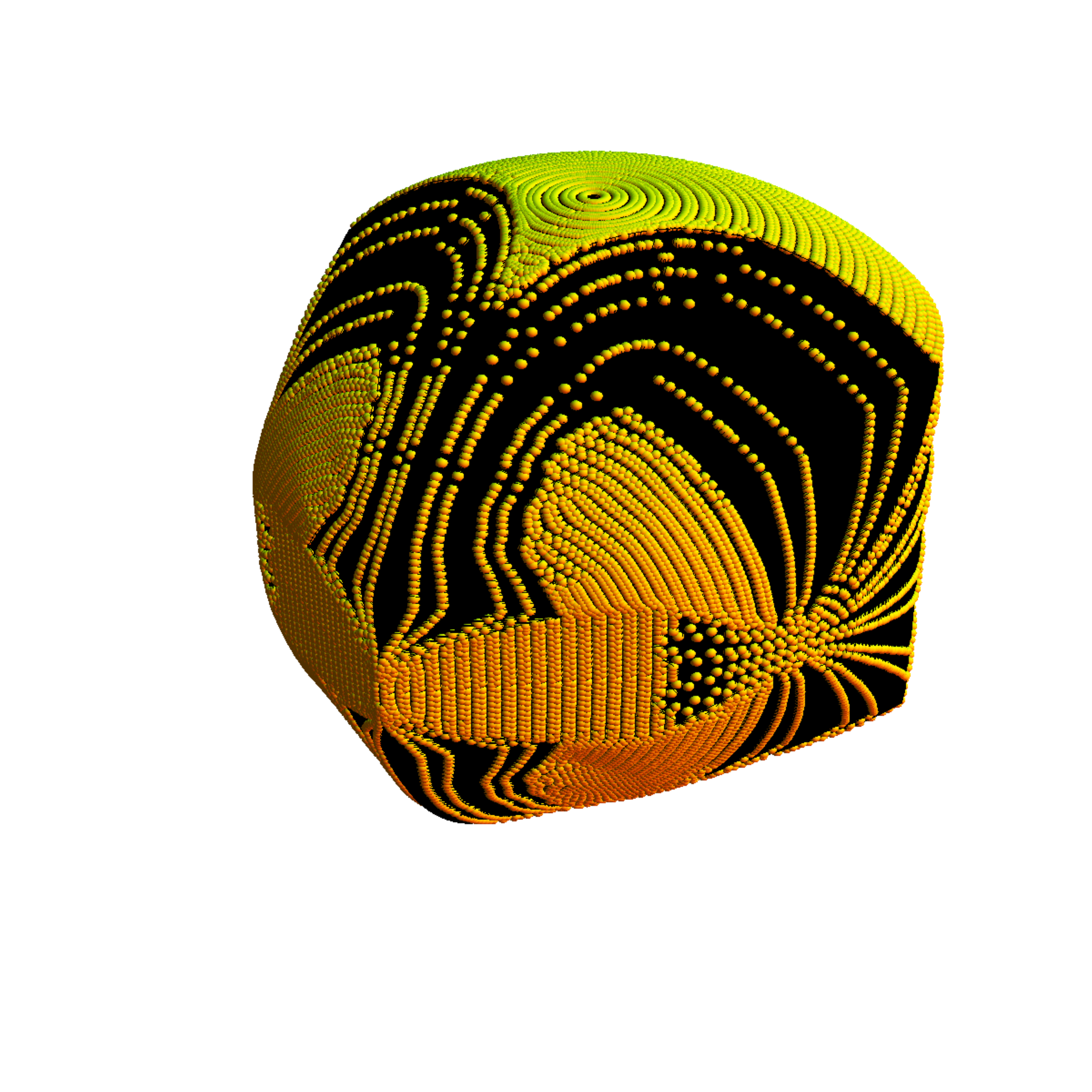}}}
}

\appendix{J) {\bf Open problems}}

Here is a list of 10 open mathematical problems which were mentioned en passant 
while illustrating some of the mathematics \\

\begin{center}
\begin{tabular}{|l|l|l|l|}   \hline
{\bf Problem}      &  {\bf Remarks}   &  {\bf References}         \\  \hline
Kissing problem    &  Unknown in dim $\geq 5$  &  \cite{Pfender}  \\  \hline
Perfect cuboid     &  Problem by Euler         &  \cite{Guy}      \\  \hline
Riemann hypothesis &  Problem by Riemann       &  \cite{BorweinRiemann}   \\  \hline
Chaos in 3D billiards & Positive entropy?      &  \cite{Zaslavsky}        \\  \hline
Local connectivity   &  for Mandelbrot set     &  \cite{princetonguide}   \\  \hline 
Random matrices      &  Almost periodic entries & \cite{battle}           \\  \hline
Isospectral drums    &  Convex examples?       &  \cite{Kac66}            \\  \hline 
Falling coin problem  & Positive entropy?      &  \cite{SGSPK}            \\  \hline
Normality of $\pi$    &  Unknown in any base   &  \cite{BBCDDY,Artacho}   \\  \hline
Last statement of Jacobi & 4 cusps for non-umbillical points  & \cite{Sinclair,SinclairTanaka} \\ \hline 
Mandelbulb            & Connectivity properties  &  \cite{White} \\ \hline
\end{tabular}
\end{center}

\appendix{K) {\bf Scanning}}

Software deveoloping tools like ``Kinect Fusion" from Microsoft, 
``KScan3D" or ``Artec Studio" use the affordable ``Microsoft Kinect" device to scan objects.
The Kinect uses a PrimeSense 3D sensing technology, where an IR light source projects light
onto a scene which is then captured back by a CMOS sensor. This produces a debth image
where pixel brightness corresponds to distance to the camera. 
Commercial tools like ``KScan3D" or ``Artec Studio" can produce textured objects. 
During the workshop, we were shown how to use the ``Kinect" together with 
the opensource software ``Processing" to scan objects to be 3D printed \cite{Borenstein}. 
This is a low cost approach which allows to produces STL files from real objects but there
is a long way to get the quality of cutting edge products like Artec, which work
well and do not need any setup. All these approaches produce triangulated meshes.

\parbox{16.8cm}{ \parbox{10cm}{
\parbox{4cm}{\scalebox{0.12}{\includegraphics{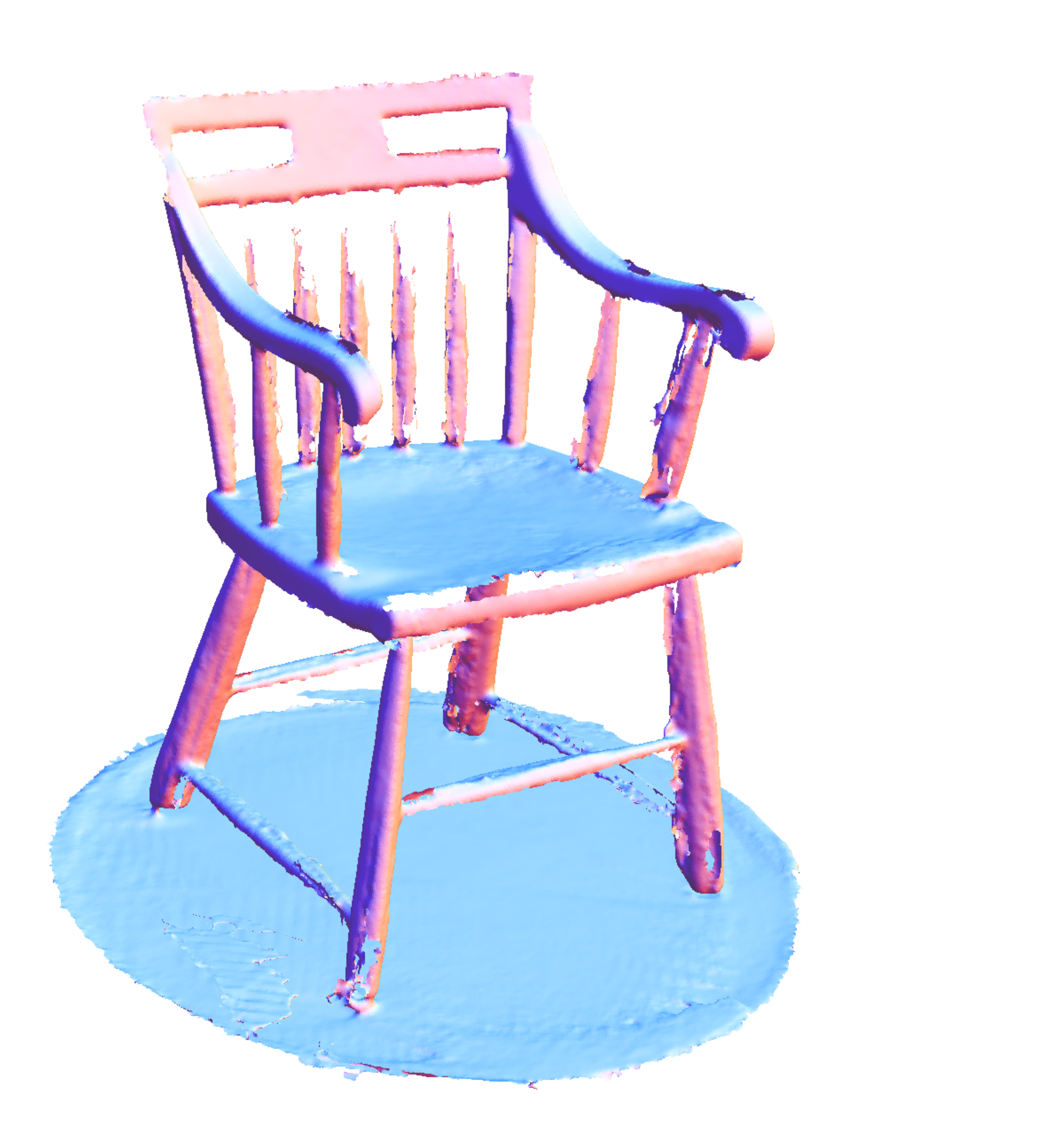}}}
\parbox{5.5cm}{\scalebox{0.19}{\includegraphics{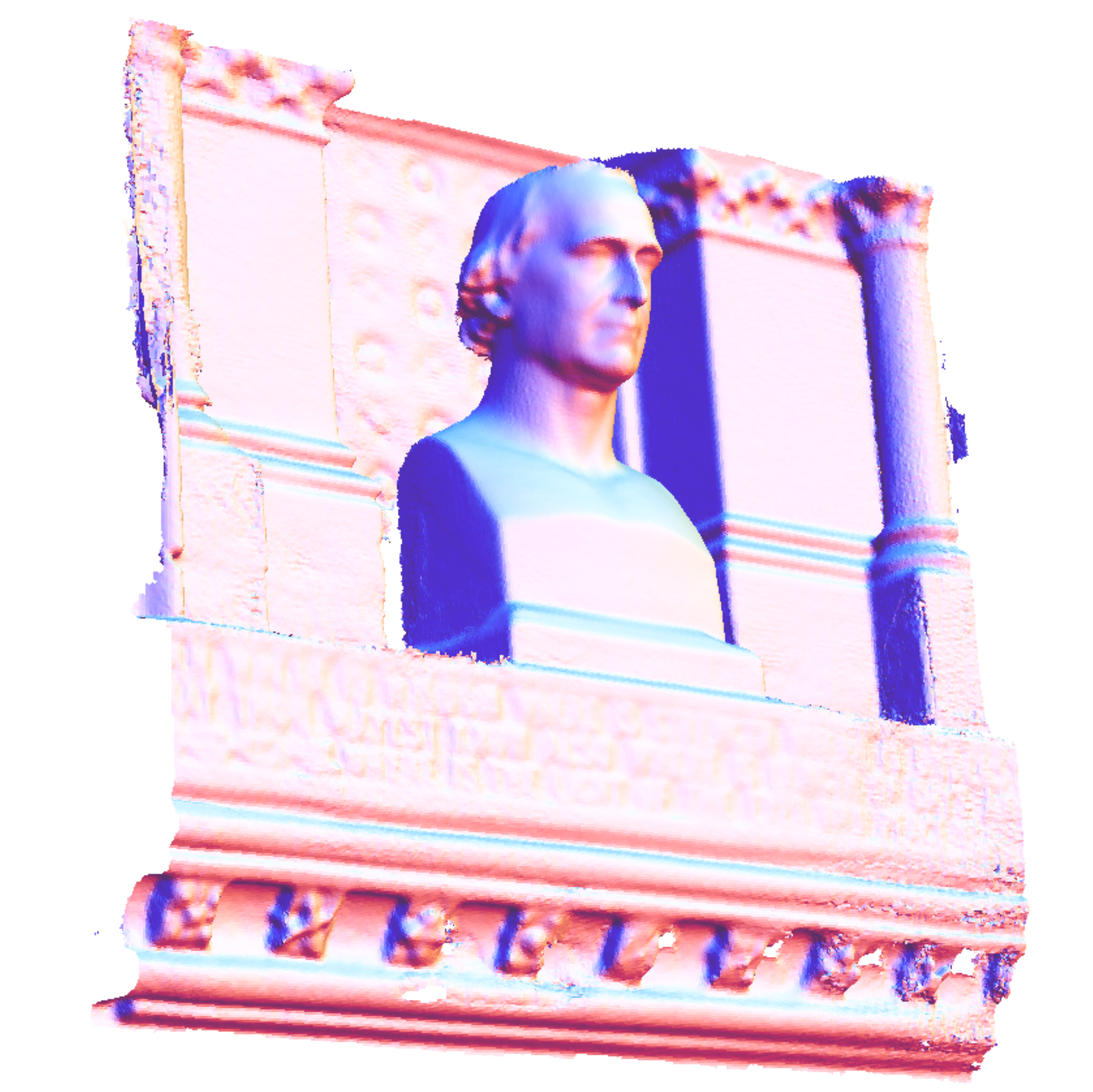}}}
} \parbox{6cm}{
To the left a scan of a Harvard office chair 
using ``Artec Studio 9.1" and kinect. 
The scanner software running on a macbook air handled over 10 frames
per second. The Artec software combined things to a 3D object, 
cleaned the data and exported the STL which was then imported into Mathematica. 
The right picture shows the face of James Walker, president of Harvard
College, 1863-1850. It got scanned at Harvard Memorial hall with the same setup. 
}}

A much needed, but difficult problem is be to get mathematical OpenScad style descriptions of the objects. 
Similar than optical character recognition or fractal image compression it is the task to get from 
high entropy formats to low entropy mathematical rules describing the object. 
All such tasks need a large amount or artificial intelligence.
In the future, an advanced system of this type would be able to get from 
scanned pictures of the ``Antikythera" a  CAD blueprint of a machine which would work the same way 
as the original ancient computer did. \\

\appendix{L) {\bf More source examples}}

Since this file will be on ArXiv, where the LaTeX source can be inspected, you can copy paste the 
source from there:

\absatz{L1) Calabi Yau surface}

The following code has been used to generate Youtube animations of Calabi-Yau manifolds
and was also used in the workshop as a demonstration: 

\hspace{-1cm}
\begin{small}
\begin{lstlisting}[language=Mathematica]
G[alpha_]:=Module[{},n = 5; R=10; 
CalabiYau[z_, k1_, k2_] := Module[{
    z1 = Exp[2Pi I k1/n]Cosh[z]^(2/n),
    z2 = Exp[2Pi I k2/n]Sinh[z]^(2/n)},
    N[{Re[z1],Re[z2],Cos[alpha]Im[z1]+Sin[alpha]Im[z2]}]];
F[k1_,k2_]:=Module[{}, XX=CalabiYau[u+I v,k1,k2];
S1=Graphics3D[Table[{Hue[Abs[v]/2],Tube[Table[XX,{u,-1,1,2/R}]]},{v,0,Pi/2,Pi/20}]];
S2=Graphics3D[Table[{Hue[Abs[u]/2],Tube[Table[XX,{v,0,Pi/2,(Pi/2)/R}]]},{u,-1,1,0.1}]]; 
Show[{S1,S2}]];
Show[Table[F[k1,k2],{k1, 0, n-1}, {k2,0,n-1}]]];
G[0.2]
\end{lstlisting}
\end{small}

The surface is featured on the cover of \cite{Happening1994}.

\absatz{L2) Escher Stairs}

This code was used to produce an object placed in google earth so that if one looks
at the object from the right place, the Escher effect appears. We mention this example because
it illustrates how one can also build things in Mathematica ``layer by layer". It was much
easier for us to build the individual ``floors" of the building using matrices and then place
bricks where the matrices have entries $1$. 

\hspace{-1cm}
\begin{small}
\begin{lstlisting}[language=Mathematica]
A0={{1,0,0,0,0,0,0},{1,0,0,0,0,0,0},{1,0,0,0,0,0,0},
    {1,0,0,0,0,0,0},{1,0,0,0,0,0,0},{1,1,1,1,1,0,0}};
A1={{0,0,0,0,0,0,1},{0,0,0,0,0,0,1},{0,0,0,0,0,0,1},
    {1,0,0,0,0,0,1},{1,0,0,0,0,0,1},{1,1,1,1,1,1,1}};
A2={{0,0,0,0,0,0,1},{0,0,0,0,0,0,1},{0,0,0,0,0,0,1},
    {0,0,0,0,0,0,1},{0,0,0,0,0,0,1},{0,0,0,0,1,1,1}};
A3={{0,0,0,0,1,1,1},{0,0,0,0,0,0,1},{0,0,0,0,0,0,1},
    {0,0,0,0,0,0,0},{0,0,0,0,0,0,0},{0,0,0,0,0,0,0}}; 
cuboid[x_]:={{Blue,  Polygon[{x+{0,0,0},x+{1,0,0},x+{1,1,0},x+{0,1,0}}]},
             {Green, Polygon[{x+{0,0,0},x+{0,1,0},x+{0,1,1},x+{0,0,1}}]},
             {Orange,Polygon[{x+{0,0,0},x+{1,0,0},x+{1,0,1},x+{0,0,1}}]},
             {Yellow,Polygon[{x+{1,1,1},x+{0,1,1},x+{0,0,1},x+{1,0,1}}]},
             {Red,  ,Polygon[{x+{1,1,1},x+{1,0,1},x+{1,0,0},x+{1,1,0}}]},
             {Pink  ,Polygon[{x+{1,1,1},x+{0,1,1},x+{0,1,0},x+{1,1,0}}]}};
S0=Table[If[A0[[i,j]]==1,cuboid[{i,j,-1}],{}],{i,Length[A0]},{j,Length[A0[[1]]]}];
S1=Table[If[A1[[i,j]]==1,cuboid[{i,j,0}],{}],{i,Length[A1]},{j,Length[A1[[1]]]}];
S2=Table[If[A2[[i,j]]==1,cuboid[{i,j,1}],{}],{i,Length[A2]},{j,Length[A2[[1]]]}];
S3=Table[If[A3[[i,j]]==1,cuboid[{i,j,2}],{}],{i,Length[A3]},{j,Length[A3[[1]]]}];
S=Graphics3D[{Red,EdgeForm[{Red,Opacity[0],Thickness[0.00001]}],
{S0,S1,S2,S3}}, Boxed->False, ViewPoint -> {-0.407193, 2.2753, 2.47128}, 
 ViewVertical -> {0.101614, -0.171295, 1.71738}]
\end{lstlisting}
\end{small}

\absatz{L3) Drinkable proof}

\hspace{-1cm}
\begin{small}
\begin{lstlisting}[language=Mathematica]
P=ParametricPlot3D; 
M=6;a=0.1;b=0.06;b1=Sin[b];b2=Cos[b];b1a=ArcSin[b1/(1+a)];b2a=(1+a)Cos[b1a];
bowl1 =P[{Cos[t] Sin[s],Sin[t] Sin[s],Cos[s]}+{0,0,2},{t,0,2Pi},{s,Pi/2,Pi-b}];
bowl2 =P[(1+a){Cos[t] Sin[s],Sin[t] Sin[s],Cos[s]}+{0,0,2},{t,0,2Pi},{s,Pi/2,Pi-b1a}];
drain =P[{b1 Cos[t],b1 Sin[t],z},{z,2-b2,2-b2a},{t,0,2Pi}];
cone1 =P[{(1-z)Cos[t],(1-z)Sin[t],z-1},{t,0,2Pi},{z,0,1}];
cone2 =P[{(1+a)(1-z) Cos[t],(1+a)(1-z) Sin[t],z-1},{t,0,2Pi},{z,0,1}];
wall1 =P[{Cos[t],Sin[t],z},{t,0,2Pi},{z,-1,0},PlotStyle->{Blue}];
wall2 =P[{(1+a)Cos[t],(1+a)Sin[t],z},{t,0,2Pi},{z,-1,0}];
top1  =P[{r Cos[t],r Sin[t],0},{r,1,1+a},{t,0,2Pi}];
top2  =P[{r Cos[t],r Sin[t],2},{r,1,1+a},{t,0,2Pi}];
connec=Graphics3D[{Table[m=2Pi*k/M;a2=a/2;{Blue,Cylinder[{{(1+a2)Cos[m],(1+a2)Sin[m],-1},
       {(1+a2)Cos[m],(1+a2)Sin[m],2}},a2]},{k,M}]}];
bottom=P[{r Cos[t],r Sin[t],-1},{r,0, 1+a},{t,0,2Pi}];
S=Show[{connec,drain,bowl1,bowl2,cone1,cone2,wall1,wall2,bottom,top1,top2}]
\end{lstlisting}
\end{small}

\absatz{L4) Lorenz Attractor}

This is the computation of the Lorentz attractor using Runge Kutta. Even so Mathematica
has sophisticated numerical differential solvers built in, doing things by hand has advantages
since it avoids the "black box" of hidden library code. The following implentation 
only uses basic arithmetic and can be translated easily into any programming language. 
Even within Mathematica, it can have advantages to go back to the basics. For Hamiltonian systems
for example, one can so assure that the evolution stays on the energy surface, or one wants to 
enforce that the code respects certain symmetries. An other reason is when wanting to implement 
motion on manifolds like a torus, where original DSolve routines would produce discontinuous
solution curves. We have even used code as given below to solve 
ordinary differential equations numerically within the open source raytracer ``Povray".
We see that the Runge-Kutta integration step can be implemented in one line. 

\hspace{-1cm}
\begin{small}
\begin{lstlisting}[language=Mathematica]
RK[f_,x_,s_]:=Module[{a,b,c,u,v,w,q},u=s*f[x];
 a=x+u/2;v=s*f[a];b=x+v/2;w=s*f[b];c=x+w;q=s*f[c];x+(u+2v+2w+q)/6];
f[{x_,y_,z_}]:={10 (y - x), -x z + 28 x - y, x y -8/3 z};
T[X_]:= RK[f,X,0.01];s=NestList[T,{0,1,0},1000]; 
S1=Graphics3D[{Yellow,Tube[s]}];
a=10; b=-7; S2=Graphics3D[{Red,Cuboid[{{-a,-a,b-2},{a,a,b}}]}];
S3=Graphics3D[{Red,Cylinder[{s[[1]],s[[1]]+{0,0,b-s[[1,3]]}},1]}];
S=Show[{S1,S2,S3},AspectRatio->1,Boxed->False]; 
\end{lstlisting}
\end{small}

\absatz{L5) Feigenbaum bifurcation}

The following example produces a Feigenbaum bifurcation diagram in space. 
We took a circle map $T_c(x) = c \sin(\pi x)$ and vary $c$ from $0.55$ to $1$.
The $c$ corresponds to the height. Depending on $c$, we see the attractor
at that specific value. After a cascade of period doubling bifurcations, the 
Feigenbaum attractor appears. The color is the Lyapunov exponent. Red and orange
colors show negative Lyapunov exponent which means stability. Blue on the top
is chaos. 

\hspace{-1cm}
\begin{small}
\begin{lstlisting}[language=Mathematica]
g=Compile[{c},({c,#} &) /@ Union[Drop[NestList[c*Sin[Pi #]&,0.3,3500],400]]];
L[g_]:=Module[{cc=g[[1,1]]}, Sum[Log[Abs[cc*Pi*Cos[Pi*g[[i,2]]]]],{i,Length[g]}]/Length[g]];
P[x_]:={Sphere[{Cos[2Pi x[[2]]] x[[1]], Sin[2Pi x[[2]]] x[[1]],5 x[[1]]},0.010]};
Graphics3D[Table[{Hue[Max[0,lya[g[c]]]],Map[P,g[c]]},{c,0.55,1,0.001}]];
\end{lstlisting}
\end{small}

\appendix{M) {\bf Math education time lines}}
The folowing slide was shown first at the ICTM Boston \cite{KnillICTM}. It has been 
updated since: new technologies like Wolfram Alpha, HTML5 or WebGL have been added. 
The time line illustrates how much has changed in the classroom in the last 40 years. 
The book \cite{Toolsteaching} documents this in-debth from other angles too. 

\begin{center}
\scalebox{0.29}{\includegraphics{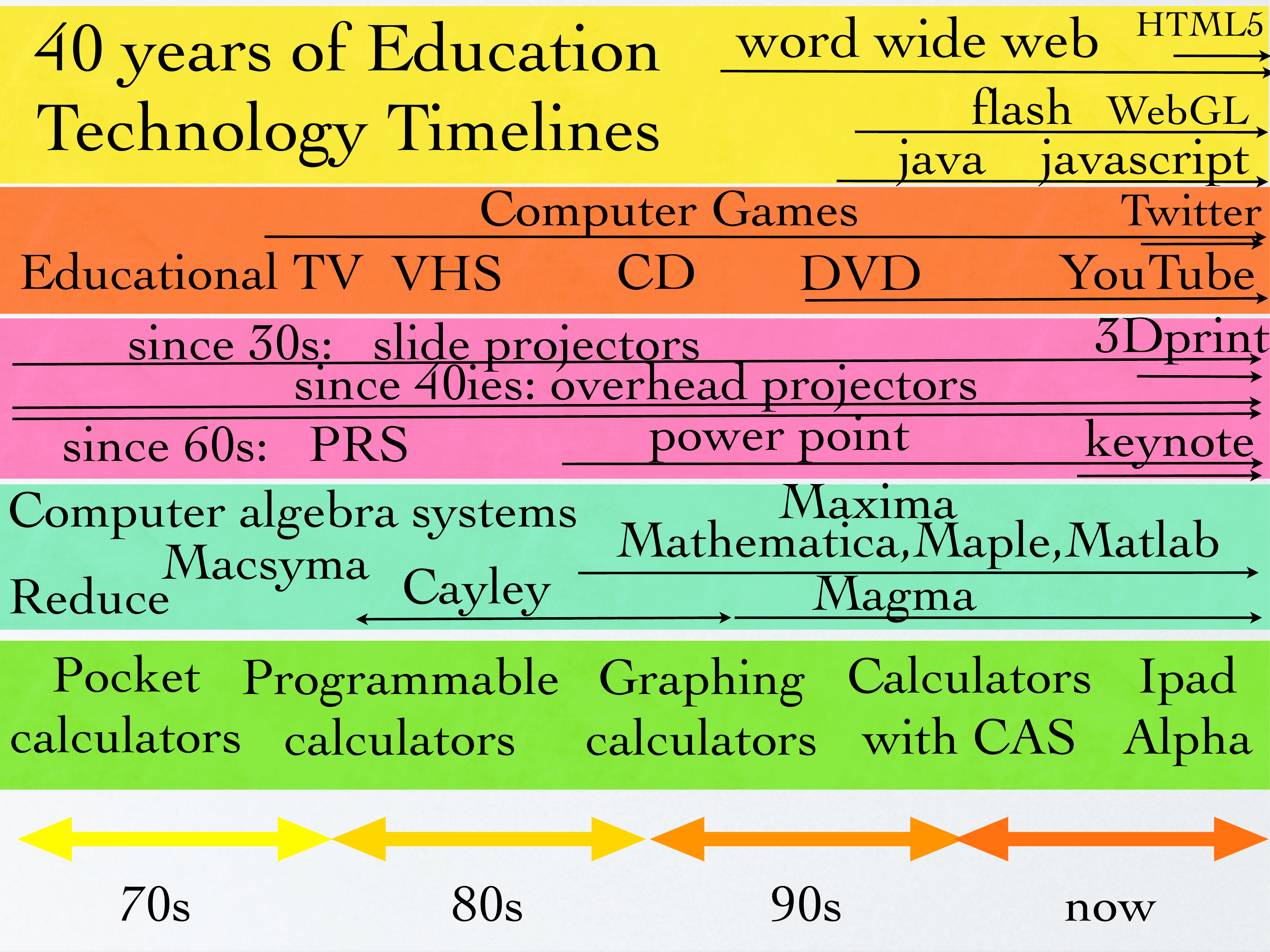}} 
\end{center}

\vspace{-4mm}
We witnessed 3 revolutions: {\bf ``new math"} brought more advanced mathematics to the classroom of
the {\bf generation X}. Calculators and new tools like the overhead projector or Xerox copying
tools amplified these changes. 20 years later came the {\bf ``math wars"}. It effected the generation Y
and also contained both content and philosophy shifts as well as technological earthquakes like 
widespread use of {\bf computer algebra systems} and the {\bf world wide web}. Again 20 years later, now for the
{\bf generation Z}, we see {\bf social media} and 
{\bf massive open online courses} changing the education landscape. We have currently no idea yet, 
where this is going, but when looking back to the other revolutions, it is likely that the soup is again 
eaten colder than cooked: social media start to show their limitations and the drop out rates for 
MOOCs courses can be enormous. Revolutions are often closer to convolutions \cite{Grattan-Guinness}.
Fortunately, future generations of teachers and students can just pick from a 
larger and larger menu and tools. What worked will survive.  But some things have not changed over thousands
of years: an example is a close student-teacher interaction. Its a bit like with 
{\bf Apollonian cones} which were used by the Greeks already in the classroom. They still are around. What has changed
is that they can now be printed fast and cheaply. For a student, creating and experimenting with a real model can make as 
big as an impression as the most fancy 3D visualization, seen with the latest virtual reality gaming gear.

\appendix{N) {\bf Revolutions}}
Here is a  simplified caricature illustrating
change in technology and information. This was done for a talk \cite{KnillBoston}: 

\begin{tabular}{|c|c|c|}   \hline
{\bf Industrial}  &
{\bf Information}    &
{\bf Perception}  \\  \hline
\begin{tabular}{llr}
"Steam"  & Textile        &   1775 \\
"Steel"  & Chemical       &   1850 \\
"Space"  & Exploration    &   1969
\end{tabular} 
&
\begin{tabular}{lll}
"Print"   & Press                &   1439 \\
"Talk"    & Telegraph            &   1880 \\ 
"Connect" & Web                  &   1970 
\end{tabular} 
&
\begin{tabular}{lll}
"Macro"   & Telescope            &   1608 \\
"Micro"   & Microscope           &   1644 \\
"Meso"    & 3DScan               &   1960 
\end{tabular} \\ \hline
\end{tabular}

\vspace{1mm}

\begin{tabular}{|c|c|c|}   \hline
{\bf Media}  &
{\bf Storage}    &
{\bf Education}  \\  \hline
\begin{tabular}{lll}
"Picture"  &   Photo              &   1827 \\
"Movie"    &   Film               &   1878 \\
"Objects"  &   3DPrint            &   1970
\end{tabular}
&
\begin{tabular}{lll}
"Optical"  &   CD, DVD            &   1958 \\
"Magnetic" &   Tape, Drives       &   1956 \\
"Electric" &   NAND Flash         &   1980
\end{tabular}
& 
\begin{tabular}{lll}
"Gen. X" & New Math         &   1970               \\
"Gen. Y" & Math Wars        &   1990              \\
"Gen. Z" & MOOC             &   2010 
\end{tabular} \\ \hline
\end{tabular} \\

Of course, these snapshots are massive simplifications.
These tables hope to illustrate the fast-paced time we live in. \\

{\bf Remarks and Acknowledgements}. The April 7 2013, version of this document has appeared in \cite{CFZ}.
We added new part F)-N) in the last section as well as the last 16 graphics examples as well as
a couple of more code examples and references. We would like to thank
{\bf Enrique Canessa}, {\bf Carlo Fonda} and {\bf Marco Zennaro}
for organizing the workshop at ICTP and inviting us, {\bf Marius Kintel}
for valuable information on OpenSCAD, {\bf Daniel Pietrosemoli} at MediaLab Prado for introducing
us to use ``Processing" for 3D scanning, {\bf Gaia Fior} for printing out demonstration models and showing how 
easy 3D printing can appear, {\bf Ivan Bortolin} for 3D printing an appollonian cone for us, 
{\bf Stefan Rossegger} (CERN) for information on file
format conversions, {\bf Gregor L\"utolf} (http://www.3drucken.ch) 
for information on relief conversion and being a true pioneer for 3D printing in the classroom. 
Thanks to {\bf Thomas Tucker} for references and pictures of the genus 2 group. 
Thanks to {\bf Anna Zevelyov}, Director of Business Development at 
the Artec group company, for providing us with an educational evaluation licence of Artec Studio 9.1
which allowed us to experiment with 3D scanning. 

\pagebreak


\end{document}